\documentclass[11pt]{article}

\usepackage{amssymb, amsmath, amsthm, latexsym}
\usepackage{fullpage, color}
\usepackage{url}
\usepackage{algorithm}
\usepackage{algpseudocode}
\usepackage{tabularx}
\usepackage{paralist}
\usepackage{mathtools}
\usepackage{tcolorbox}
\usepackage{xcolor}

\usepackage{bbm} 

\usepackage{makecell}
\usepackage{multirow}
\usepackage{booktabs}

\usepackage[flushleft]{threeparttable} 

\usepackage{caption}
\usepackage{multirow}
\usepackage{colortbl}
\definecolor{bgcolor}{rgb}{0.8,1,1}
\definecolor{bgcolor2}{rgb}{0.8,1,0.8}
\definecolor{bgcolor3}{rgb}{0.67, 0.94, 0.82}

\usepackage{nicefrac}       

\definecolor{niceblue}{rgb}{0.0,0.19,0.56}

\usepackage{hyperref}
\hypersetup{colorlinks,linkcolor={blue},citecolor={niceblue},urlcolor={blue}}

\usepackage[numbers]{natbib}
\usepackage{sidecap}

\usepackage{pifont}
\usepackage{subfigure}

\newcommand{\R}{\mathbb{R}}

\def\<#1,#2>{\left\langle #1,#2\right\rangle}

\usepackage{mdframed} 
\usepackage{thmtools}

\definecolor{shadecolor}{gray}{0.9}
\declaretheoremstyle[
headfont=\normalfont\bfseries,
notefont=\mdseries, notebraces={(}{)},
bodyfont=\normalfont,
postheadspace=0.5em,
spaceabove=1pt,
mdframed={
  skipabove=8pt,
  skipbelow=8pt,
  hidealllines=true,
  backgroundcolor={shadecolor},
  innerleftmargin=4pt,
  innerrightmargin=4pt}
]{shaded}

\declaretheorem[style=shaded,within=section]{definition}
\declaretheorem[style=shaded,sibling=definition]{theorem}

\declaretheorem[style=shaded,sibling=definition]{corollary}

\declaretheorem[style=shaded,sibling=definition]{lemma}
\declaretheorem[style=shaded,sibling=definition]{remark}

\usepackage[colorinlistoftodos,bordercolor=orange,backgroundcolor=orange!20,linecolor=orange,textsize=scriptsize]{todonotes}

\newcommand{\cO}{{\cal O}}

\newcommand{\mA}{{\bf A}}

\newcommand{\mC}{{\bf C}}

\newcommand{\mI}{{\bf I}}

\newcommand{\mM}{{\bf M}}

\newcommand{\mT}{{\bf T}}
\newcommand{\mU}{{\bf U}}

\newcommand{\mLambda}{{\bf \Lambda}}

\def\R{\mathbb{R}}

\def\R{\mathbb R}

\def\la{\langle}
\def\ra{\rangle}

\usepackage{xspace}

\newcommand{\algname}[1]{{\sf #1}\xspace}

\usepackage{hyperref}

\usepackage{accents}
\newlength{\dhatheight}

\begin{document}
\title{Averaged Heavy-Ball Method}
\author{Marina Danilova$^{1,2}$ \quad Grigory Malinovsky$^3$
}
\date{$^1$ Institute of Control Sciences of Russian Academy of Sciences, Russia\\
$^2$ Moscow Institute of Physics and Technology, Russia\\
$^3$ King Abdullah University of Science and Technology, Saudi Arabia}
\maketitle
\begin{abstract}
Heavy-Ball method (\algname{HB}) is known for its simplicity in implementation and practical efficiency. However, as with other momentum methods, it has non-monotone behavior, and for optimal parameters, the method suffers from the so-called \textit{peak effect}. To address this issue, in this paper, we consider an averaged version of Heavy-Ball method (\algname{AHB}). We show that for quadratic problems \algname{AHB} has a smaller maximal deviation from the solution than \algname{HB}. Moreover, for general convex and strongly convex functions, we prove non-accelerated rates of global convergence of \algname{AHB} and its weighted version. We conduct several numerical experiments on minimizing quadratic and non-quadratic functions to demonstrate the advantages of using averaging for \algname{HB}.
\end{abstract}

\tableofcontents

\section{Introduction}\label{sec:intro}

First-order optimization methods have good convergence guarantees and are simple to implement. Therefore they are widely used in various applications. In particular, accelerated or first-order momentum methods such as Nesterov's method \cite{nesterov1983method} and Heavy-Ball method \cite{polyak1964some} and their various extensions are prevalent in some practically essential tasks, e.g., training of deep neural networks.

Due to its efficiency in solving non-convex optimization problems \cite{danilova2020recent}, Heavy-Ball method gained significant attention in recent years. As a result, a number of its modifications were proposed, including stochastic \cite{yang2016unified, taylor2019stochastic, defazio2020understanding}, zeroth-order \cite{gorbunov2020stochastic}, and distributed variants \cite{yu2019linear, mishchenko2019distributed}, to mention a few.

However, even for simple (strongly) convex problems, accelerated/momentum methods have non-monotone behavior. For example, in the recent paper \cite{danilova2018non}, the authors show that Heavy-Ball method (\algname{HB}) with optimal parameters has so-called \textit{peak-effect} even for simple quadratic minimization problems. This means that in this case the distance to the solution during the initial iterations of \algname{HB}. Moreover, the maximal distance is proportional to $\sqrt{\varkappa}$ \cite{danilova2018non, mohammadi2021transient}, where $\varkappa$ is the condition number of the problem. Therefore, for ill-conditioned problems ($\varkappa \gg 1$) peak-effect can be significant.

\paragraph{Contributions.} To address this issue, in this work, we consider an averaged version of the Heavy-Ball method called Averaged Heavy-Ball method (\algname{AHB}). We study the maximal deviation of this method for quadratic functions and prove the global convergence guarantees in the convex and strongly convex (not necessarily quadratic) cases for \algname{AHB} and its version based on the weighted averaging (\algname{WAHB}). For quadratic functions with a specific property of the spectrum, our theoretical results show that there exists a choice of parameters for \algname{AHB} such that momentum parameter $\beta$ is sufficiently large but the maximal deviation is significantly smaller than for \algname{HB} with optimal parameters. We derive global complexity results for \algname{AHB} and \algname{WAHB} matching the best-known ones for \algname{HB}. To the best of our knowledge, we prove the first global convergence results for \algname{HB} with averaging in the strongly convex case (see the summary in Table~\ref{tab:known_results}). Moreover, our numerical experiments corroborate our theoretical observations and show that \algname{HB} with a properly adjusted averaging scheme converges faster than \algname{HB} without averaging and has smaller oscillations.

\subsection{Preliminaries}
We focus on the following minimization problem
\begin{equation}
    \min\limits_{x\in \R^n}f(x),\label{eq:main_problem}
\end{equation}
where $f: \R^n \to \R$ is $L$-smooth and $\mu$-strongly convex function.

\begin{definition}[$L$-smoothness]\label{def:L_smooth}
    Differentiable function $f: \R^n \to \R$ is called $L$-smooth for some constant $L > 0$, if its gradient is $L$-Lipschitz, i.e., for all $x,y\in \R^n$
    \begin{equation}
        \|\nabla f(x) - \nabla f(y)\|_2 \le L\|x - y\|_2. \label{eq:L_smooth_def}
    \end{equation}
\end{definition}

\begin{definition}[$\mu$-strong convexity]\label{def:mu_str_cvx}
    Differentiable function $f: \R^n \to \R$ is called $\mu$-strongly convex for some constant $\mu \ge 0$, if for all $x,y\in \R^n$ the following inequality holds:
    \begin{equation}
        f(y) \ge f(x) + \langle \nabla f(x), y-x \rangle + \frac{\mu}{2}\|y-x\|_2^2. \label{eq:mu_str_cvx_def}
    \end{equation}
\end{definition}

Throughout the paper we use standard notation for optimization literature \cite{polyak1987introduction, nesterov2018lectures}, e.g., $x^*$ denotes the solution of \eqref{eq:main_problem}, $R_0 = \|x_0 - x^*\|_2$ is the distance from the starting point to the solution, $\varkappa = \nicefrac{L}{\mu}$ is the condition number of the problem.

\subsection{Related work}

\begin{algorithm}[h]
\caption{Heavy-Ball method (\algname{HB})}
\label{alg:HB_m}   
\begin{algorithmic}[1]
\Require starting points $x_0$, $x_1$ (by default $x_0 = x_1$), number of iterations $N$, stepsize $\alpha > 0$, momentum parameter $\beta \in [0,1]$
\For{$k=0,\ldots, N-1$}
\State  $x_{k+1} = x_k - \alpha \nabla f(x_k) + \beta (x_k - x_{k-1})$
\EndFor
\Ensure $x_k$ 
\end{algorithmic}
\end{algorithm}

\paragraph{Convergence guarantees for Heavy-Ball method.} Heavy-Ball method \cite{polyak1964some} (\algname{HB}, Algorithm~\ref{alg:HB_m}) is the first optimization method with momentum proposed in the literature. In \cite{polyak1964some}, the author proves \textit{local} $\cO(\sqrt{\nicefrac{L}{\mu}}\log(\nicefrac{1}{\varepsilon}) )$ convergence rate for twice continuously
differentiable $L$-smooth and $\mu$-strongly convex functions. The first global convergence results for \algname{HB} are obtained in \cite{ghadimi2015global}, where the authors derive \textit{global} $\cO\left(\nicefrac{LR_0^2}{\varepsilon}\right)$ convergence rate of \algname{HB} and \algname{AHB} for $L$-smooth convex ($\mu = 0$) functions and $\cO(\nicefrac{L}{\mu}\log(\nicefrac{1}{\varepsilon}))$ convergence rate of \algname{HB} for $L$-smooth and $\mu$-strongly convex functions. Although these results establish the global convergence of \algname{HB} (and \algname{AHB} in the convex case), the rates are non-accelerated, i.e., they are not optimal \cite{nemirovsky1983problem} unlike the local convergence rate derived in \cite{polyak1964some}. This issue is partially resolved in \cite{lessard2016analysis}, where the authors prove that \algname{HB} converges with the asymptotically accelerated rate for strongly convex quadratic functions. Moreover, they also show that there exists a non-twice differentiable strongly convex function such that \algname{HB} does not converge for this objective. Next, using Performance Estimation Problem tools \cite{taylor2017performance, taylor2018lyapunov, taylor2019stochastic}, one can show that for standard choices of parameters \algname{HB} has the non-accelerated rate of convergence. However, the following question remains open: \textit{does there exist a choice of parameters for \algname{HB} such that the method converges globally with the accelerated rate for twice differentiable $L$-smooth and (strongly) convex functions?} Although we do not address this question in our work, we highlight it here due to its theoretical importance.

\paragraph{Non-monotone behavior of Heavy-Ball method.} From the classical analysis of \algname{HB} \cite{polyak1964some}, it is known that the following choice of parameters $\alpha$ and $\beta$ ensures the best convergence rate for \algname{HB} up to the numerical constant factors:
\begin{equation}
\alpha = \alpha^* = \frac{4}{(\sqrt{L} + \sqrt{\mu})^2}, \quad 
\beta = \beta^* = \left(\frac{\sqrt{L} - \sqrt{\mu}}{\sqrt{L} + \sqrt{\mu}}\right)^2. \label{eq:optimal_params}
\end{equation}
However, recently it was shown \cite{danilova2018non} that \algname{HB} with optimal parameters suffers from the so-called \textit{peak effect} at the beginning of the convergence. In particular, the maximal deviation can be of the order $\sqrt{\varkappa} = \sqrt{\nicefrac{L}{\mu}}$. Similar results were also derived in \cite{mohammadi2021transient}. However, in practice, it is worth mentioning that the optimal parameters from \eqref{eq:optimal_params} are rarely used and, as a result, the non-monotonicity of \algname{HB} is not that significant.

\section{Maximal Deviations on Quadratic Problems}\label{sec:quadratic_case}
In this section, we consider the instance of \eqref{eq:main_problem} with $f(x)$ being a quadratic function. That is, we assume that $f(x) = \frac12 x^{\top} \mA x$, where $\mA \in \mathbb{S}^n_{++}$ is a $n \times n$ positive definite matrix. For this problem, we prove that Averaged Heavy-Ball method with a certain choice of parameters has a smaller deviation of the iterates from the optimum at initial iterations than the Heavy-Ball method with optimal parameters.

\subsection{Heavy-Ball Method}
Recently it was shown \cite{danilova2018non} that \algname{HB} with optimal parameters \eqref{eq:optimal_params} suffers from so-called \textit{peak effect} at the beginning of the convergence. In particular, according to the following theorem, the maximal deviation can be of the order $\sqrt{\varkappa}$.
\begin{theorem}[Theorem 1 from \cite{danilova2018non}]\label{thm:HB_bad_example}
    Consider $f(x) = \frac12 x^{\top} \mA x, \; \mA = \mathrm{diag}\left(\mu, \ \lambda_2, \ \ldots, \ \lambda_{n-1}, L\right)$, where $\mu \le \lambda_2 \le \lambda_3 \le \ldots \le \lambda_{n-1} \le L$ . Then for $x^0=x^1 = (1,1,\ldots,1)^{\top}$ the iterates $\{x_k\}_{k \ge 0}$ produced by \algname{HB} with $\alpha = \alpha^*,$ $\beta = \beta^*$ satisfy
    \begin{equation}
        \max\limits_{k} \|x_k\|_{\infty} \ge \frac{\sqrt{\varkappa}}{2e}.\label{eq:HB_peak}
    \end{equation}
\end{theorem}

\begin{algorithm}[t]
\caption{Averaged Heavy-Ball method (\algname{AHB})}
\label{alg:AHB_m}   
\begin{algorithmic}[1]
\Require  starting points $x_0$, $x_1$ (by default $x_0 = x_1$), number of iterations $N$, stepsize $\alpha > 0$, momentum parameter $\beta \in [0,1]$
\For{$k=1,\ldots, N-1$}
\State  $x_{k+1} = x_k - \alpha \nabla f(x_k) + \beta (x_k - x_{k-1})$
\State  $\overline{x}_{k+1} = \frac{1}{k+2}\sum\limits_{i=0}^{k+1} x_i$ \Comment{One can recurrently implement this step: $\overline{x}_{k+1} = \frac{k\overline{x}_{k} + x_{k+1}}{k+1}$}
\EndFor
\Ensure $\overline{x}_k$ 
\end{algorithmic}
\end{algorithm}

\subsection{Averaged Heavy-Ball method}\label{sec:heavy_ball_method_average}
In this subsection, we consider the modification of \algname{HB} that returns the average of the iterates produced by \algname{HB}. We call the resulting method Averaged Heavy-Ball method (\algname{AHB}, see Algorithm~\ref{alg:AHB_m}).

We start with showing that for the same initialization, \algname{AHB} with $\alpha = \nicefrac{1}{L}$ and not too large $\beta$ has significantly more minor deviations than \algname{HB} with optimal parameters when $\varkappa$ is sufficiently large under some assumptions on the spectrum of $\mA$. 
\begin{theorem}\label{thm:AHB_bad_example_deviations}
    Consider $f(x) = \frac12 x^{\top} \mA x$ with $ \mA = \mathrm{diag}\left(\mu, \ \lambda_2, \ \ldots, \ \lambda_{n-1}, L\right)$, where $\mu \le \lambda_2 \le \lambda_3 \le \ldots \le \lambda_{n-1} \le L$ and $\lambda_2 \ge 10\mu$, $L \ge 100\mu$. Then for $x^0=x^1 = (1,1,\ldots,1)^{\top}$ and for all $k \ge 0$ the iterates $\{\overline{x}_k\}_{k \ge 0}$ generated by \algname{AHB} with $\alpha = \nicefrac{1}{L},$ $\beta \in [(1-3\sqrt{\nicefrac{\mu}{L}})^2,  (1-2\sqrt{\nicefrac{\mu}{L}})^2]$ satisfy 
    \begin{equation}
        \max\limits_{k} \|\overline{x}_k\|_{\infty} \le 2.\label{eq:AHB_peak}
    \end{equation}
\end{theorem}

That is, comparing bounds \eqref{eq:HB_peak} and \eqref{eq:AHB_peak} for $\varkappa \gg 1$, we conclude \algname{AHB} with the parameters from Theorem~\ref{thm:AHB_bad_example_deviations} has much smaller deviations then \algname{HB} with parameters from \eqref{eq:optimal_params}. However, Theorem~\ref{thm:AHB_bad_example_deviations} works only for the particular initialization. The guarantees independent of $x^0, x^1$ are much more valuable and that is what we derive in the next subsection.

\subsection{Maximal Deviation of \algname{AHB} for Arbitrary Initialization}\label{sec:max_dev_arbitrary_init}
Consider the matrix representation of \algname{HB} update rule:
\begin{equation}
    \begin{bmatrix}
        x_{k+1} - x_*\\
        x_k - x_*
    \end{bmatrix} = \mT \cdot \begin{bmatrix}
        x_{k} - x_*\\
        x_{k-1} - x_*
    \end{bmatrix} = \ldots = \mT^k \cdot \begin{bmatrix}
        x_{1} - x_*\\
        x_{0} - x_*
    \end{bmatrix} , \label{eq:HB_matrix_update}
\end{equation}
where
\begin{equation}
    \mT = \left[\begin{array}{@{}c|c@{}}
     \begin{matrix}
     (1 + \beta) \mI - \alpha \mA
     \end{matrix}
     & -\beta \mI \\
    \hline
     \mI &
     {\bf 0}
    \end{array}\right] \in \R^{2n\times 2n},\quad \begin{bmatrix}
        x_{k+1} - x_*\\
        x_k - x_*
    \end{bmatrix} \in \R^{2n}. \label{eq:matrix_T_definition}
\end{equation}
Therefore, we have
\begin{equation}
    x_{k} - x^* = \underbrace{\begin{bmatrix}
         {\bf 0} & \mI
    \end{bmatrix}}_{\mC} \mT^k \begin{bmatrix}
        x_{1} - x_*\\
        x_{0} - x_*
    \end{bmatrix}.\label{eq:HB_new_matrix_represntation}
\end{equation}
For convenience, we also introduce the following notation:
\begin{equation*}
    z_k = \begin{bmatrix}
         x_{k+1} - x_*\\
         x_k - x_*
    \end{bmatrix}.
\end{equation*}
Following \cite{mohammadi2021transient}, we study the worst case deviation $\|x_k - x_*\|_2$ in the relation to $\|z_0\|_2$, i.e., we focus on the following quantity
\begin{equation*}
    \max\limits_{k\ge 0}\sup\limits_{z_0 \neq 0} \frac{\|x_k - x_*\|_2}{\|z_0\|_2} \overset{\eqref{eq:HB_new_matrix_represntation}}{=} \max\limits_{k\ge 0}\sup\limits_{z_0 \neq 0} \frac{\|\mC\mT^kz_0\|_2}{\|z_0\|_2} = \max\limits_{k\ge 0}\|\mC\mT^k\|_2 
\end{equation*}
that is the largest spectral norm of the matrices $\mC \mT^k$ for $k \ge 0$. Clearly, one can choose $z_0$, i.e., starting points $x_0$ and $x_1$, in such a way that $z_0$ is in the direction of the principal right singular vector of $\mC \mT^k$ implying $\|x_k - x_*\|_2 = \|\mC \mT^k\|_2\|z_0\|_2$. Therefore, $\|\mC \mT^k\|_2$ is a tight and natural measure of the worst case deviation of the iterates produced by \algname{HB}. Since this quantity depends on the choice of $\alpha$ and $\beta$ we denote it as $\text{dev}_{\algname{HB}}(\alpha,\beta) := \max\limits_{k\ge 0}\|\mC\mT^k(\alpha,\beta)\|_2$.

For \algname{AHB} we know
\begin{equation}
    \overline{x}_{k} - x_* = \frac{1}{k+1}\sum\limits_{t=0}^k(x_k - x_*) = \frac{1}{k+1}\sum\limits_{t=0}^{k}\mC\mT^t \begin{bmatrix}
        x_{1} - x_*\\
        x_{0} - x_*
    \end{bmatrix}.\notag
\end{equation}
We introduce new notation:
\begin{equation*}
    \text{dev}_{\algname{AHB}}(\alpha,\beta) := \max\limits_{k\ge 0}\left\|\frac{1}{k+1}\sum\limits_{t=0}^k\mC\mT^t(\alpha,\beta)\right\|_2.
\end{equation*}
As for \algname{HB}, $ \text{dev}_{\algname{AHB}}(\alpha,\beta)$ is also a reasonable measure of the worst case deviation of the iterates produced by \algname{AHB}. Moreover, due to the Jensen's inequality and convexity of $\|\cdot\|_2$ we have $\text{dev}_{\algname{AHB}}(\alpha,\beta) \le \text{dev}_{\algname{HB}}(\alpha,\beta)$.

\begin{theorem}\label{thm:AHB_deviation_arbitrary_init}
    Consider $f(x) = \frac12 x^{\top} \mA x$ with $\mA = \mA^\top \succ 0$ with eigenvalues $\lambda_1 \le \ldots \le \lambda_n$, $\lambda_2 \ge F^2\lambda_1$, $F > 14$, $F \le \sqrt{\nicefrac{\lambda_n}{\lambda_1}}$, $\lambda_n \ge 10000\lambda_1$. Then the maximal deviation of \algname{AHB} and \algname{HB} with $\alpha = \nicefrac{1}{L}$ and $(1 - \sqrt{\nicefrac{\lambda_2}{\lambda_n}})^2< \beta \le (1 - F\sqrt{\nicefrac{\lambda_1}{\lambda_n}})^2$ is at least $\nicefrac{(\sqrt{F^2-1})}{2e\sqrt{6}}$ times smaller then the maximal deviation of \algname{HB} with $\alpha = \alpha^*$ and $\beta = \beta^*$ given in \eqref{eq:optimal_params}:
    \begin{equation}
        \text{dev}_{\algname{AHB}}(\alpha,\beta) \le \text{dev}_{\algname{HB}}(\alpha,\beta) \le \frac{2e\sqrt{6}}{\sqrt{F^2-1}}\text{dev}_{\algname{HB}}(\alpha^*,\beta^*). \label{eq:HB_deviation_init_independent}
    \end{equation}
\end{theorem}

The constant $\nicefrac{2e\sqrt{6}}{\sqrt{F^2-1}}$ can be sufficiently small and $\beta$ can be sufficiently large at the same time when the condition number $\varkappa$ is large enough. For example, for $\varkappa = 10^8$ and $F = 200$ one can choose $\beta = \left(1 - \nicefrac{F}{\sqrt{\varkappa}}\right)^2 \approx 0.96$ and get $\nicefrac{2e\sqrt{6}}{\sqrt{F^2-1}} \approx 0.067$.

\section{Convergence Guarantees for Non-Quadratics}\label{sec:non_quadratics}

In this section, we study the convergence of \algname{AHB} for problems \eqref{eq:main_problem} with (strongly) convex and smooth objectives. First global convergence guarantees for \algname{HB} and \algname{AHB} in the convex case were obtained in \cite{ghadimi2015global}. In the same paper, the authors derived the convergence rate for \algname{HB} in the strongly convex case. See the summary of known results in Table~\ref{tab:known_results}. 

In contrast, for \algname{HB} with averaging, there are no convergence results in the strongly convex case. Below we consider two options to derive such results.

\begin{table*}[h]
    \centering
    \small
	\caption{\small Summary of known and new results on the maximal deviation and complexity bounds for \algname{HB} and its variants with averaging. Column ``Max.\ deviation'' contains the results on the maximal deviation of the methods on quadratic minimization problems (see the details in Section~\ref{sec:max_dev_arbitrary_init}), columns ``Complexity, $\mu = 0$'' and ``Complexity, $\mu > 0$'' show iteration complexity bounds for the methods applied to \eqref{eq:main_problem} with $f$ being $L$-smooth and convex / $\mu$-strongly convex but not necessarily quadratic, i.e., number of iterations needed to guarantee that the output of the method $\widehat{x}$ satisfies $f(\widehat x) - f(x_*) \le \varepsilon$ where $x_*$ is the solution of \eqref{eq:main_problem}. Our results are highlighted in green. Notation: $\varkappa = \nicefrac{L}{\mu}$ (condition number), $\Delta_0 = f(x_0) - f(x_*)$, $R_0 = \|x_0 - x_*\|_2$.}
    \label{tab:known_results}    
   \begin{threeparttable}
    \begin{tabular}{|c|c|c c c|}
         \hline
         Method & Citation & Max.\ deviation & Complexity, $\mu = 0$ & Complexity, $\mu > 0$ \\ 
\hline\hline
    \algname{HB} & \cite{danilova2018non,ghadimi2015global} & $\frac{\sqrt{\varkappa}}{2e}$\tnote{\color{red} (1)} & $\frac{LR_0^2}{\varepsilon}$\tnote{\color{red} (2)} & $\frac{\varkappa}{1-\beta}\log\frac{\Delta_0}{\varepsilon}$\tnote{\color{red} (3)} \\
    \algname{AHB} & \cite{ghadimi2015global} & N/A & $\frac{LR_0^2}{\varepsilon} + \frac{\beta LR_0^2}{(1-\beta)\varepsilon}$ & N/A \\
    \rowcolor{bgcolor3}\algname{AHB} & Thm.~\ref{thm:AHB_deviation_arbitrary_init} \& \ref{thm:WAHB_main_result} \& \ref{thm:R-AHB} & $\frac{\sqrt{6\varkappa}}{\sqrt{F^2-1}}$\tnote{\color{red}(4)} & $\frac{LR_0^2}{\varepsilon} + \frac{LR_0^2\sqrt{\beta}}{(1-\beta)\varepsilon}$ & $\left(\varkappa + \frac{\varkappa\sqrt{\beta}}{1-\beta}\right)\log\frac{\mu R_0^2}{\varepsilon}$\tnote{\color{red}(5)} \\
    \rowcolor{bgcolor3}\algname{WAHB} & Thm.~\ref{thm:WAHB_main_result} & $\frac{\sqrt{6\varkappa}}{\sqrt{F^2-1}}$\tnote{\color{red}(6)} & $\frac{LR_0^2}{\varepsilon} + \frac{LR_0^2\sqrt{\beta}}{(1-\beta)\varepsilon}$ & $\left(\varkappa + \frac{\varkappa\sqrt{\beta}}{1-\beta}\right)\log\frac{LR_0^2\left(1 + \frac{\sqrt{\beta}}{(1-\beta)}\right)}{\varepsilon}$ \\
    \hline
    \end{tabular}
    \begin{tablenotes}
      {\scriptsize
        \item [{\color{red}(1)}] This result is obtained for \algname{HB} with optimal parameters from \eqref{eq:optimal_params} (see Theorem~\ref{thm:HB_bad_example}).
        \item [{\color{red}(2)}] The complexity bound is obtained for iteration-dependent parameters: $\beta_k = \frac{k}{k+2}$, $\alpha_k = \frac{1}{L(k+1)}$.
        \item [{\color{red}(3)}] This result holds for $\alpha \in (0,\nicefrac{1}{L})$, $\beta \in [0, \sqrt{(1-\alpha L)(1-\alpha\mu)}]$. When $\varkappa \gg 1$ this assumption implies that $\beta \le 0.75$. In practical applications, e.g., training deep neural networks, much larger values for parameter $\beta$ are usually used.
        \item [{\color{red}(4)}] The result holds for a special class of quadratic functions described in Theorem~\ref{thm:AHB_deviation_arbitrary_init}. Parameters $\alpha$ and $\beta$ for \algname{AHB} are given there as well. Here $F$ is such that $\lambda_2 \ge F^2\mu$, $F> 14$, $F \le \sqrt{\varkappa}$, where $\lambda_2$ is the second smallest eigenvalue of the Hessian matrix. For large enough $\varkappa$ and $F$ one can guarantee that maximal deviation for \algname{AHB} with parameters from Theorem~\ref{thm:AHB_deviation_arbitrary_init} is much smaller than for \algname{HB} with optimal parameters from \eqref{eq:optimal_params}.
        \item [{\color{red}(5)}] The complexity bound is proven Restarted version of \algname{AHB} (\algname{R-AHB}, Algorithm~\ref{alg:R-AHB}).
        \item [{\color{red}(6)}] See {\color{red}(4)} and Remark~\ref{rem:WAHB_max_dev}.
      }
    \end{tablenotes}
    \end{threeparttable}
\end{table*}

\subsection{Weighted Averaged Heavy-Ball Method}
One way to obtain them is to change the averaging weights, see Weighted Averaged Heavy-Ball method (\algname{WAHB}, Algorithm~\ref{alg:WAHB}).
\begin{algorithm}[h]
\caption{Weighted Averaged Heavy-Ball method (\algname{WAHB})}
\label{alg:WAHB}   
\begin{algorithmic}[2]
\Require number of iterations $N$, stepsize $\alpha > 0$, momentum parameter $\beta \in [0,1]$, starting points $x_0$, $x_1$ (by default $x_1 = x_0 - \alpha \nabla f(x_0)$), weights for the averaging $\{w_k\}_{k=0}^N > 0$
\For{$k=1,\ldots, N-1$}
\State  $x_{k+1} = x_k - \alpha \nabla f(x_k) + \beta (x_k - x_{k-1})$
\State  $\overline{x}_{k+1} = \frac{1}{W_{k+1}}\sum\limits_{i=0}^{k+1} w_ix_i$, where $W_{k+1} = \sum\limits_{i=0}^{k+1}w_i$ \Comment{Recurrent analog: $\overline{x}_{k+1} = \frac{W_k\overline{x}_{k} + w_{k+1}x_{k+1}}{W_{k+1}}$}
\EndFor
\Ensure $\overline{x}_N$ 
\end{algorithmic}
\end{algorithm}
When $w_k = 1$ for all $k\ge 0$ \algname{WAHB} recovers \algname{AHB}. However, it is natural to choose larger $w_k$ for larger $k$: for such a choice of $w_k$ the method gradually ``forgets'' about the early iterates that should lead to faster convergence. Guided by this intuition we provide a rigorous analysis of \algname{WAHB} with gradually increasing $w_k$.

\begin{remark}\label{rem:WAHB_max_dev}
    We emphasize that the proof of Theorem~\ref{thm:AHB_deviation_arbitrary_init} holds for non-uniform averaging as well. That is, under assumptions of Theorem~\ref{thm:AHB_deviation_arbitrary_init} we have
    \begin{equation*}
        \text{dev}_{\algname{WAHB}}(\alpha,\beta) \le \text{dev}_{\algname{HB}}(\alpha,\beta) \le \frac{2e\sqrt{6}}{\sqrt{F^2-1}}\text{dev}_{\algname{HB}}(\alpha^*,\beta^*),
    \end{equation*}
    where 
    \begin{equation*}
        \text{dev}_{\algname{WAHB}}(\alpha,\beta) := \max\limits_{k\ge 0}\left\|\frac{1}{W_k}\sum\limits_{t=0}^kw_t\mC\mT^t(\alpha,\beta)\right\|_2.
    \end{equation*}
\end{remark}

In our derivations, we rely on the following representation of the update rule of \algname{HB} with $x_1 = x_0 - \alpha\nabla f(x_0)$:
\begin{equation}
    x_{k+1} = x_k - m_k, \quad m_{k} = \beta m_{k-1} + \alpha \nabla f(x_k),\quad m_{-1} = 0. \label{eq:HB_equiv} 
\end{equation}
Indeed, since $m_{k-1} = x_{k-1} - x_{k}$ for all $k \ge 0$ (for convenience, we use the notation $x_{-1} = x_0$) we have
\begin{equation*}
    x_{k+1} = x_k - m_k = x_k -\alpha \nabla f(x_k) -\beta m_{k-1} = x_k - \alpha \nabla f(x_k) + \beta(x_k - x_{k-1}).
\end{equation*}
Next, following \cite{mania2017perturbed,yang2016unified} we consider \textit{perturbed} or \textit{virtual} iterates:
\begin{equation}
    \widetilde{x}_k = x_k - \frac{\beta}{1-\beta}m_{k-1},\; k\ge 0. \label{eq:virtual_iterates_HB}
\end{equation}
We notice, that these iterates are not computed explicitly in the method. However, they turn out to be useful in the analysis because of the following relation: for all $k \ge 0$
\begin{eqnarray}
    \widetilde{x}_{k+1} &=& x_{k+1} - \frac{\beta}{1-\beta}m_{k} = x_k - \frac{1}{1-\beta}m_k = \widetilde{x}_k + \frac{\beta}{1-\beta}m_{k-1} -\frac{1}{1-\beta}\left(\beta m_{k-1} + \alpha \nabla f(x_k)\right) \notag\\
    &=& \widetilde{x}_k - \frac{\alpha}{1-\beta}\nabla f(x_k). \label{eq:virtual_iter_recurrence}
\end{eqnarray}

Using this notation, we we derive the following lemma measuring one iteration progress of \algname{HB}.

\begin{lemma}\label{lem:one_iter_progress_HB}
    Assume that $f$ is $L$-smooth and $\mu$-strongly convex. Let $\alpha$ and $\beta$ satisfy
    \begin{equation}
        0 < \alpha \le \frac{1-\beta}{4L}, \quad \beta \in [0,1). \label{eq:params_WAHB_1}
    \end{equation}
    Then, for all $k\ge 0$
    \begin{equation}
        \frac{\alpha}{2(1-\beta)}\left(f(x_k) - f(x_*)\right) \le \left(1 - \frac{\alpha\mu}{2(1-\beta)}\right)\|\widetilde{x}_k - x_*\|_2^2 - \|\widetilde{x}_{k+1} - x_*\|_2^2 + \frac{3L\alpha\beta^2}{(1-\beta)^3}\|m_{k-1}\|_2^2. \label{eq:one_iter_progress_HB}
    \end{equation}
\end{lemma}

As the next step, it is natural to sum up inequalities \eqref{eq:one_iter_progress_HB} for $k = 0,1,2,\ldots, K$ with weights $w_k = \left(1 - \frac{\alpha\mu}{2(1-\beta)}\right)^{-(k+1)}$ to get the bound on $f(\overline{x}_K) - f(x_*)$. However, in this case, we obtain
\begin{equation*}
    \frac{3L\alpha\beta^2}{(1-\beta)^3}\sum\limits_{k=0}^K w_k\|m_{k-1}\|^2
\end{equation*}
in the upper bound for $f(\overline{x}_K) - f(x_*)$. Therefore, we need to estimate this sum and this is exactly what the next lemma is about.

\begin{lemma}\label{lem:weighted_sum_of_momentums}
    Assume that $f$ is $L$-smooth and $\mu$-strongly convex. Let $\alpha$ and $\beta$ satisfy
    \begin{equation}
        0 < \alpha \le \frac{(1-\beta)^2}{4L\sqrt{3\beta}}, \quad \beta \in [0,1). \label{eq:params_WAHB_2}
    \end{equation}
    Then, for all $k\ge 0$
    \begin{equation}
        \frac{3L\alpha\beta^2}{(1-\beta)^3}\sum\limits_{k=0}^K w_k\|m_{k-1}\|^2 \le \frac{\alpha}{4(1-\beta)}\sum\limits_{k=0}^{K}w_k\left(f(x_k) - f(x_*)\right), \label{eq:weighted_sum_of_momentums}
    \end{equation}
    where $w_k = \left(1 - \frac{\alpha\mu}{2(1-\beta)}\right)^{-(k+1)}$.
\end{lemma}

Combining Lemmas~\ref{lem:one_iter_progress_HB}~and~\ref{lem:weighted_sum_of_momentums} we obtain the following result.
\begin{theorem}\label{thm:WAHB_main_result}
    Assume that $f$ is $L$-smooth and $\mu$-strongly convex. Let $\alpha$ and $\beta$ satisfy
    \begin{equation}
        0 < \alpha \le \min\left\{\frac{1-\beta}{4L}, \frac{(1-\beta)^2}{4L\sqrt{3\beta}} \right\}, \quad \beta \in [0,1). \label{eq:params_WAHB_3}
    \end{equation}
    Then, after $K \ge 0$ iterations of \algname{WAHB} we have
    \begin{equation}
        f(\overline{x}_K) - f(x_*) \le \frac{4(1-\beta)\|x_0 - x_*\|_2^2}{\alpha W_K}, \label{eq:WAHB_main_result}
    \end{equation}
    where $w_k = \left(1 - \frac{\alpha\mu}{2(1-\beta)}\right)^{-(k+1)}$. That is, if $\mu > 0$, then
    \begin{equation}
        f(\overline{x}_K) - f(x_*) \le \left(1 - \frac{\alpha\mu}{2(1-\beta)}\right)^K\frac{4(1-\beta)\|x_0 - x_*\|_2^2}{\alpha}, \label{eq:WAHB_main_result_str_cvx}
    \end{equation}
    and if $\mu = 0$, we have
    \begin{equation}
        f(\overline{x}_K) - f(x_*) \le \frac{4(1-\beta)\|x_0 - x_*\|_2^2}{\alpha K}. \label{eq:WAHB_main_result_cvx}
    \end{equation}
\end{theorem}

The following complexity results trivially follow from this theorem.
\begin{corollary}\label{cor:WAHB_complexity}
    Let the assumptions of Theorem~\ref{thm:WAHB_main_result} hold and
    \begin{equation*}
        \alpha = \min\left\{\frac{1-\beta}{4L}, \frac{(1-\beta)^2}{4L\sqrt{3\beta}} \right\}.
    \end{equation*}
    Then, to achieve $f(\overline{x}_K) - f(x_*) \le \varepsilon$ for $\varepsilon > 0$ \algname{WAHB} requires
    \begin{equation}
        \cO\left(\left(\frac{L}{\mu} + \frac{L\sqrt{\beta}}{\mu(1-\beta)}\right)\log\frac{LR_0^2\left(1 + \nicefrac{\sqrt{\beta}}{(1-\beta)}\right)}{\varepsilon}\right) \label{eq:WAHB_compl_str_cvx}
    \end{equation}
    iterations when $\mu > 0$, and
    \begin{equation}
        \cO\left(\frac{LR_0^2}{\varepsilon} + \frac{LR_0^2\sqrt{\beta}}{(1-\beta)\varepsilon}\right) \label{eq:WAHB_compl_cvx}
    \end{equation}
    iterations when $\mu = 0$, where $R_0 \ge \|x_0 - x_*\|_2$.
\end{corollary}

When $\mu = 0$ \algname{WAHB} recovers \algname{AHB} since $w_k = 1$ by definition. Therefore, in the convex case, this result establishes the complexity of \algname{AHB}.

\subsection{Restarted Averaged Heavy-Ball Method}
An alternative way to achieve linear convergence in the strongly convex case for Heavy-Ball method with averaging is to use the restarts technique. That is, consider Restarted Averaged Heavy-Ball method (\algname{R-AHB}, Algorithm~\ref{alg:R-AHB}).
\begin{algorithm}[h]
\caption{Restarted Averaged Heavy-Ball method (\algname{R-AHB})}
\label{alg:R-AHB}   
\begin{algorithmic}[1]
\Require number of restarts $\tau$, numbers of iterations $\{N_t\}_{t=1}^\tau$, stepsizes $\{\alpha_t\}_{t=1}^\tau > 0$, momentum parameters $\{\beta_t\}_{t=1}^{\tau} \in [0,1]$, starting point $x_0$
\State $\widehat x_0 = x_0$
\For{$t=1,\ldots, \tau$}
\State  Run \algname{AHB} (Algorithm~\ref{alg:AHB_m}) for $N_t$ iterations with stepsize $\alpha_t$, momentum parameter $\beta_t$, and starting points $\widehat x_{t-1}$, $\widehat x_{t-1} - \alpha_t \nabla f(\widehat x_{t-1})$. Define the output of \algname{AHB} by $\widehat x_t$.
\EndFor
\Ensure $\widehat{x}_\tau$ 
\end{algorithmic}
\end{algorithm}
The work of the method is split into stages. Each stage is the run of \algname{AHB} from the point obtained at the previous stage, the first stage initializes at the given point.

Based on the convergence result for \algname{AHB} in the convex case, one can get the convergence rate of \algname{R-AHB} in the strongly convex case.
\begin{theorem}\label{thm:R-AHB}
    Assume that $f$ is $L$-smooth and $\mu$-strongly convex. Let $\alpha_t = \alpha$, $\beta_t = \beta$, $N_t = N$ for all $t=1,\ldots,\tau$ and
    \begin{equation}
        0 < \alpha \le \min\left\{\frac{1-\beta}{4L}, \frac{(1-\beta)^2}{4L\sqrt{3\beta}} \right\}, \quad \beta \in [0,1),\quad N = \left\lceil \frac{16(1-\beta)}{\alpha\mu} \right\rceil. \label{eq:params_R-AHB}
    \end{equation}
    Then, after $\tau = \max\{\left\lceil\log_2(\nicefrac{\mu R_0^2}{\varepsilon})\right\rceil - 1,1\}$ iterations with $R_0 \ge \|x_0 - x_*\|_2$ \algname{R-AHB} produces such point $\widehat{x}_\tau$ that $f(\widehat{x}_\tau) - f(x_*) \le \varepsilon$. Furthermore, if
    \begin{equation*}
        \alpha = \min\left\{\frac{1-\beta}{4L}, \frac{(1-\beta)^2}{4L\sqrt{3\beta}} \right\},
    \end{equation*}
    then the total number of \algname{AHB} iterations equals
    \begin{equation}
        \cO\left(\left(\frac{L}{\mu} + \frac{L\sqrt{\beta}}{\mu(1-\beta)}\right)\log\frac{\mu R_0^2}{\varepsilon}\right). \label{eq:R-AHB_compl}
    \end{equation}
\end{theorem}

\section{Numerical Experiments}\label{sec:numerical_exp}

We conducted several numerical experiments to compare the behavior of \algname{HB} with and without averaging applied to minimize quadratic functions and solve logistic regression problem. The code was written in Python 3.7 using standard libraries.

\subsection{Quadratic Functions}
In this section, we consider three quadratic functions:
\begin{eqnarray}
    f_{\text{random}}(x) &=& \frac{1}{2}x^\top \mA_{\text{rand}} x - (x^{*})^\top \mA_{\text{rand}} x, \label{eq:random_quadratic_function}\\
    f_{\text{Nesterov}}(x) &=& \frac{L-\mu}{8}\left(x_1^2 + \sum\limits_{i=1}^{n-1}(x_i - x_{i+1})^2 - 2x_1\right) + \frac{\mu}{2}\|x\|^2, \label{eq:nesterov_function}\\
    f_{\text{Toeplitz}}(x) &=& \frac{1}{2}x^\top \mA_{\text{Toeplitz}} x, \label{eq:toeplitz_matrix}
\end{eqnarray}
where matrix $\mA_{\text{rand}} = \hat{\mA}^\top \hat{\mA}$, the elements of matrix $\hat\mA \in \R^{n\times n}$ are independently sampled from the standard Gaussian distribution, and $\mA_{\text{Toeplitz}} \in \R^{n\times n}$ is a Toeplitz with a first row $(2,-1,1,0,\ldots,0)$. Function from \eqref{eq:nesterov_function} is a classical function used to derive lower bounds for the complexity of first-order methods applied to minimize smooth strongly convex functions \cite{nesterov2018lectures}.

\begin{figure}[H]
    \centering
    \includegraphics[width=0.325\textwidth]{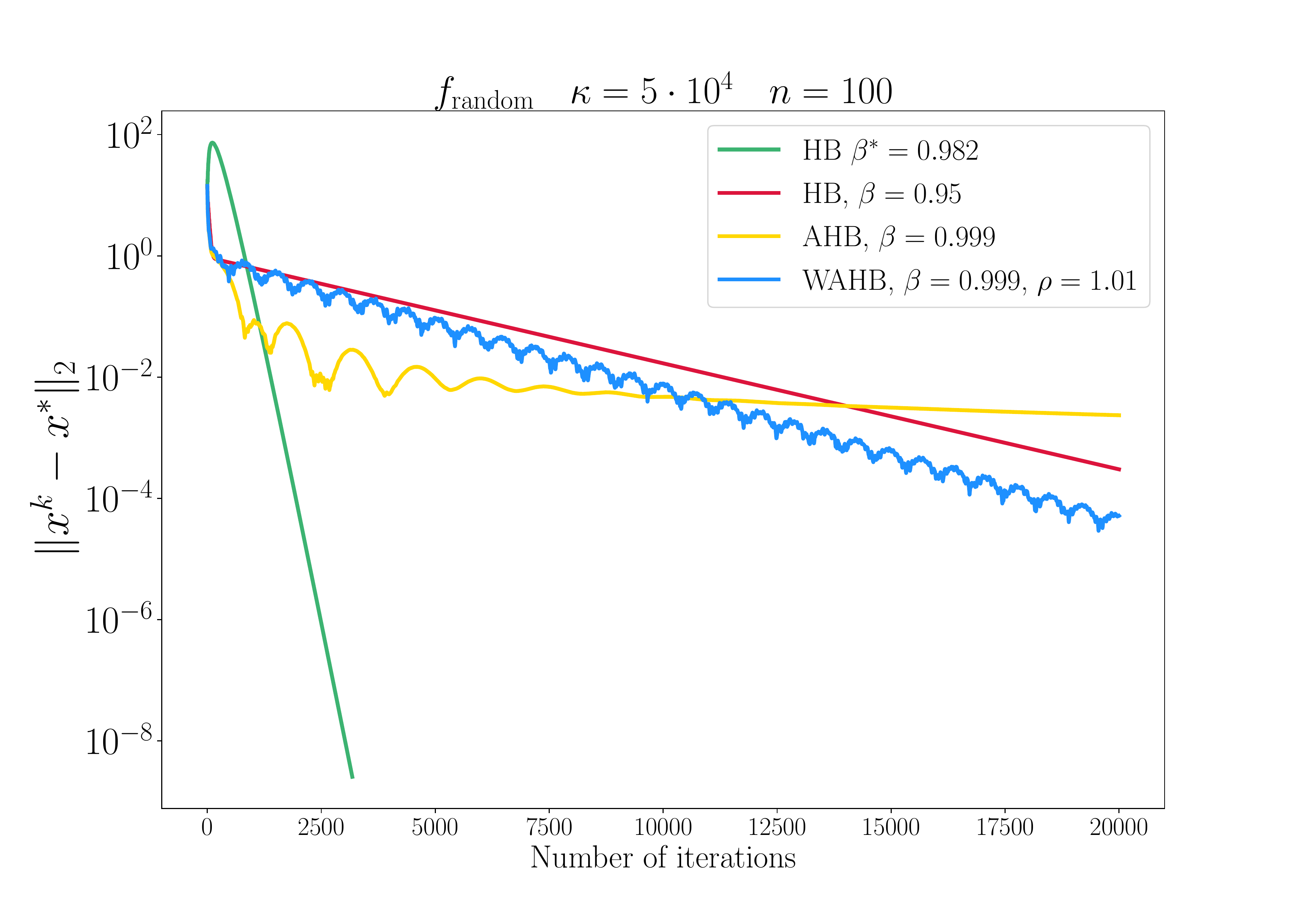}
    \includegraphics[width=0.325\textwidth]{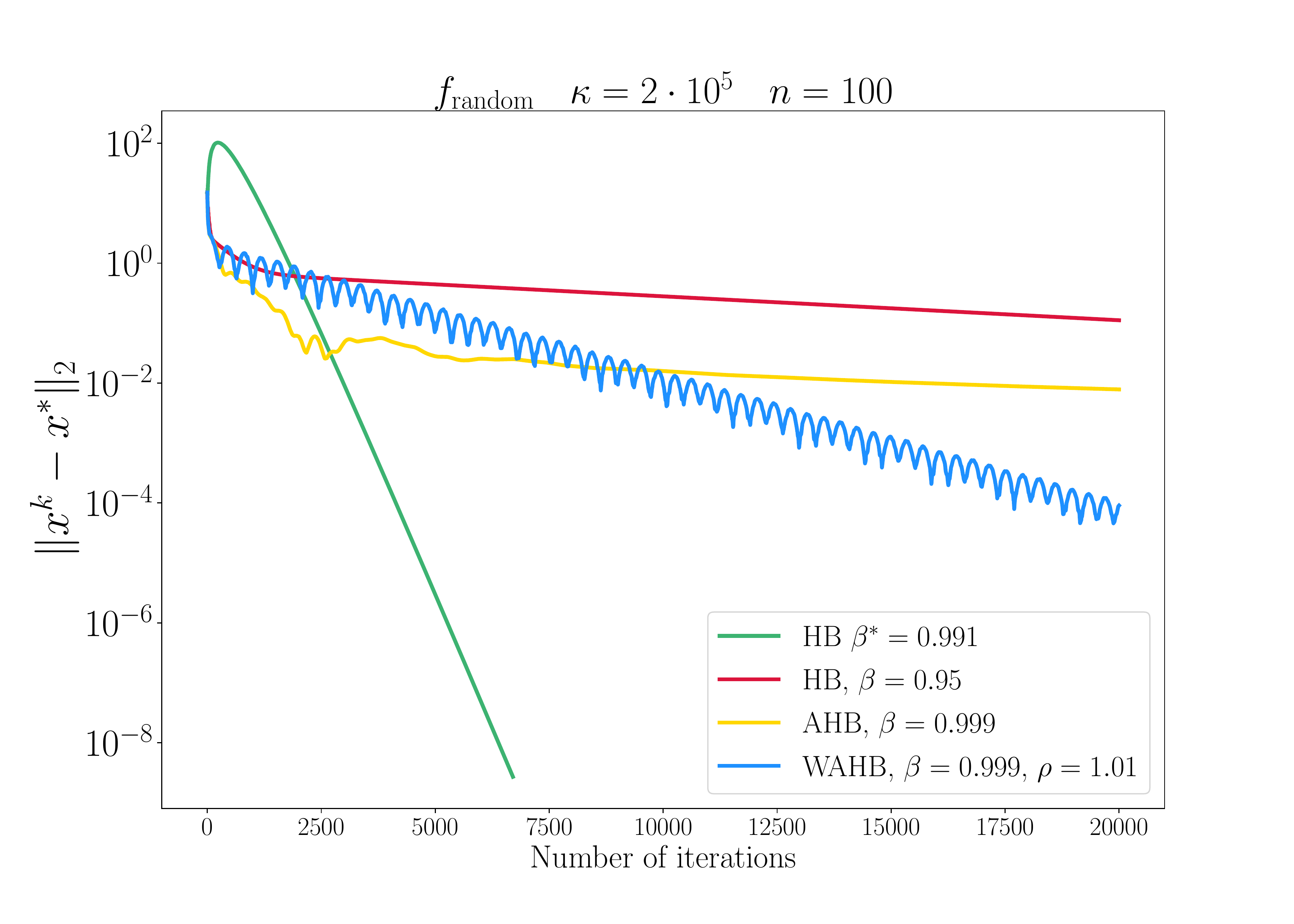}
    \includegraphics[width=0.325\textwidth]{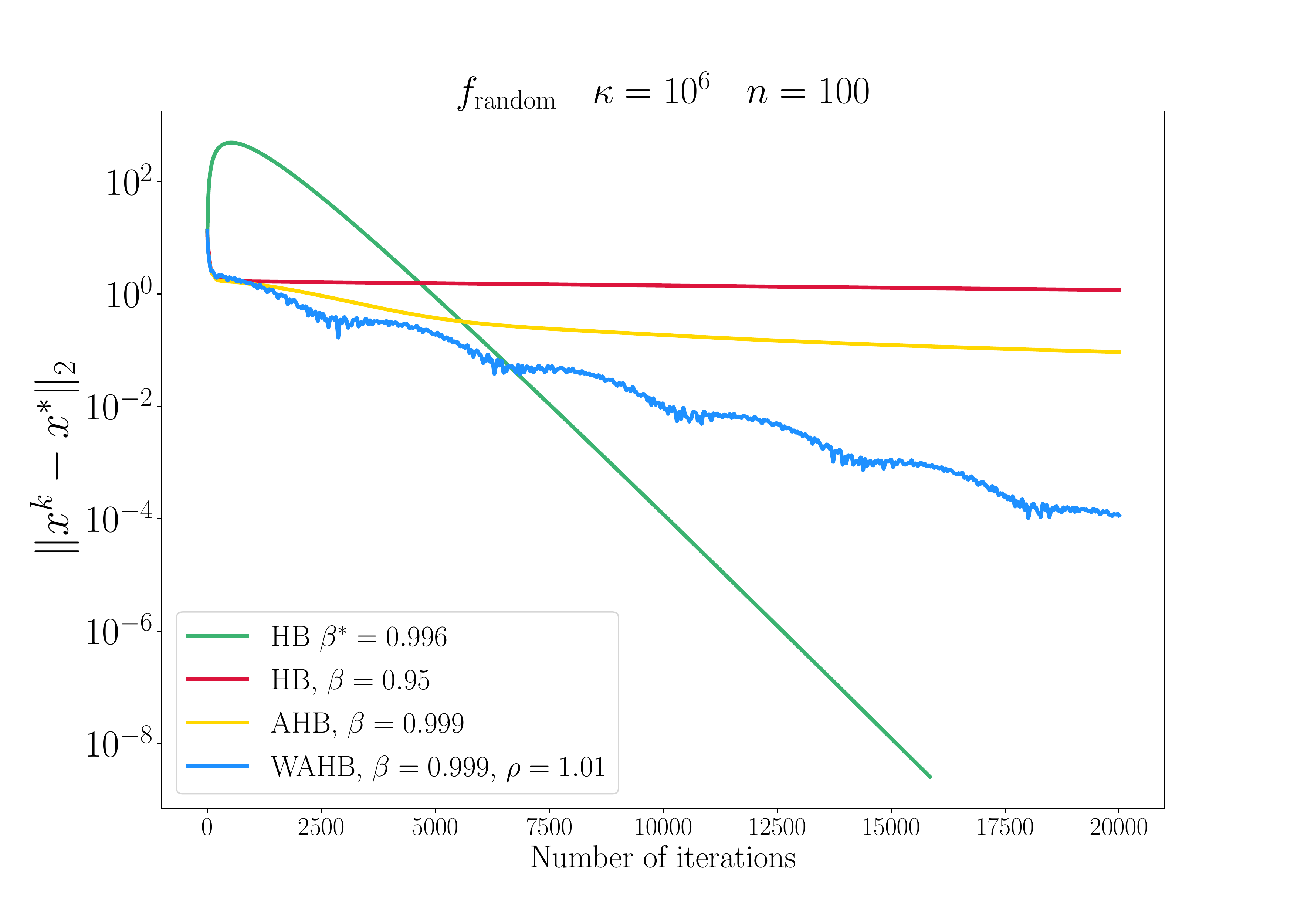}\\
    \includegraphics[width=0.325\textwidth]{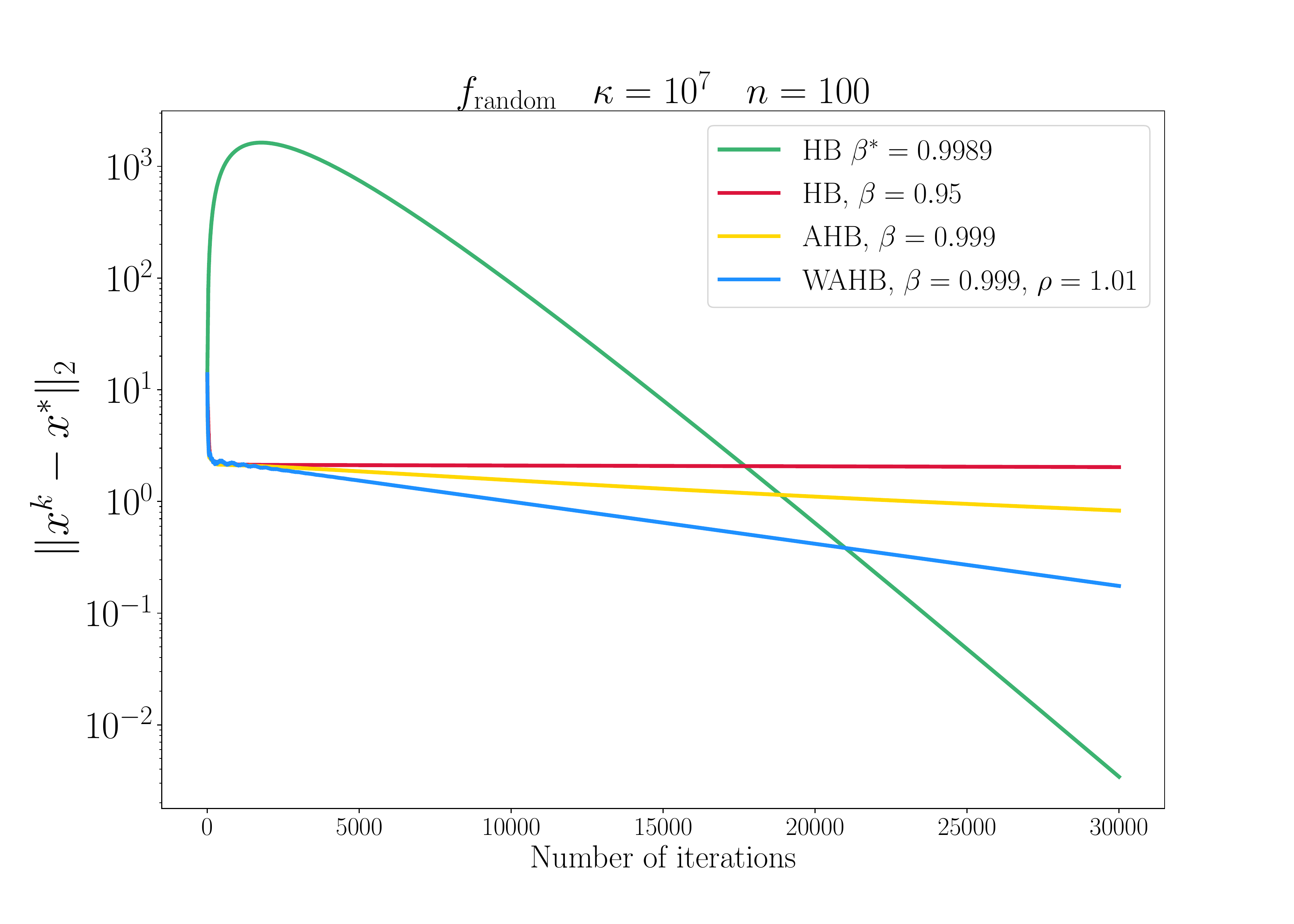}
    \includegraphics[width=0.325\textwidth]{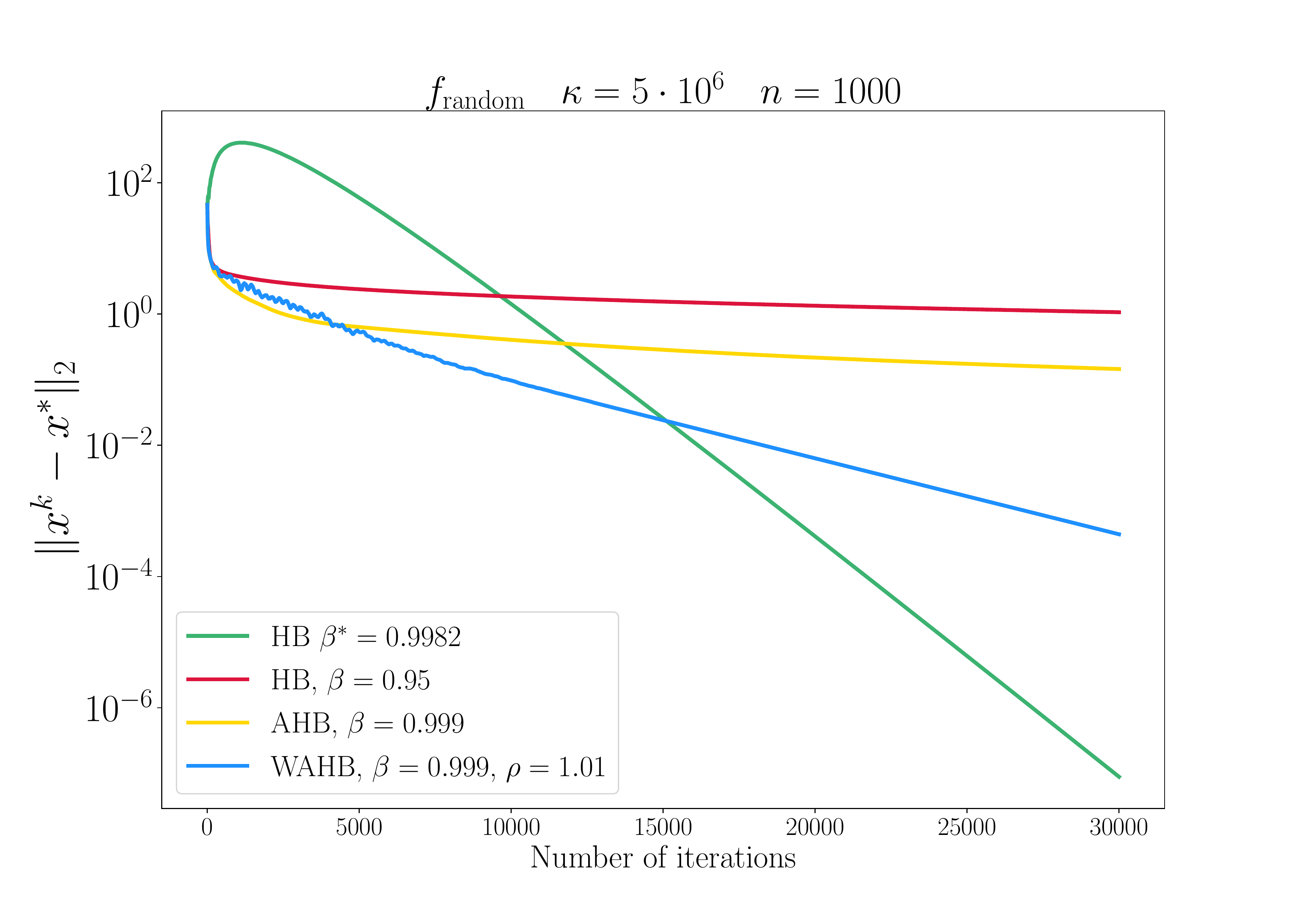}
    \includegraphics[width=0.325\textwidth]{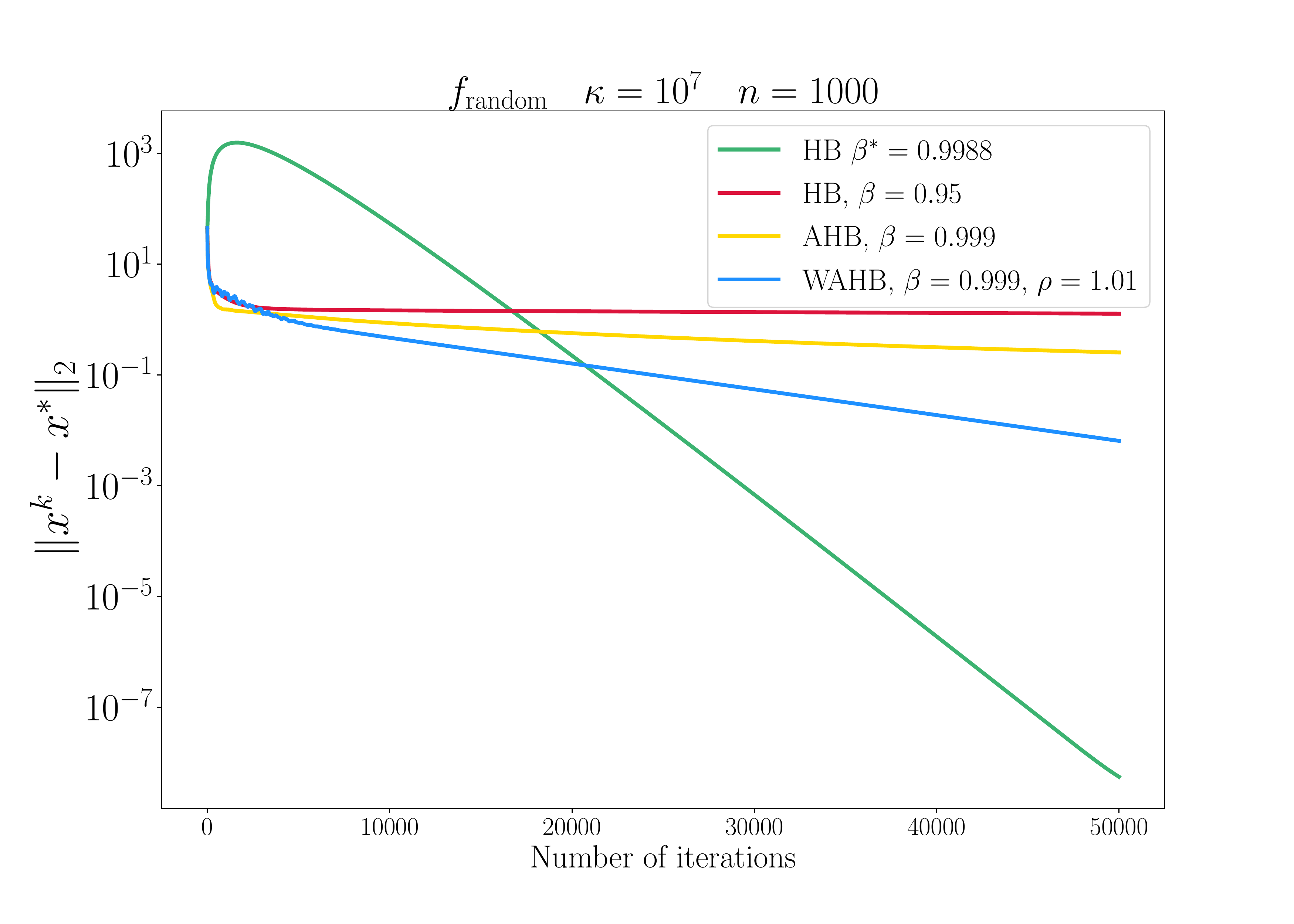}
    \caption{Trajectories of \algname{HB}, \algname{AHB}, and \algname{WAHB} applied to minimize a quadratic function from \eqref{eq:random_quadratic_function} with different condition numbers $\varkappa$ and dimension $n$.}
    \label{fig:random_quadratic}
\end{figure}

\begin{figure}[H]
    \centering
    \includegraphics[width=0.49\textwidth]{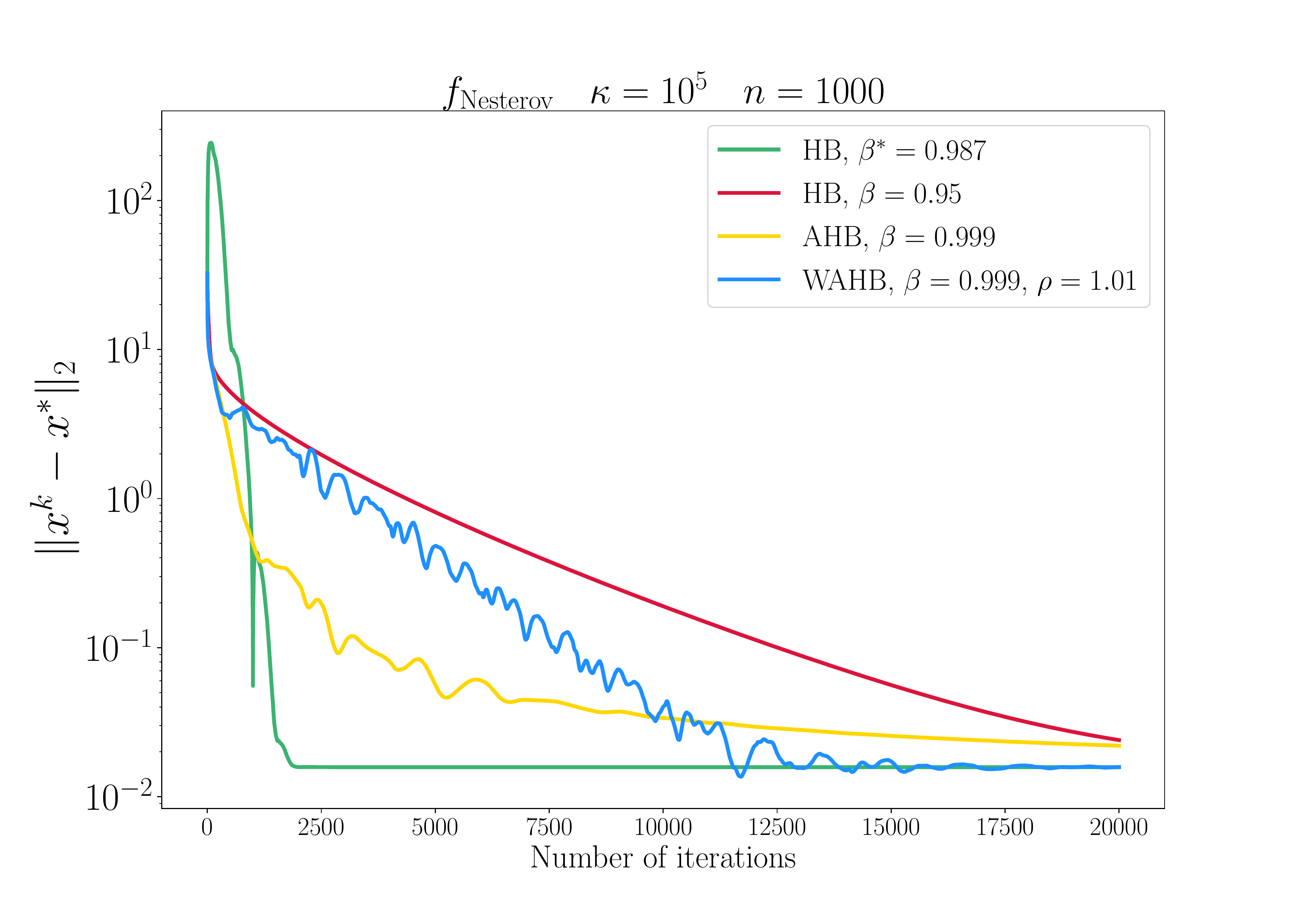}
    \includegraphics[width=0.49\textwidth]{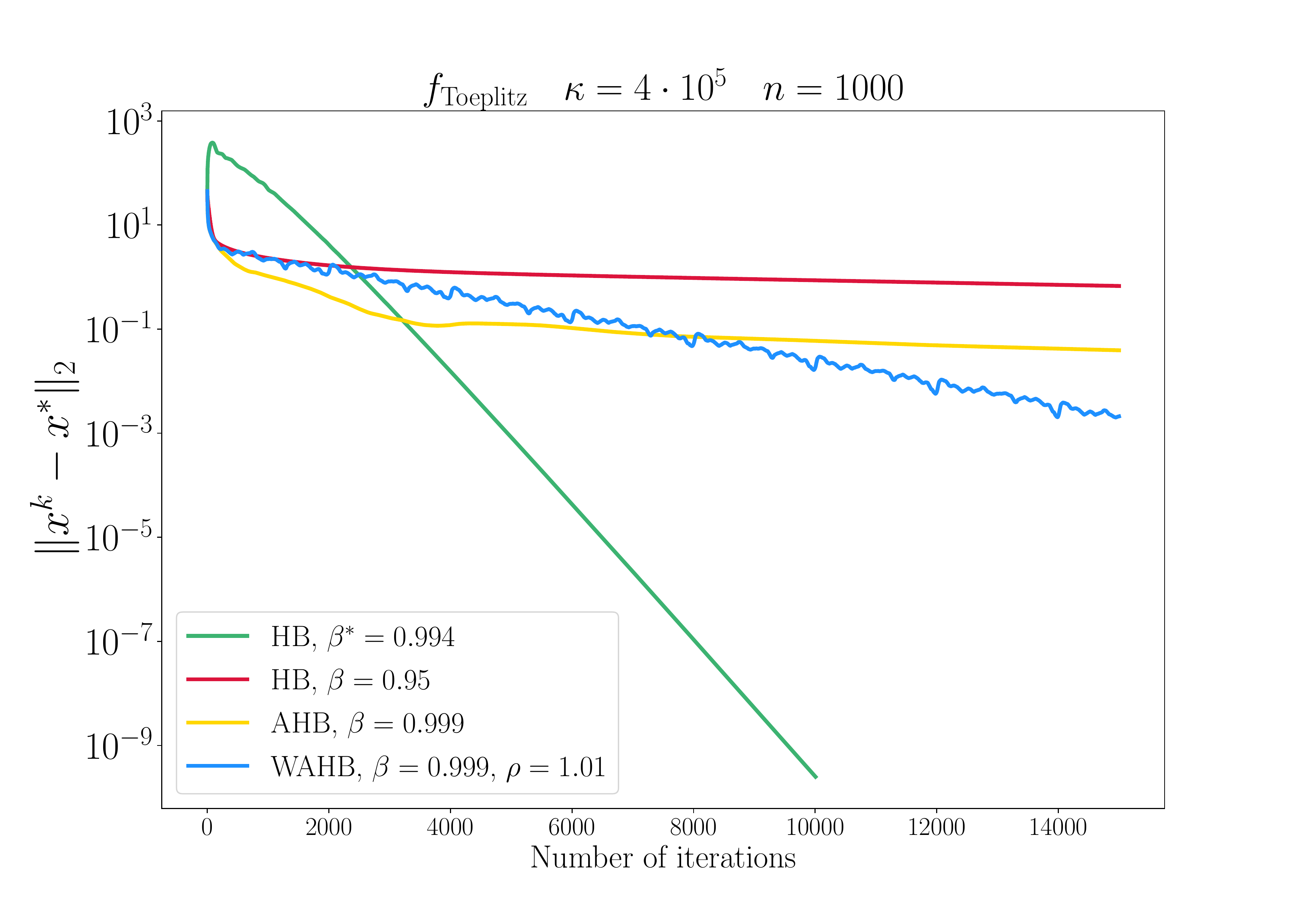}
    \caption{Trajectories of \algname{HB}, \algname{AHB}, and \algname{WAHB} applied to minimize a quadratic functions from \eqref{eq:nesterov_function} and \eqref{eq:toeplitz_matrix} with condition numbers $\varkappa\sim 10^5$ and dimension $n = 1000$.}
    \label{fig:nesterov_toeplitz}
\end{figure}

We run \algname{HB} with $\beta = 0.95$ (standard choice of $\beta$), \algname{AHB} and \algname{WAHB} with $\beta = 0.999$ (large $\beta$) to minimize each of these functions. For these methods we used stepsize $\alpha = \nicefrac{1}{L}$. The weights for \algname{WAHB} were chosen as $w_k = \rho^k$ for $\rho = 1.01$. Moreover, we also tested \algname{HB} with optimal parameters from \eqref{eq:optimal_params}. One can find the results in Figures~\ref{fig:random_quadratic}~and~\ref{fig:nesterov_toeplitz}.

These results show that methods with averaging (\algname{AHB} and \algname{WAHB}) converge reasonably well during the first iterations of the method even with large $\beta = 0.999$, which was larger than the optimal $\beta^*$ in all our experiments. Moreover, unlike \algname{HB} with optimal parameters, \algname{AHB} and \algname{WAHB} do not suffer from the peak effect. The absence of peak effect allows us to use \algname{HB} with averaging for the first iterates and then restart the method. Finally, we emphasize that \algname{HB} with $\beta = 0.95$ converges slower than \algname{WAHB} with $\beta = 0.999$ in all our experiments and slower than \algname{AHB} with $\beta = 0.999$ in almost all experiments (except the first one shown in Figure~\ref{fig:random_quadratic}). We also tested \algname{HB} with $\beta = 0.999$ and observed very slow convergence for the method in this case.

To conclude, our experiments on quadratic functions highlight the benefits of using \algname{AHB} and \algname{WAHB} with large $\beta$ and standard $\alpha = \nicefrac{1}{L}$.

\subsection{Logistic Regression with $\ell_2$-Regularization}
Next, we also consider logistic regression with $\ell_2$-regularization:
\begin{equation}
    \min\limits_{x\in \R^n}\left\{f(x) = \frac{1}{m}\sum\limits_{i=1}^m\log\left(1 + \exp\left(-y_i\cdot (\mA x)_i\right)\right) + \frac{\ell_2}{2}\|x\|_2^2\right\}, \label{eq:logreg_loss}
\end{equation}
where $m$ is the total number of data points/samples, $y_i \in \{-1,1\}$ is a label of $i$-th datapoint, and $\mA \in \R^{m\times d}$ is a feature matrix. This function is known to be $\ell_2$-strongly convex and $(L+\ell_2)$-smooth with $L = \nicefrac{\sigma_{\max}^2(\mA)}{4m}$, where $\sigma_{\max}(A)$ is the maximal singular value of matrix $\mA$. We take the datasets, i.e., pairs of $(\mA,\{y_i\}_{i=1}^m)$, from LIBSVM library \cite{chang2011libsvm}, see the summary of the considered datasets in Table~\ref{tab:summary_data}.

\begin{table}[h]
    \centering
    \caption{Summary of the considered datasets for the logistic regression.}\label{tab:summary_data}
    \begin{tabular}{|l|c|c|c|}
        \hline
        & {\tt a9a} & {\tt phishing} & {\tt w8a}\\
        \hline
        $m$ (\# of data points) & 32 561 & 11 055 & 49 749 \\
        \hline
        $d$ (\# of features) & 123 & 68 & 300\\
        \hline
    \end{tabular}
\end{table}

\begin{algorithm}[h]
\caption{Tail-Averaged Heavy-Ball method (\algname{TAHB})}
\label{alg:TAHB}   
\begin{algorithmic}[1]
\Require  starting points $x_0$, $x_1$ (by default $x_0 = x_1$), number of iterations $N$, stepsize $\alpha > 0$, momentum parameter $\beta \in [0,1]$, tail size $s \ge 0$
\For{$k=1,\ldots, N-1$}
\State  $x_{k+1} = x_k - \alpha \nabla f(x_k) + \beta (x_k - x_{k-1})$
\State  $\overline{x}_{k+1} = \begin{cases}\frac{1}{k+2}\sum\limits_{i=0}^{k+1} x_i,& \text{if } k+1 < s,\\ \frac{1}{s}\sum\limits_{i=0}^{s-1} x_{k+1-i},& \text{if } k+1 \ge s\end{cases}$ \Comment{It is required to store the last $s$ iterates}
\EndFor
\Ensure $\overline{x}_k$ 
\end{algorithmic}
\end{algorithm}

\begin{figure}[H]
    \centering
    \includegraphics[width=0.325\textwidth]{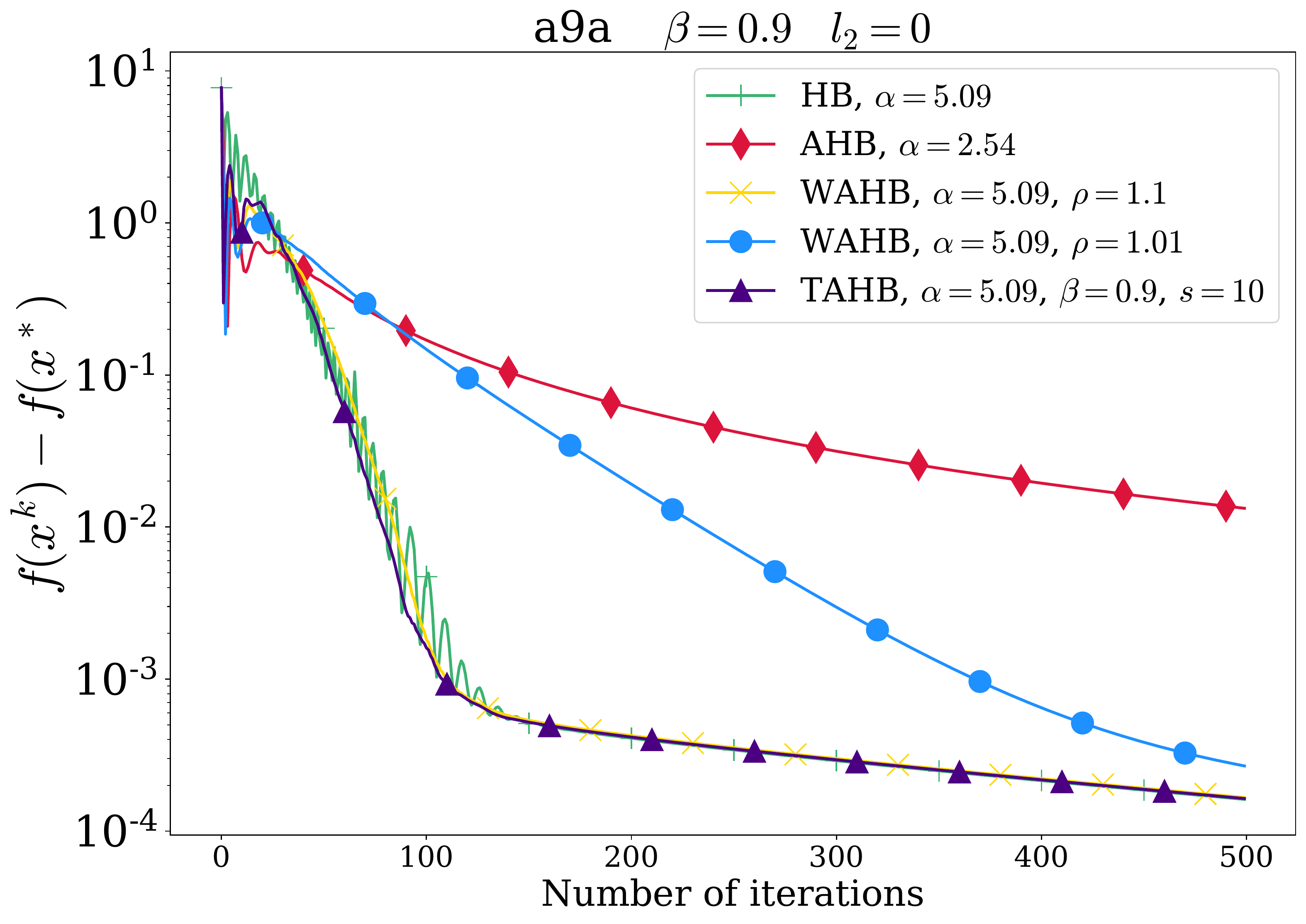}
    \includegraphics[width=0.325\textwidth]{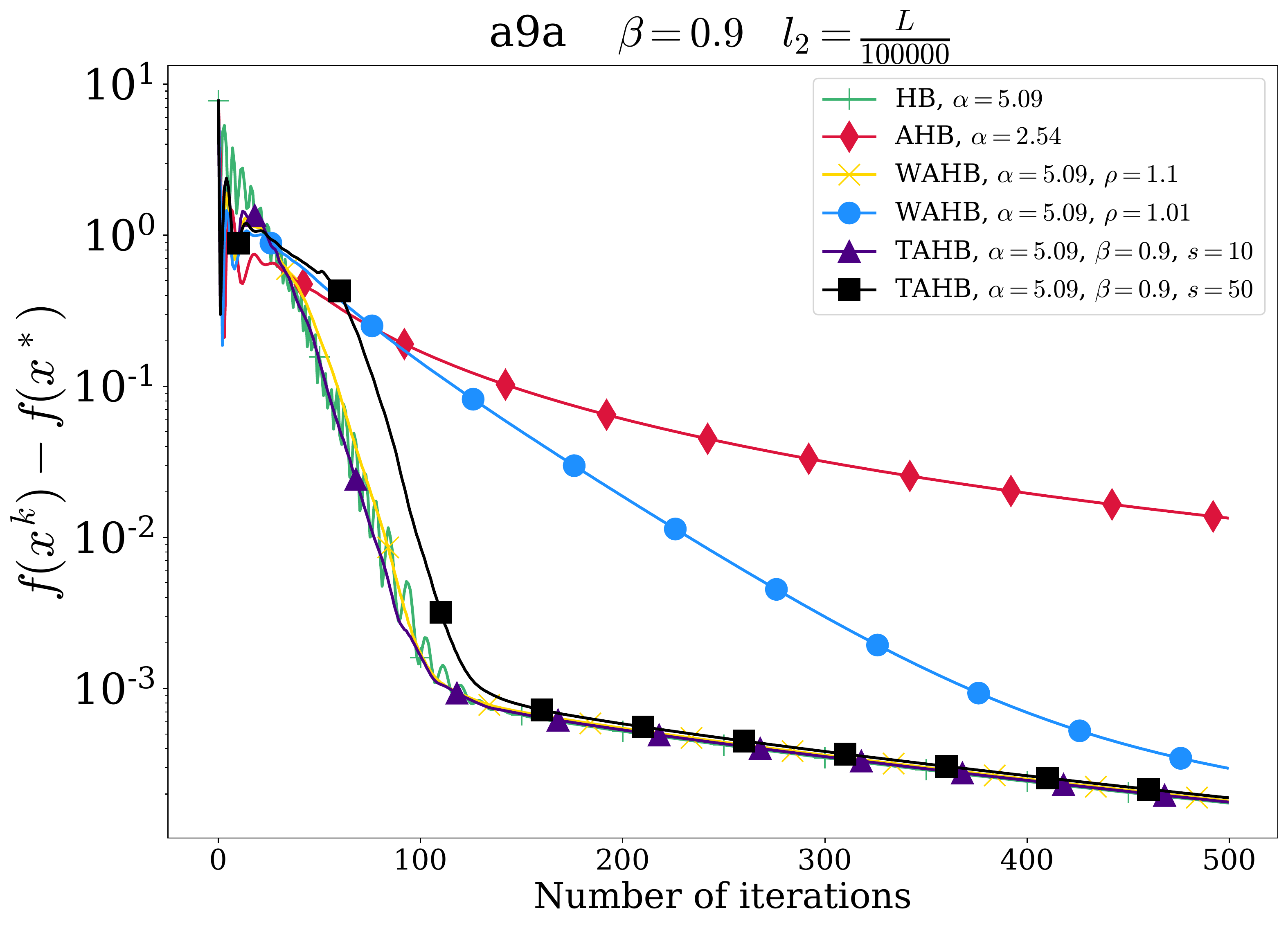}
    \includegraphics[width=0.325\textwidth]{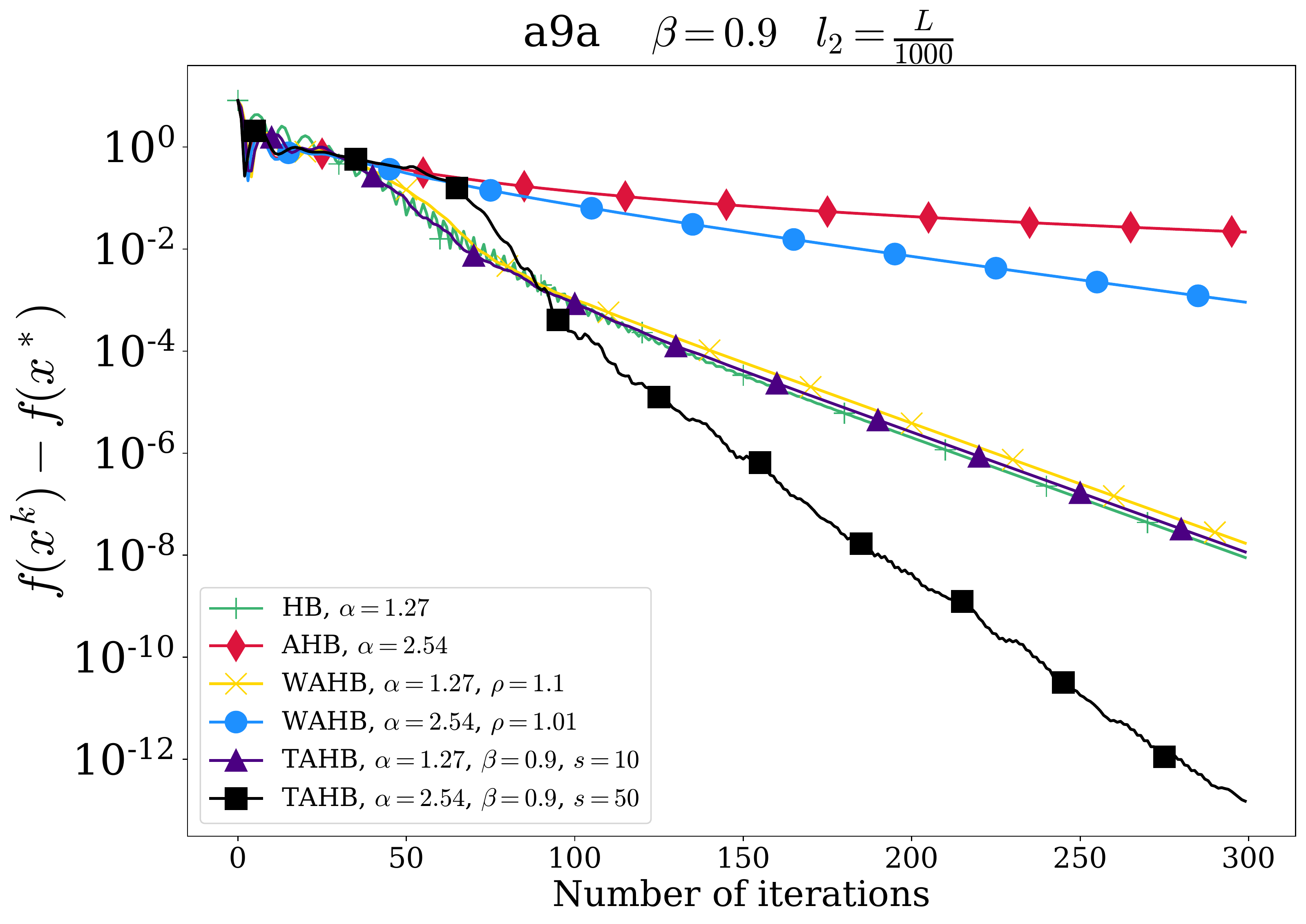}\\
    \includegraphics[width=0.325\textwidth]{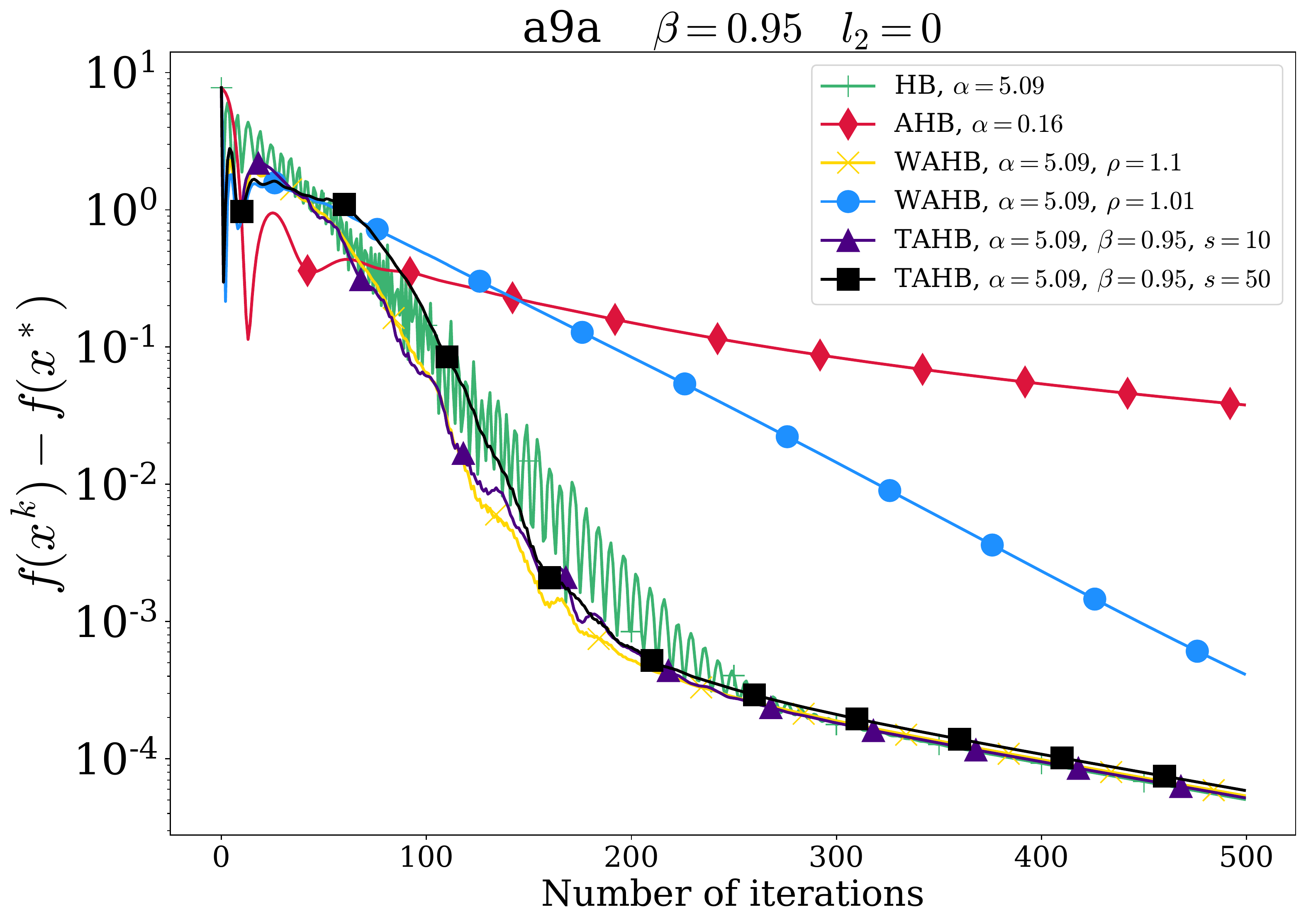}
    \includegraphics[width=0.325\textwidth]{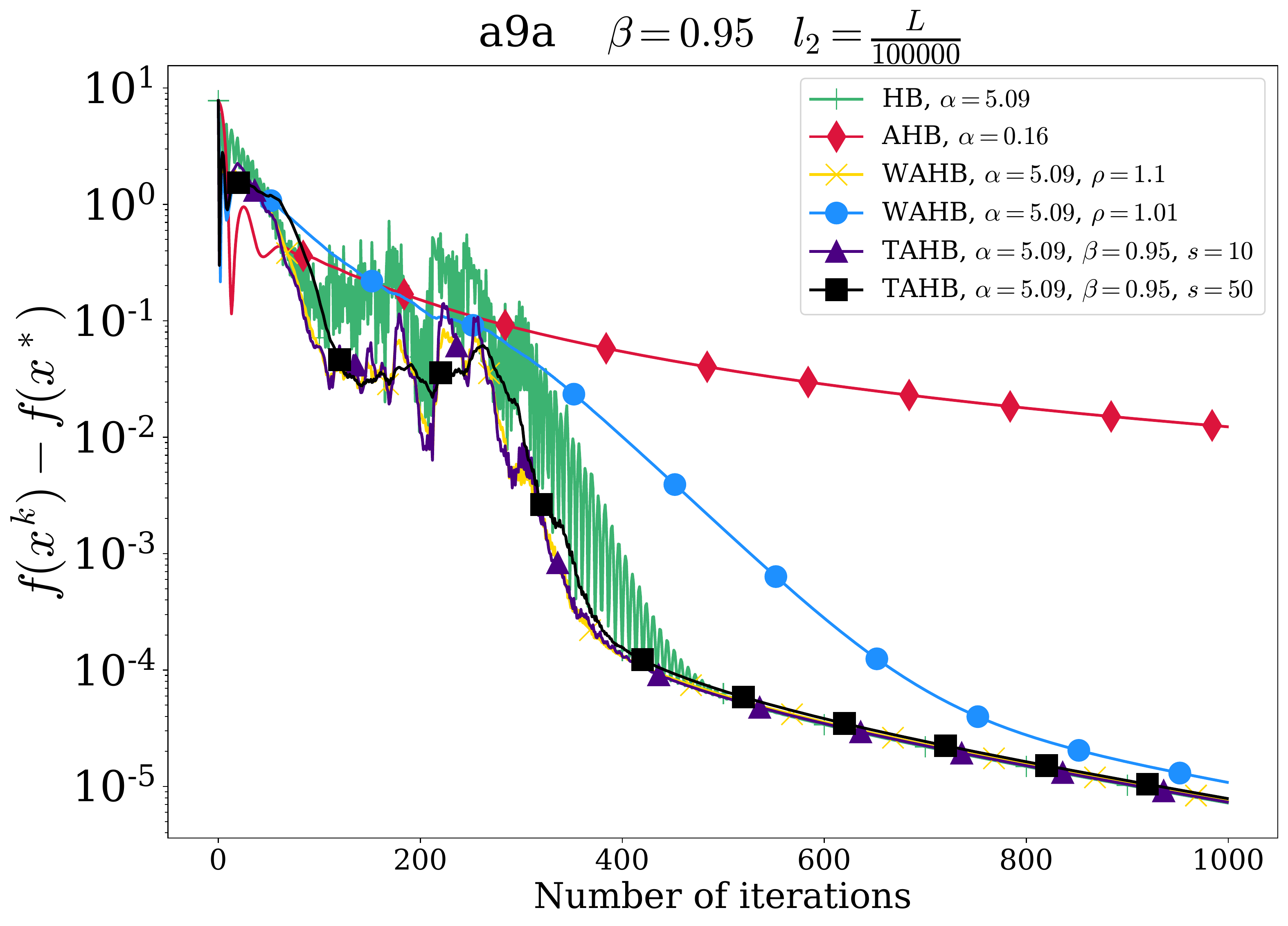}
    \includegraphics[width=0.325\textwidth]{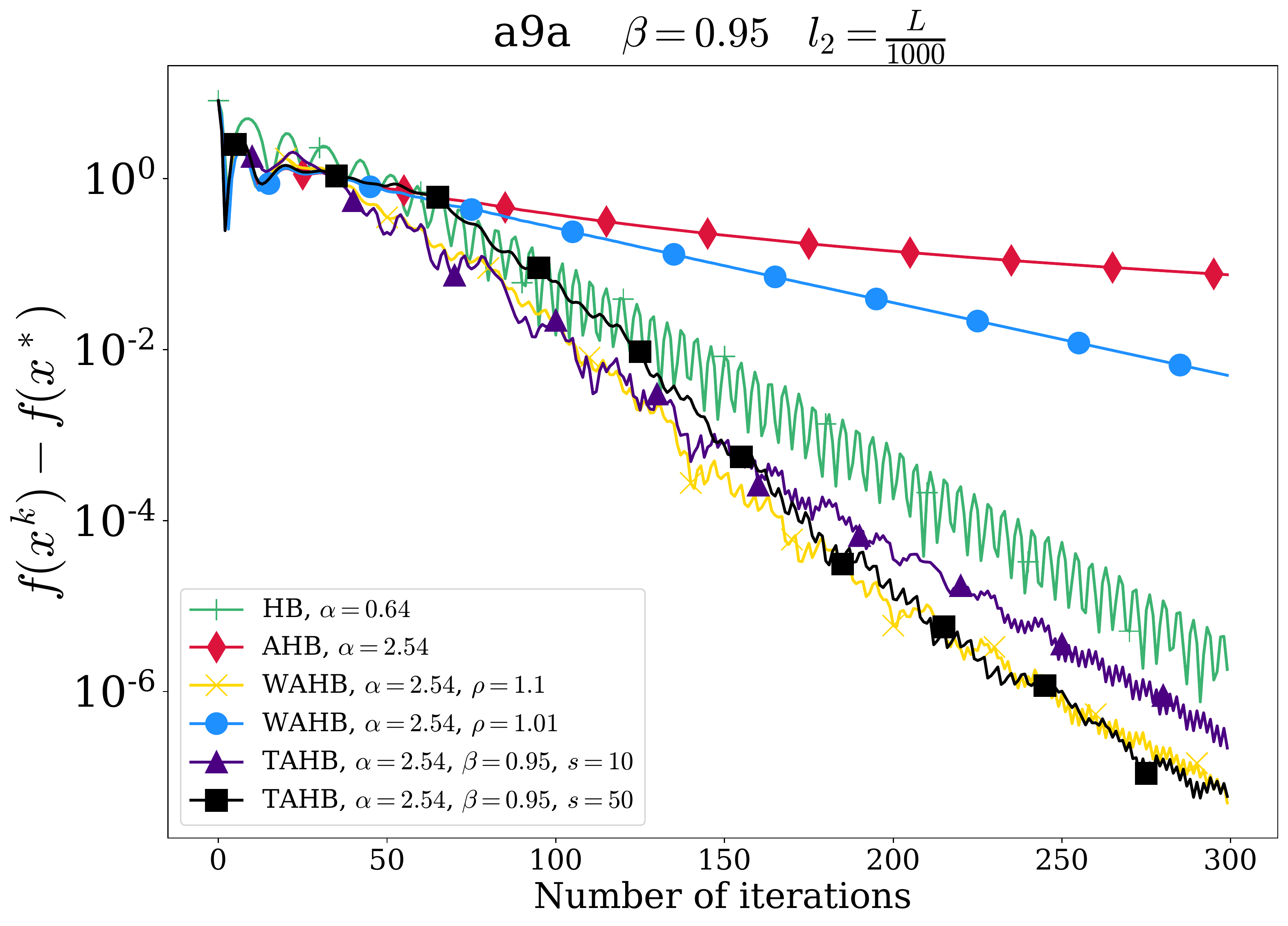}\\
    \includegraphics[width=0.325\textwidth]{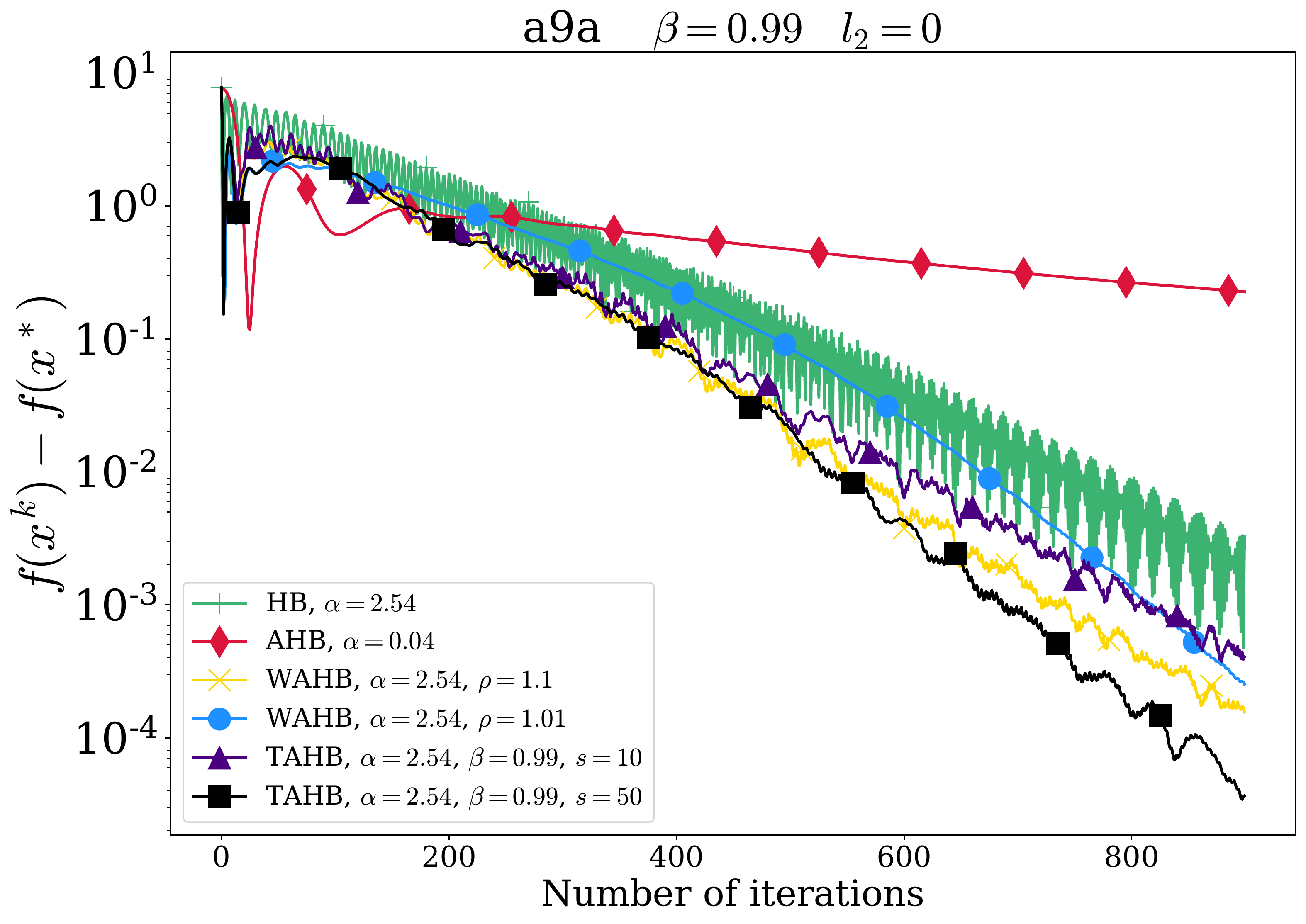}
    \includegraphics[width=0.325\textwidth]{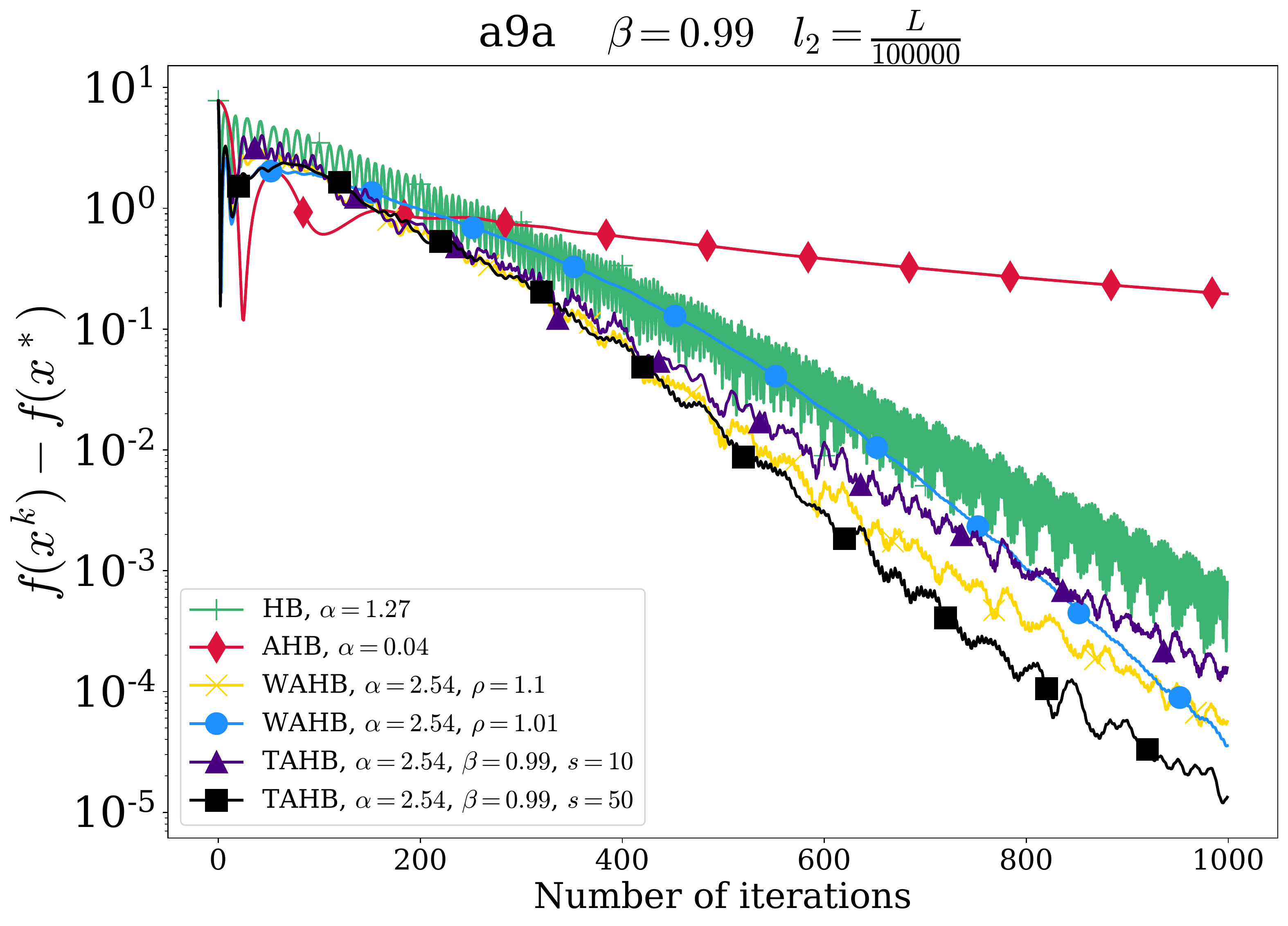}
    \includegraphics[width=0.325\textwidth]{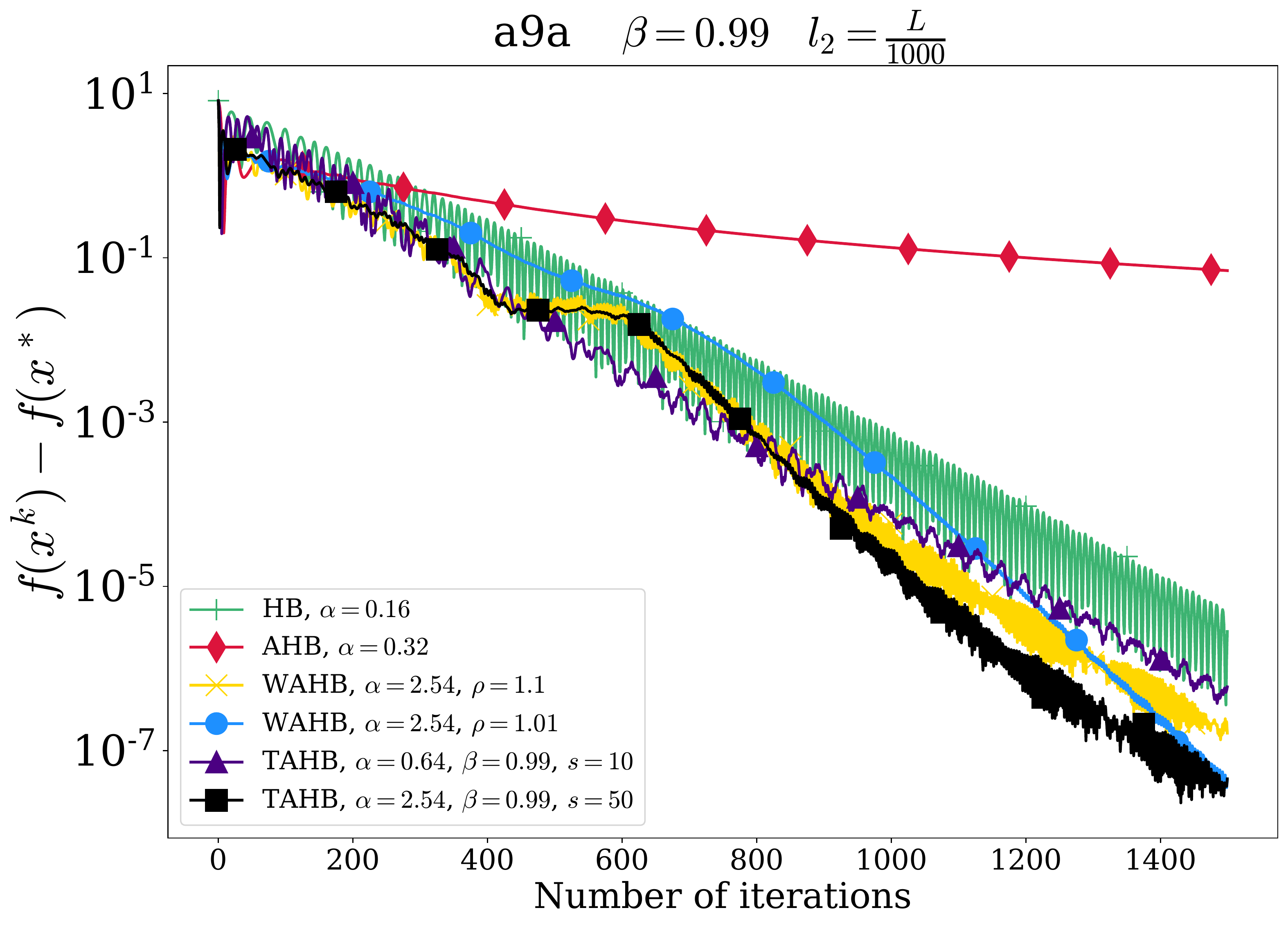}\\
    \includegraphics[width=0.325\textwidth]{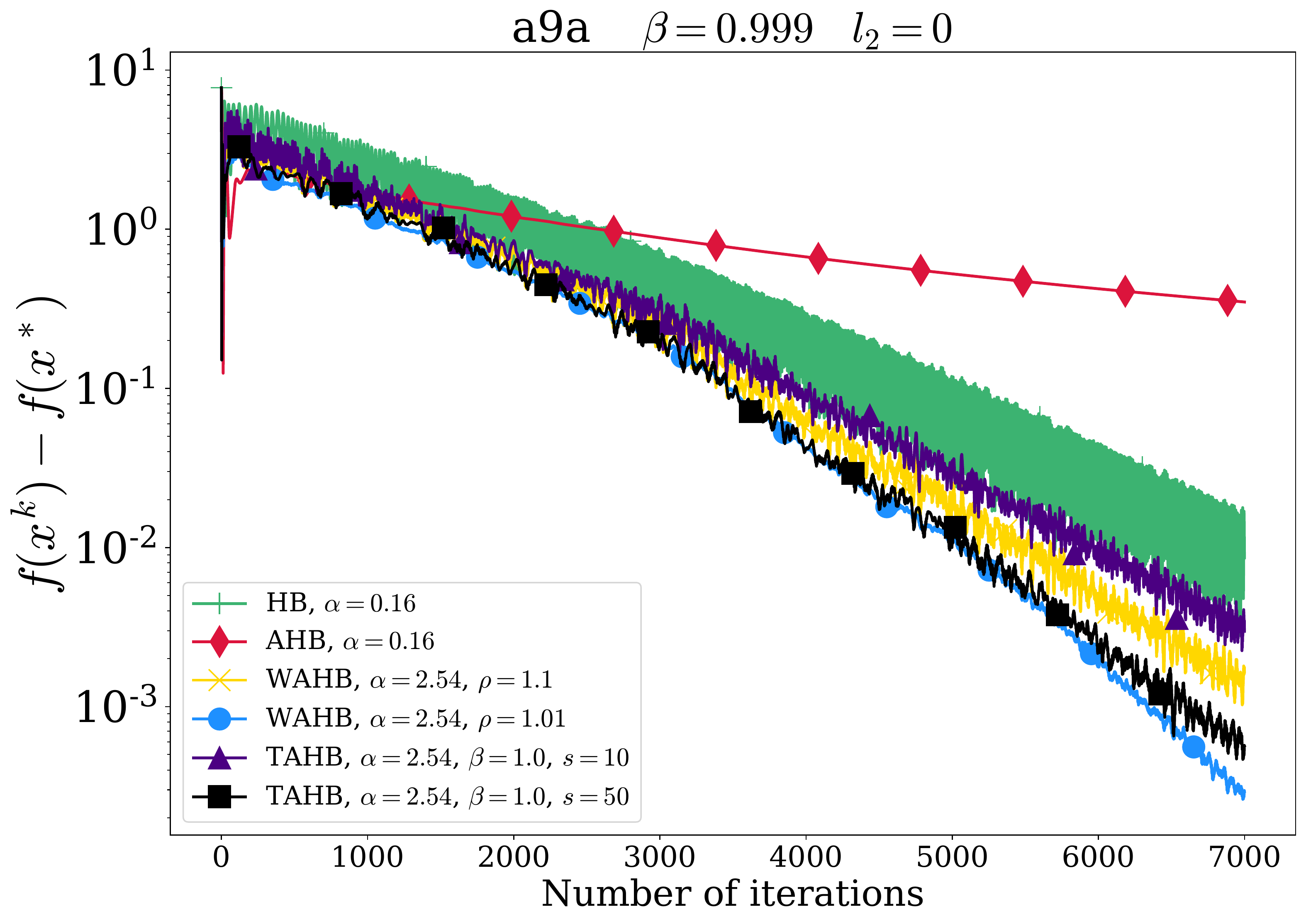}
    \includegraphics[width=0.325\textwidth]{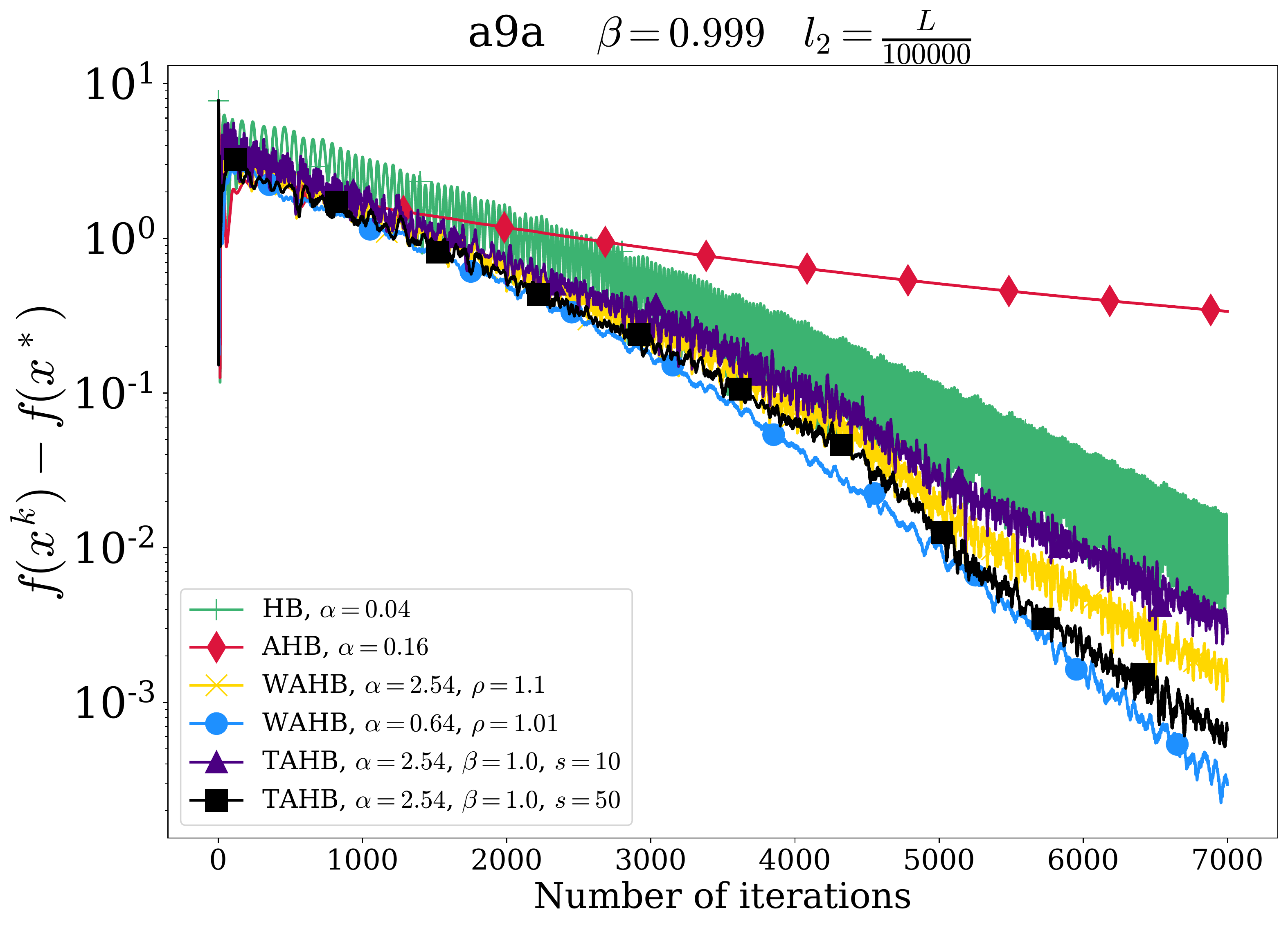}
    \includegraphics[width=0.325\textwidth]{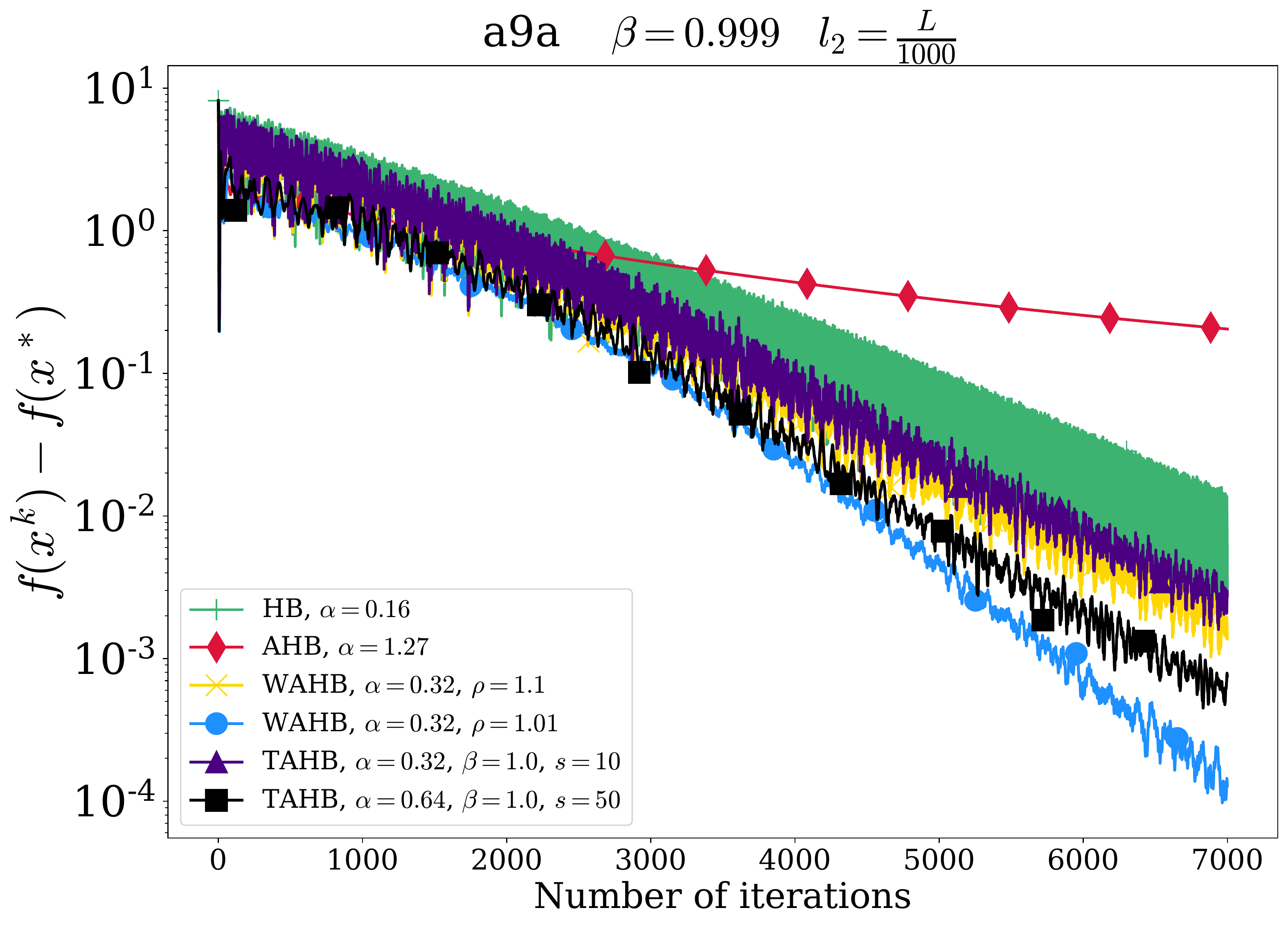}
    \caption{Trajectories of \algname{HB}, \algname{AHB}, \algname{WAHB}, and \algname{TAHB} with different momentum parameters $\beta$ applied to solve logistic regression problem with $\ell_2$-regularization for {\tt a9a} dataset. Stepsize $\alpha$ was tuned for each method and each choice of $\beta$ (and $\rho$, $s$) separately.}
    \label{fig:a9a_different_betas}
\end{figure}

\begin{figure}[H]
    \centering
    \includegraphics[width=0.325\textwidth]{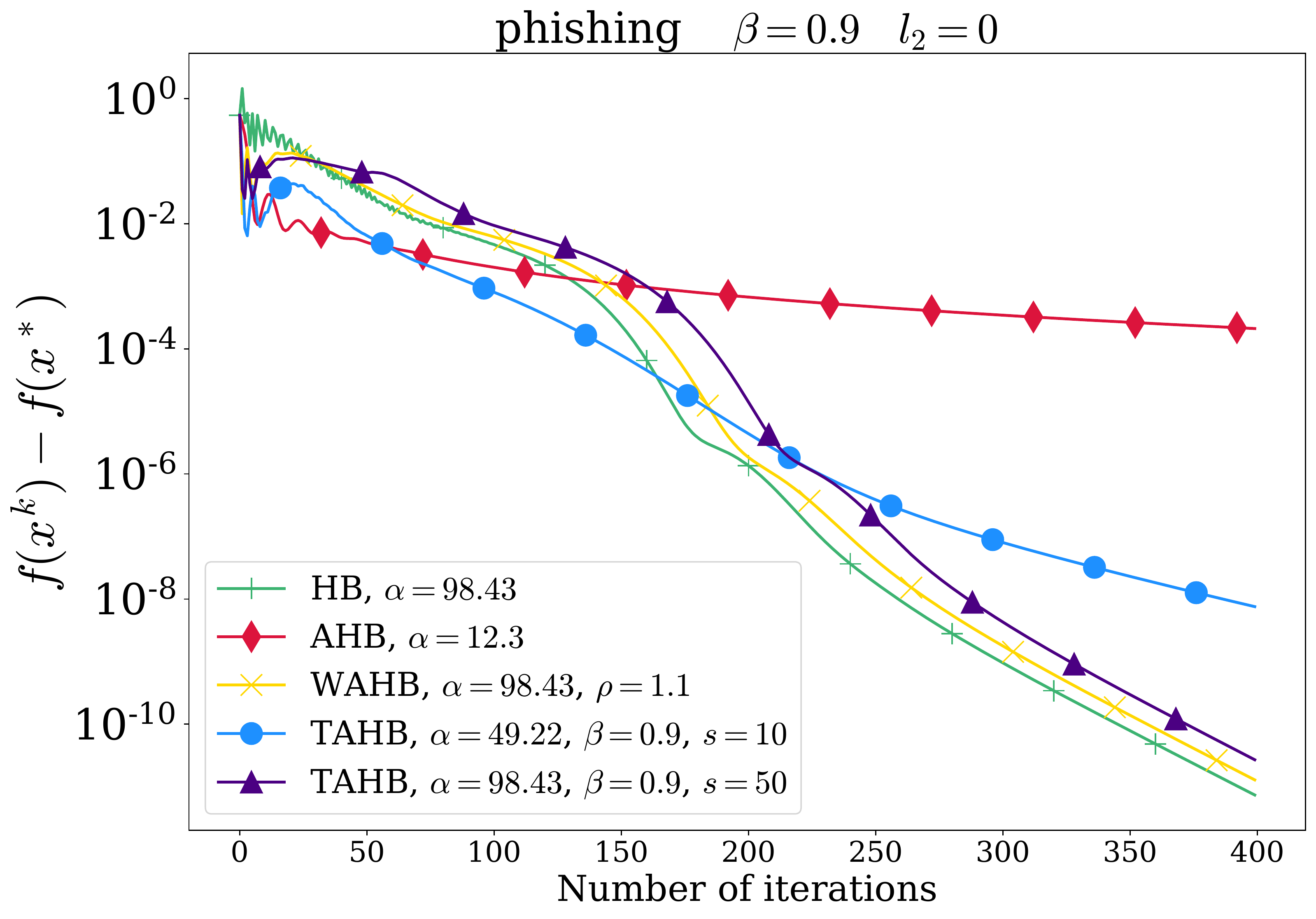}
    \includegraphics[width=0.325\textwidth]{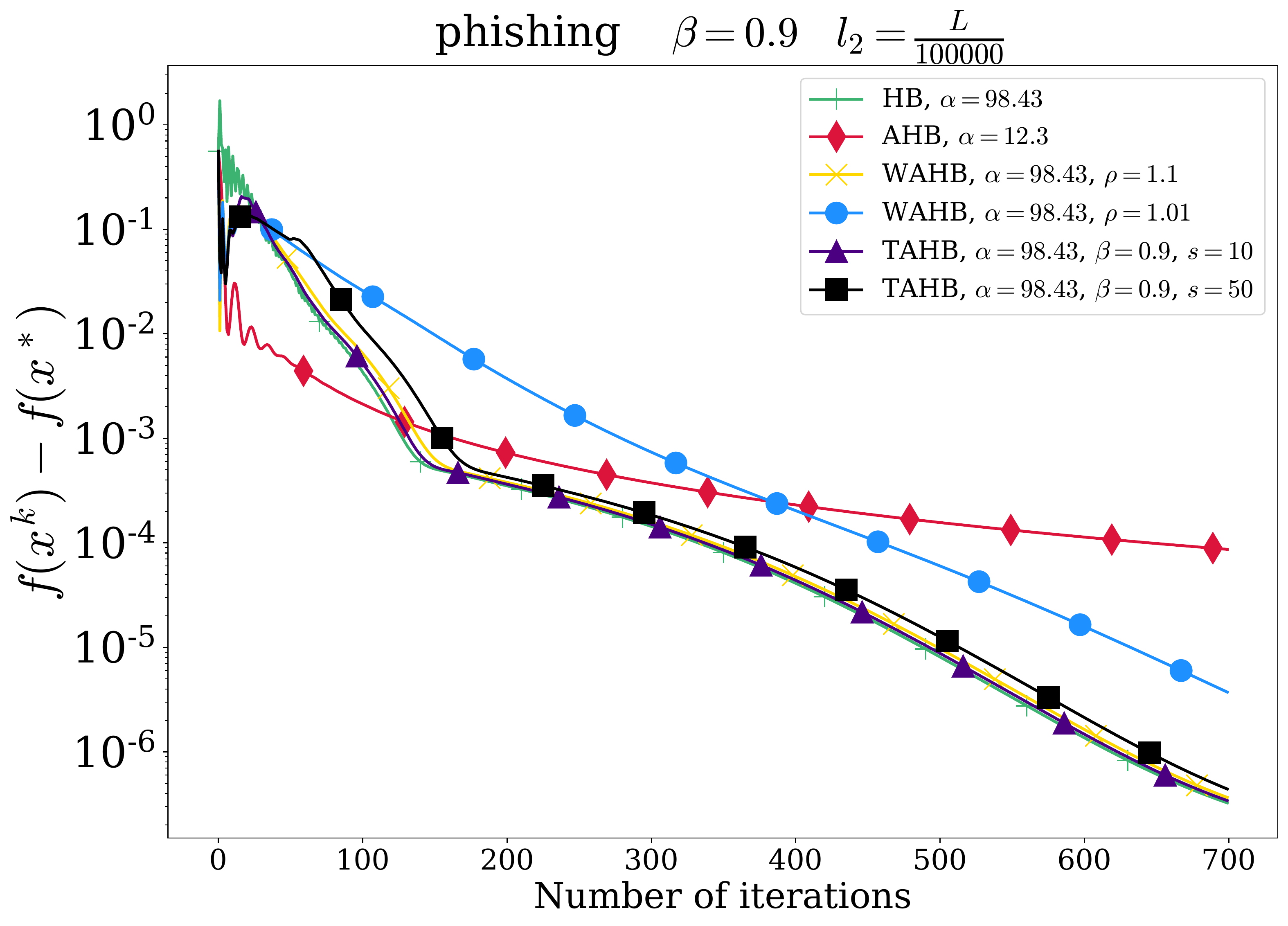}
    \includegraphics[width=0.325\textwidth]{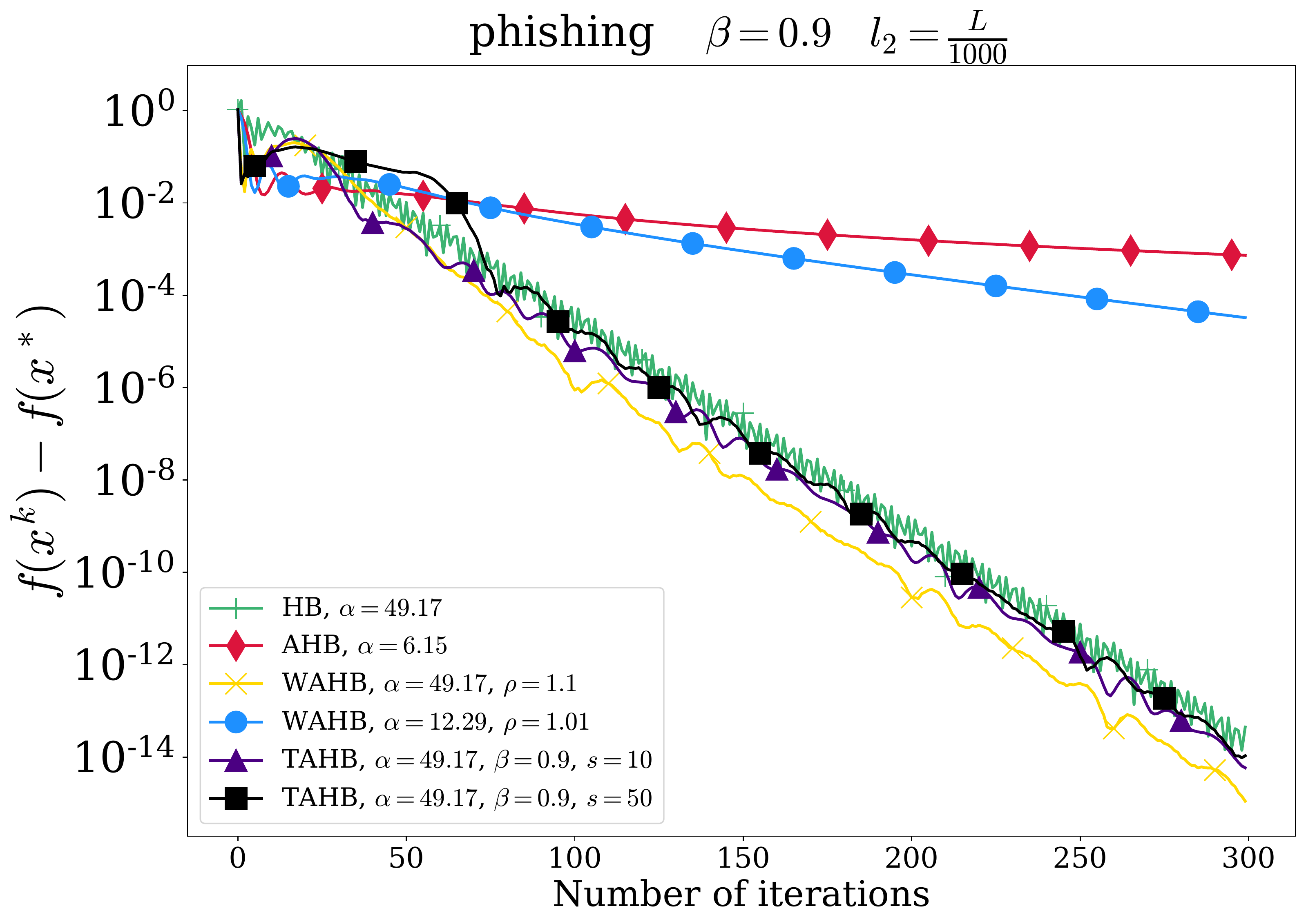}\\
    \includegraphics[width=0.325\textwidth]{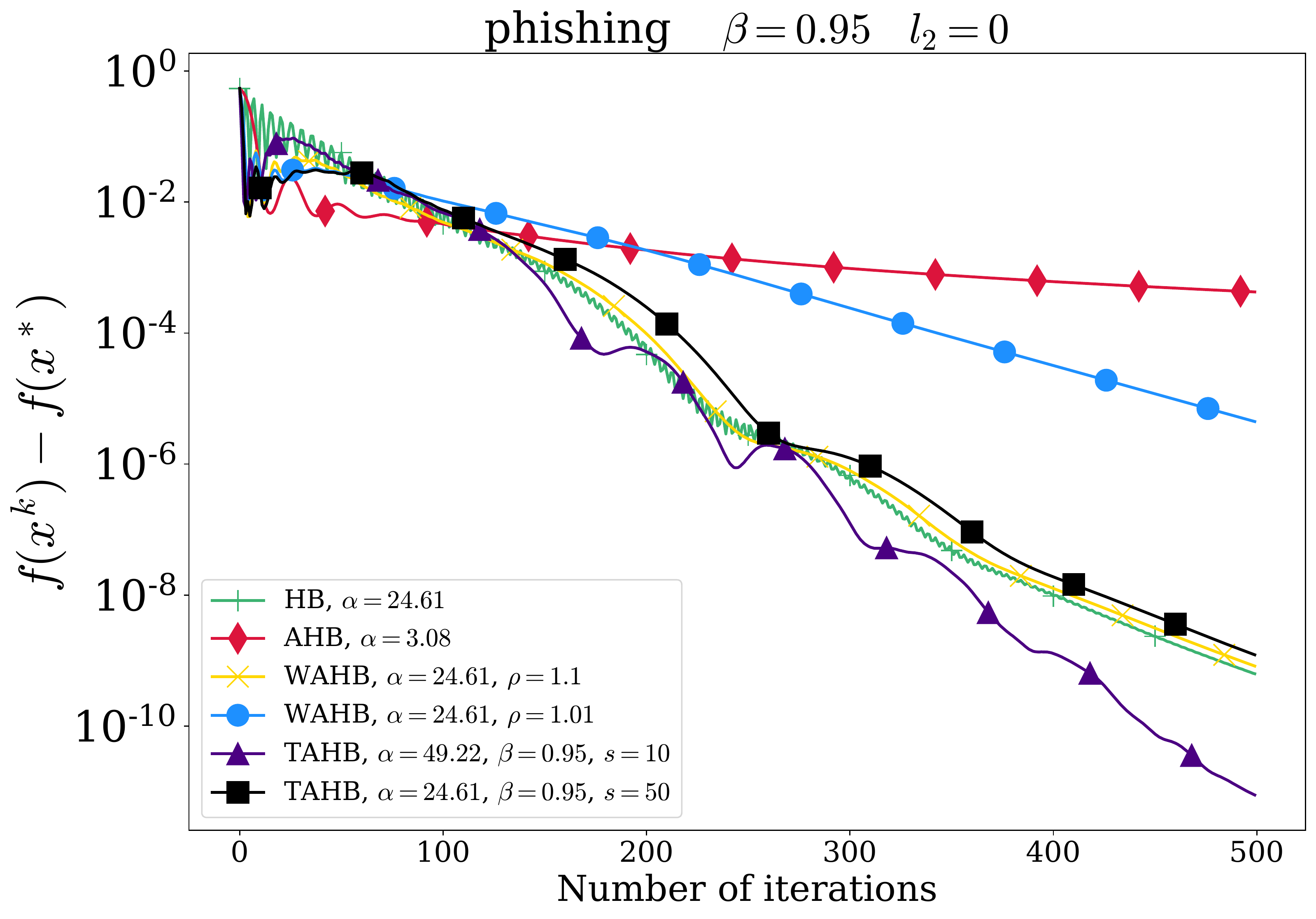}
    \includegraphics[width=0.325\textwidth]{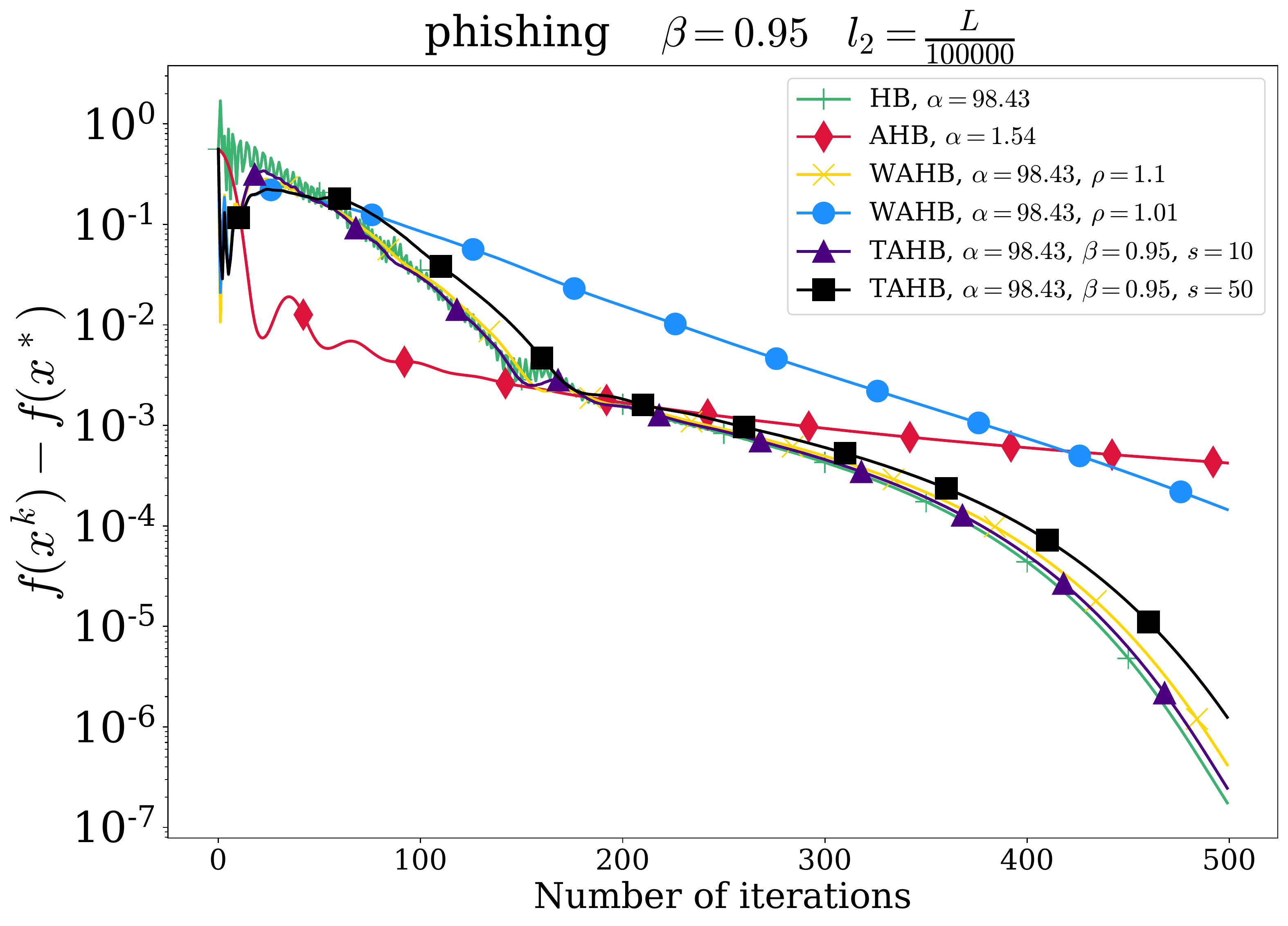}
    \includegraphics[width=0.325\textwidth]{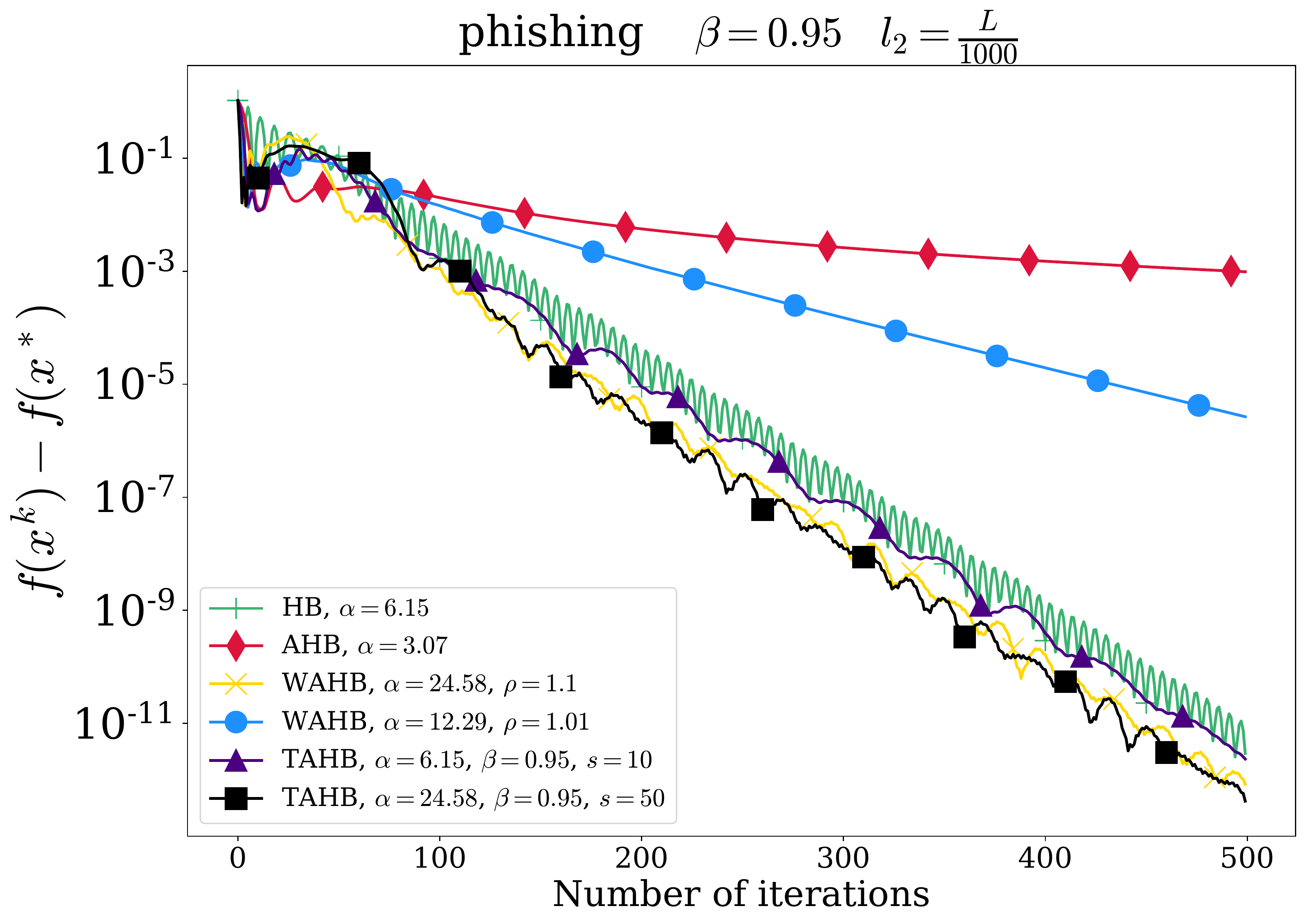}\\
    \includegraphics[width=0.325\textwidth]{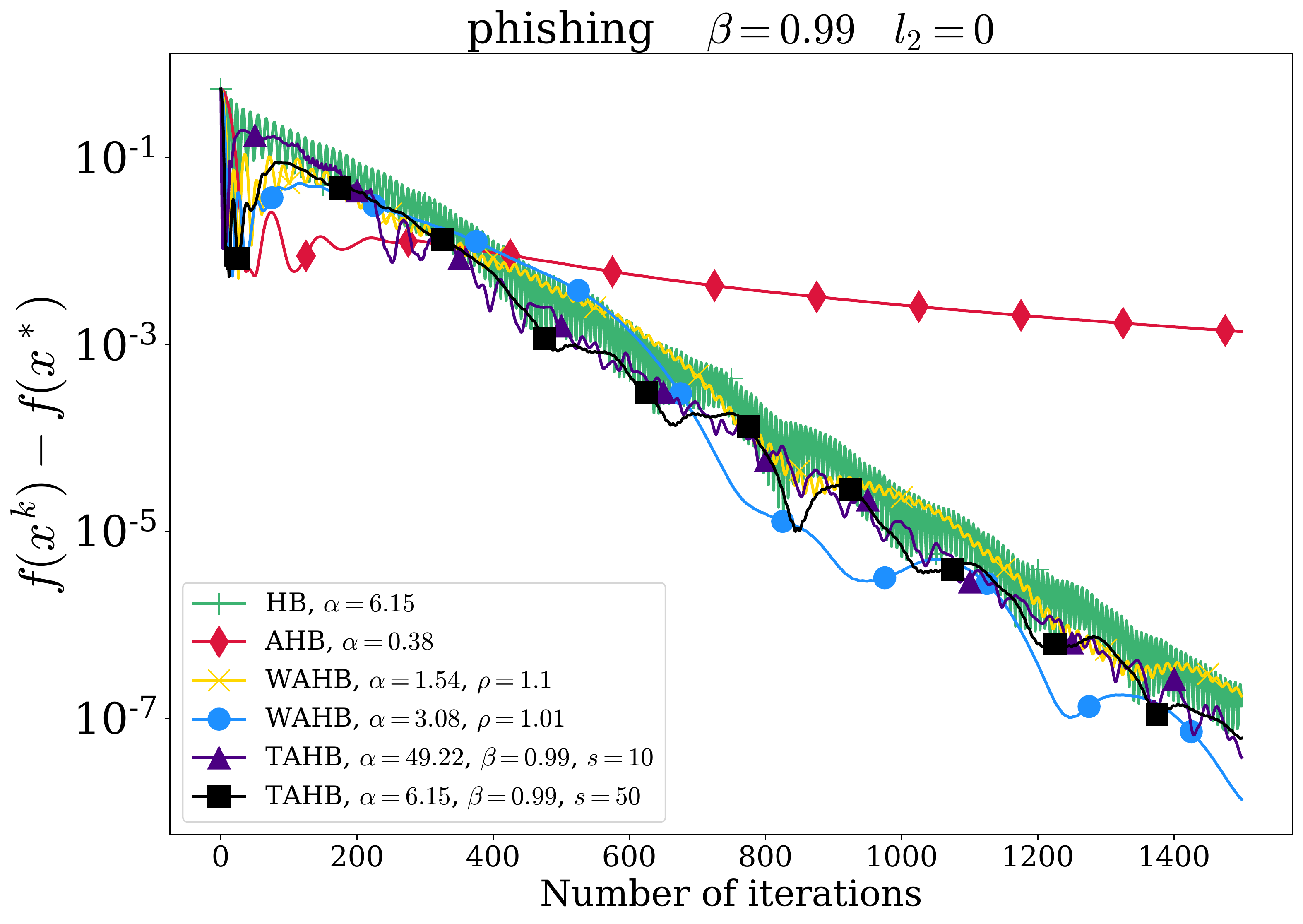}
    \includegraphics[width=0.325\textwidth]{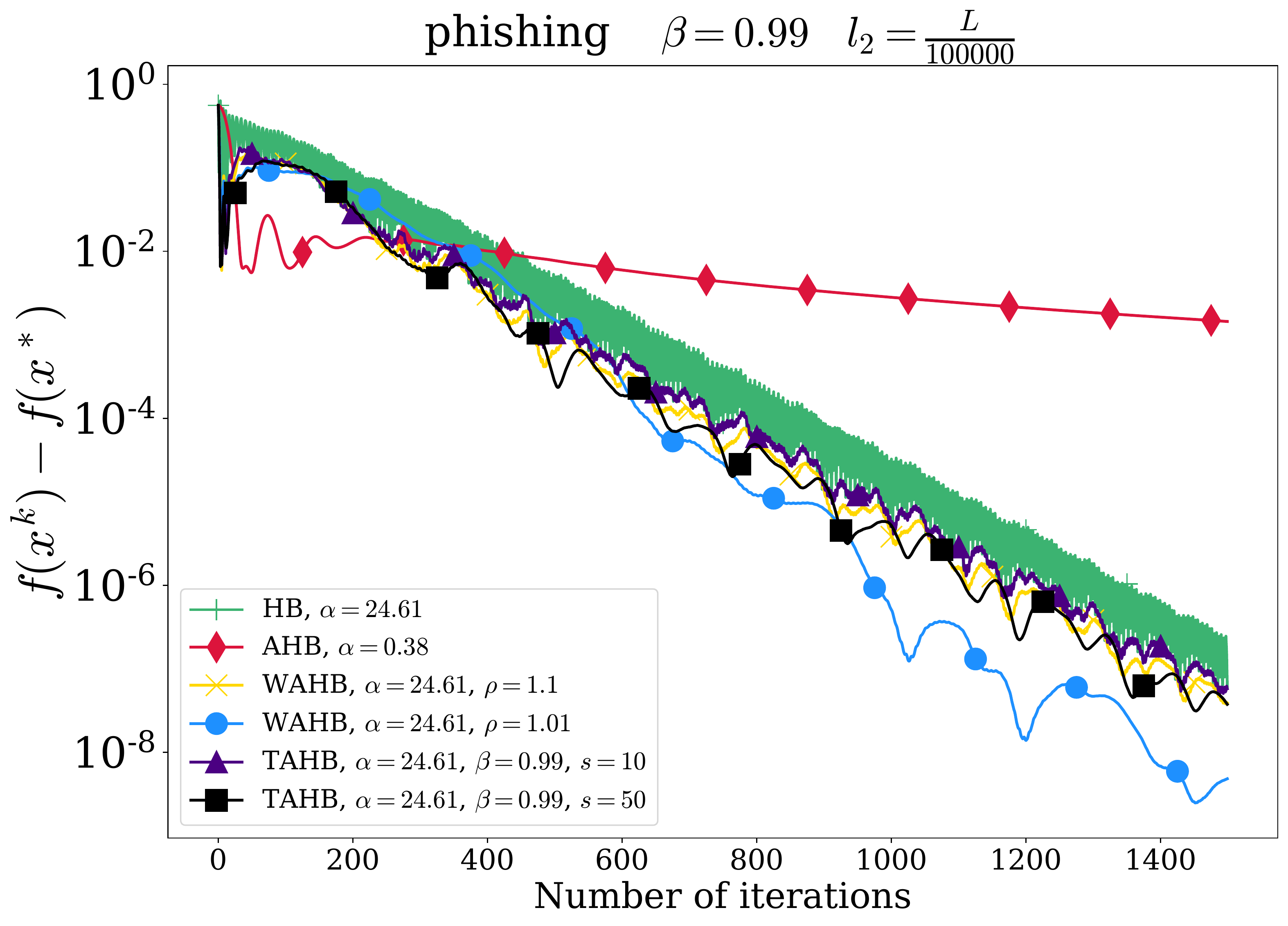}
    \includegraphics[width=0.325\textwidth]{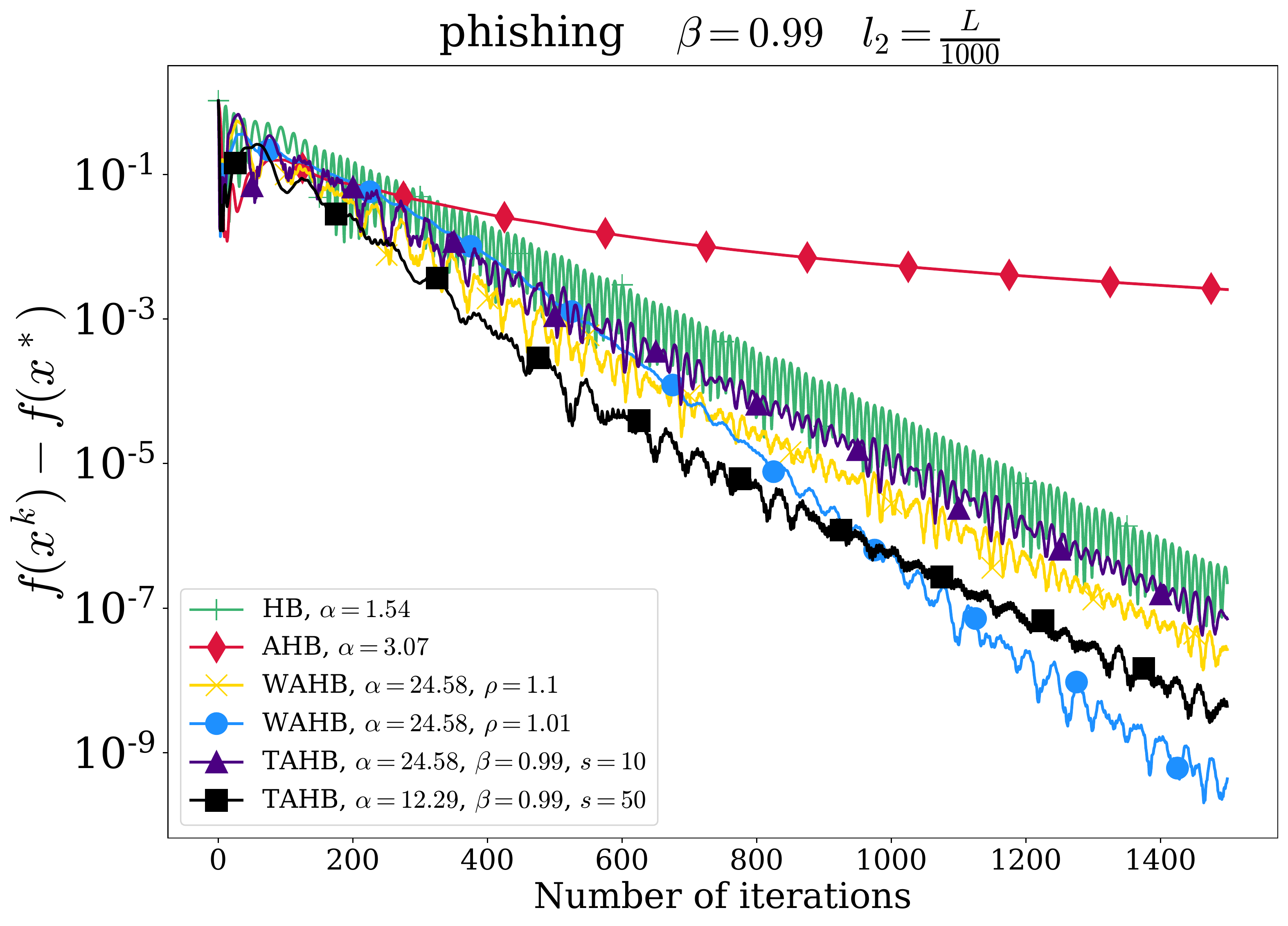}\\
    \includegraphics[width=0.325\textwidth]{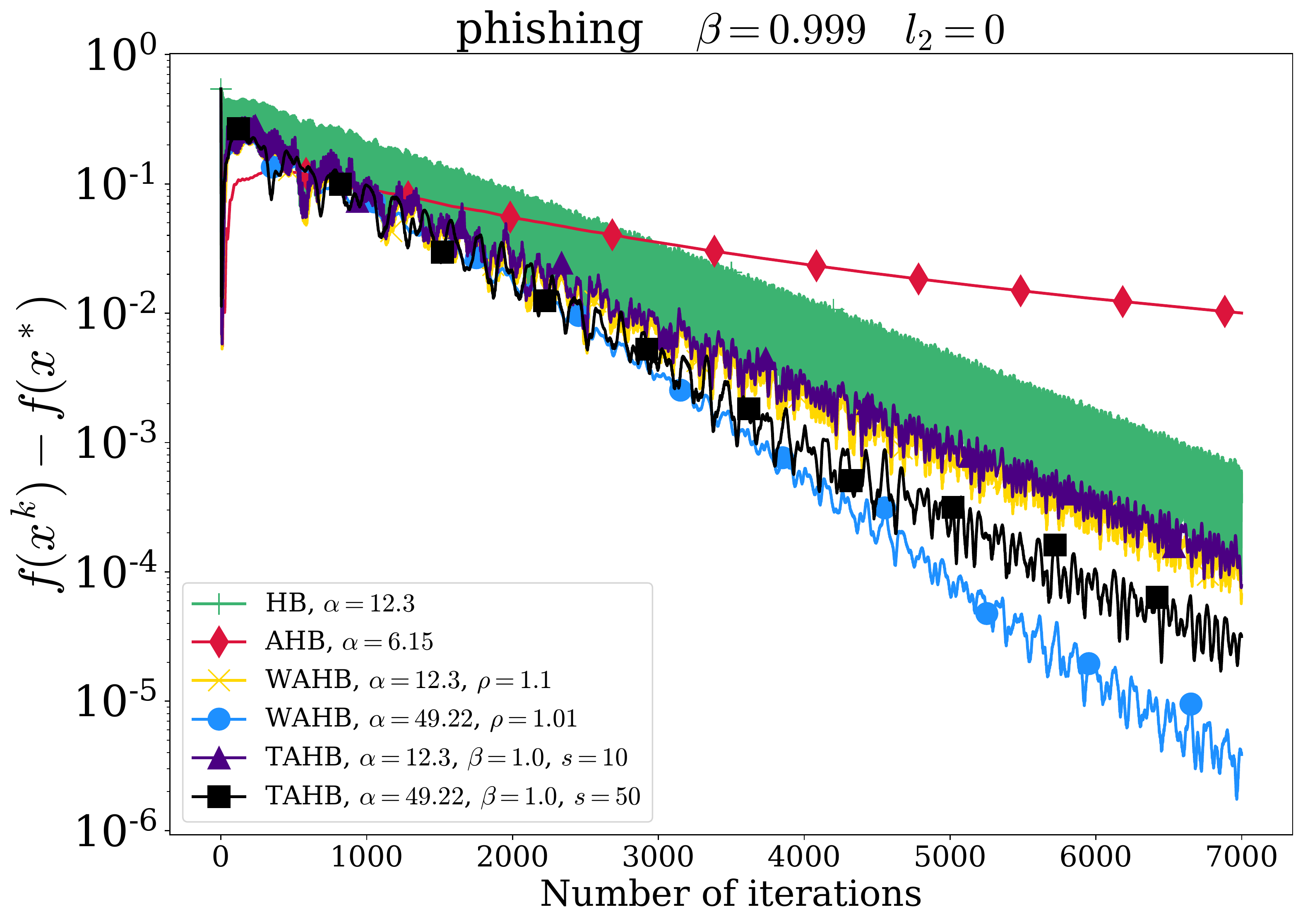}
    \includegraphics[width=0.325\textwidth]{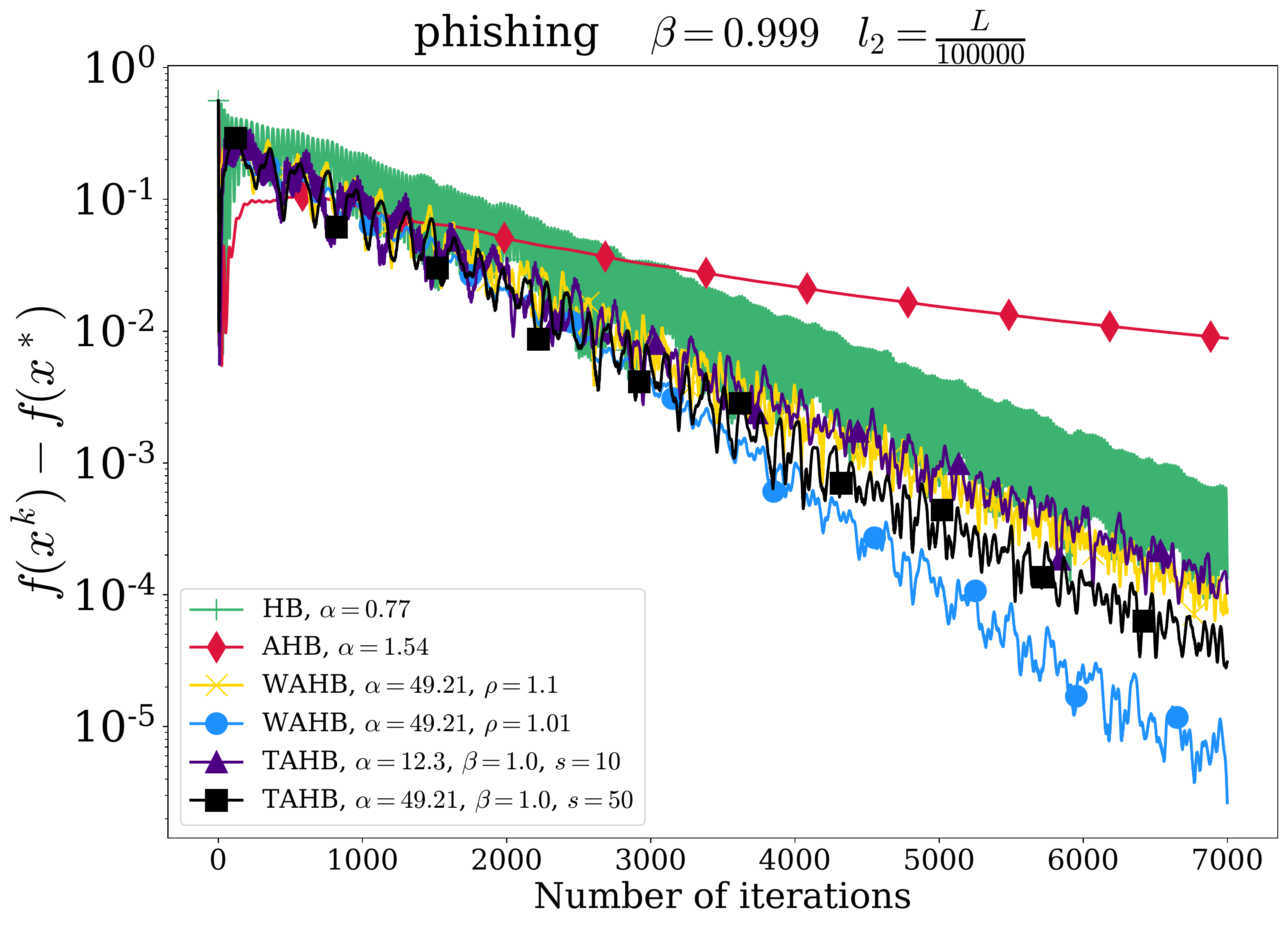}
    \includegraphics[width=0.325\textwidth]{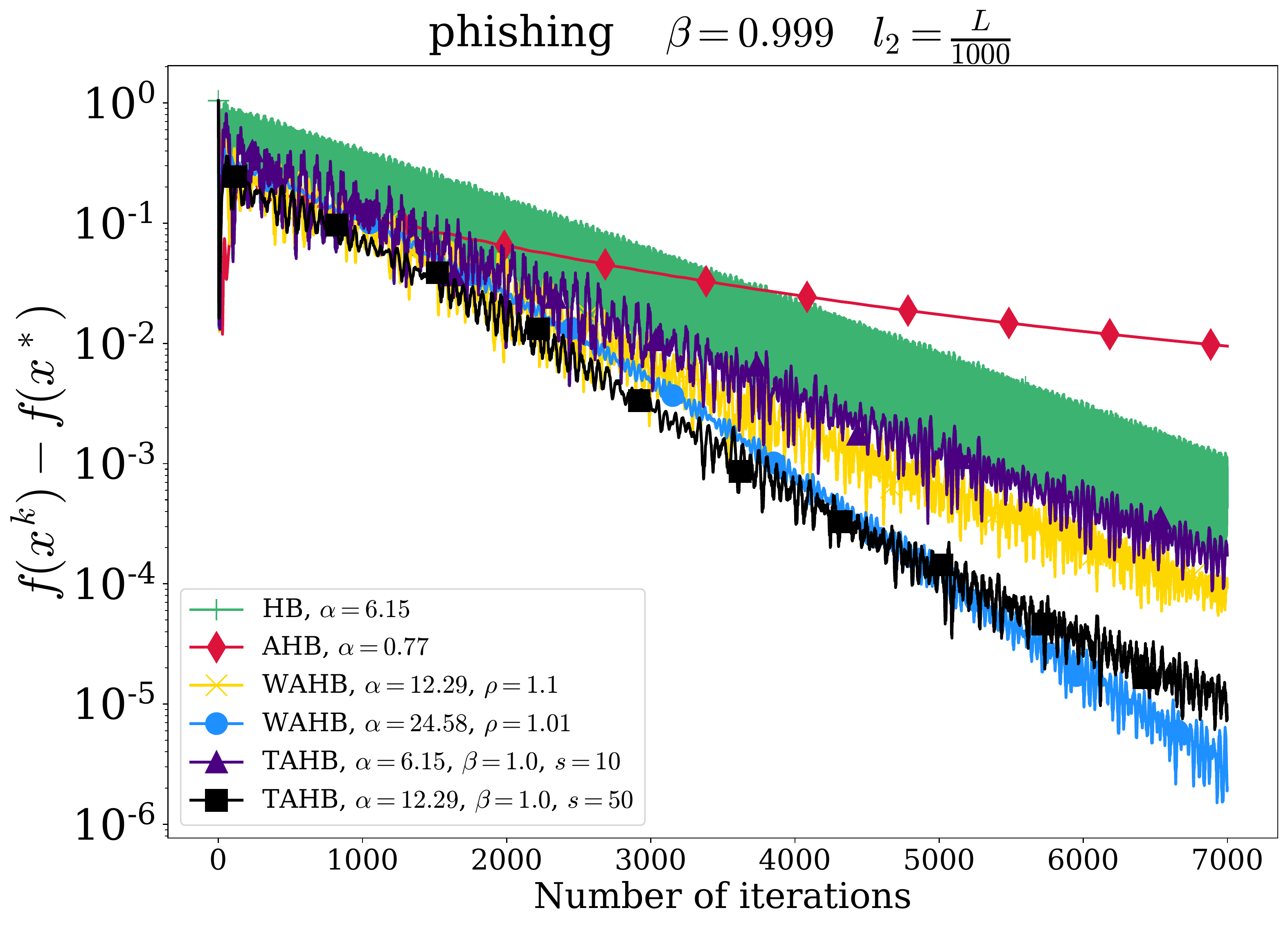}
    \caption{Trajectories of \algname{HB}, \algname{AHB}, \algname{WAHB}, and \algname{TAHB} with different momentum parameters $\beta$ applied to solve logistic regression problem with $\ell_2$-regularization for {\tt phishing} dataset. Stepsize $\alpha$ was tuned for each method and each choice of $\beta$ (and $\rho$, $s$) separately.}
    \label{fig:phishing_different_betas}
\end{figure}

\begin{figure}[H]
    \centering
    \includegraphics[width=0.325\textwidth]{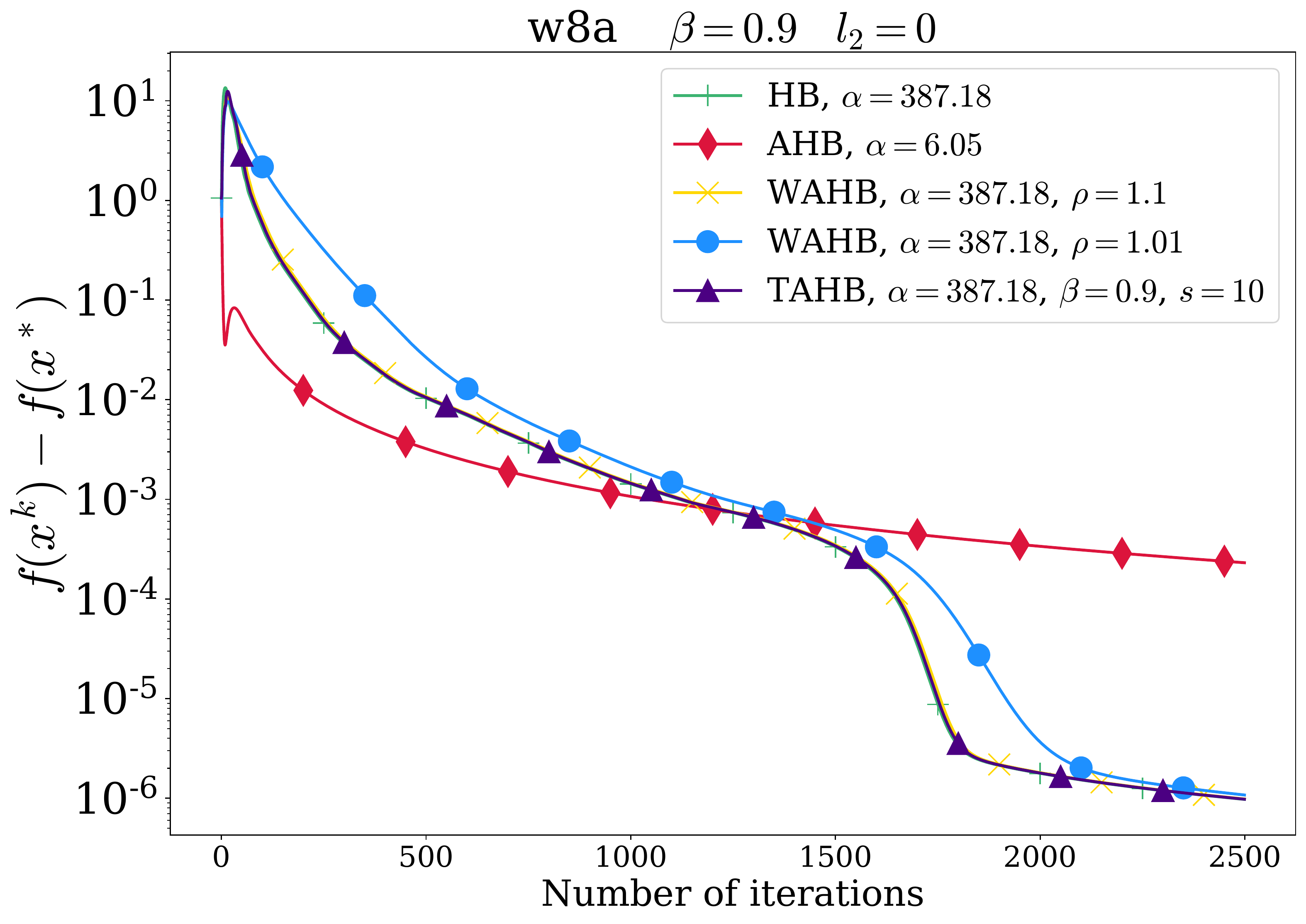}
    \includegraphics[width=0.325\textwidth]{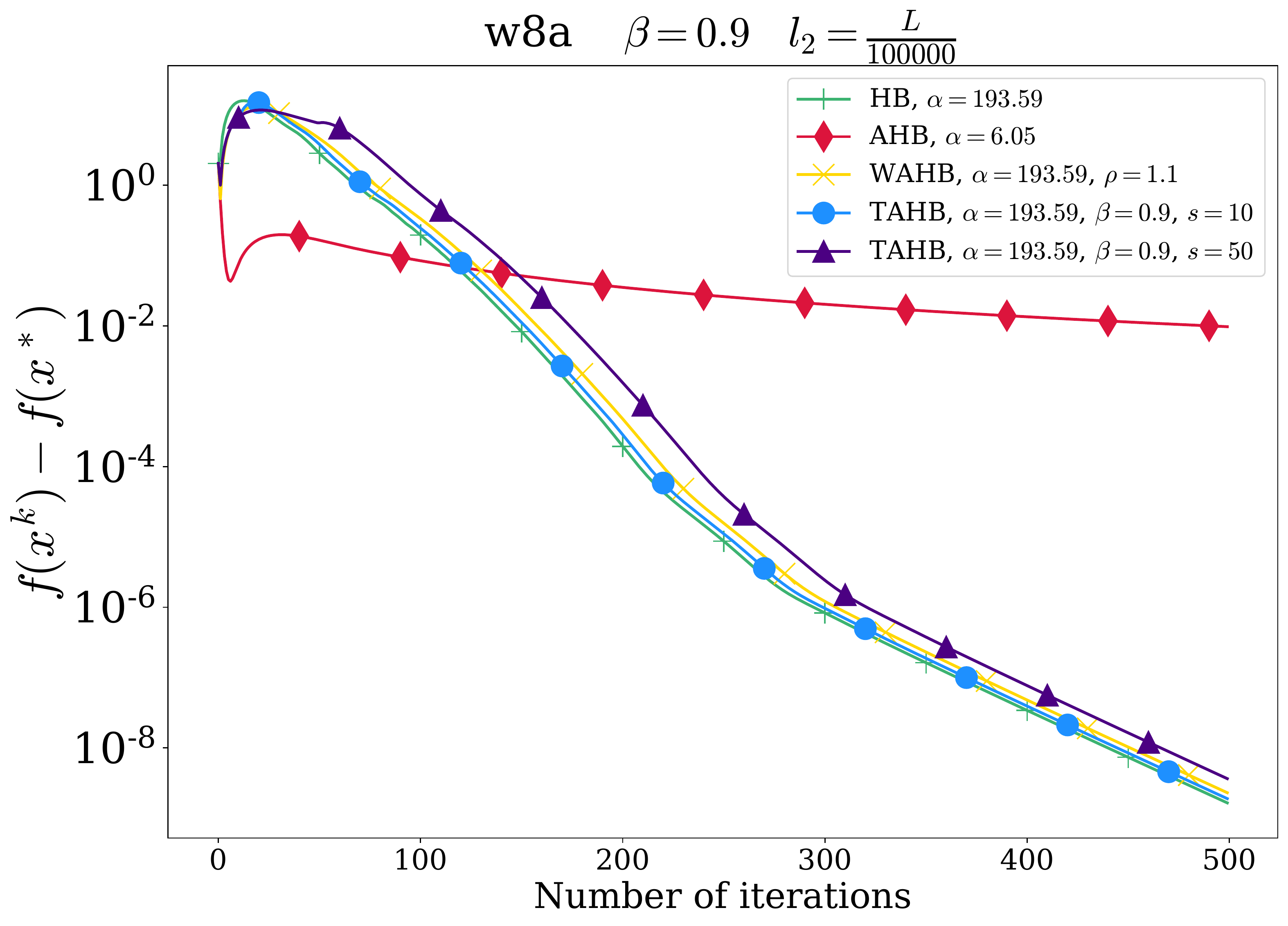}
    \includegraphics[width=0.325\textwidth]{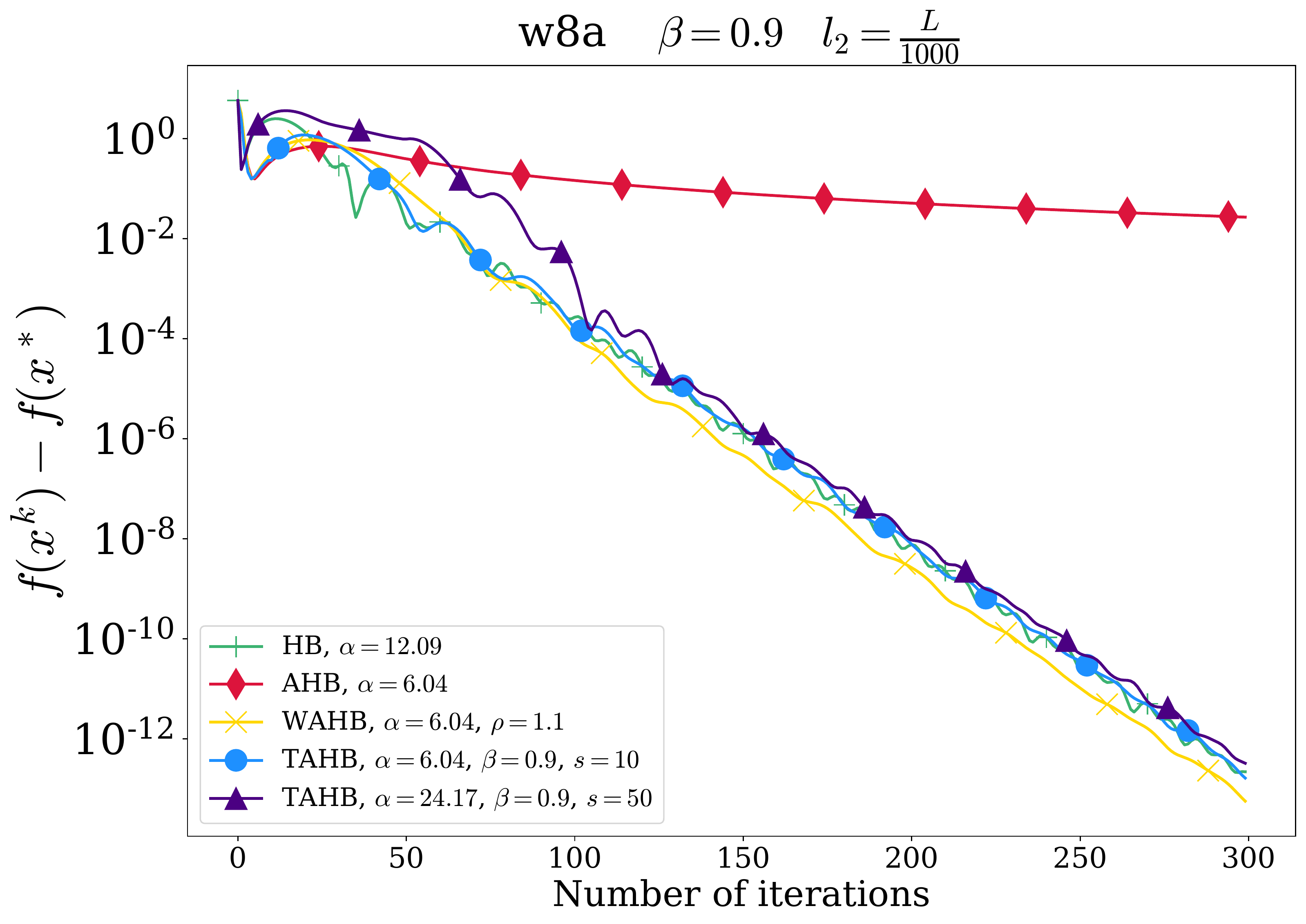}\\
    \includegraphics[width=0.325\textwidth]{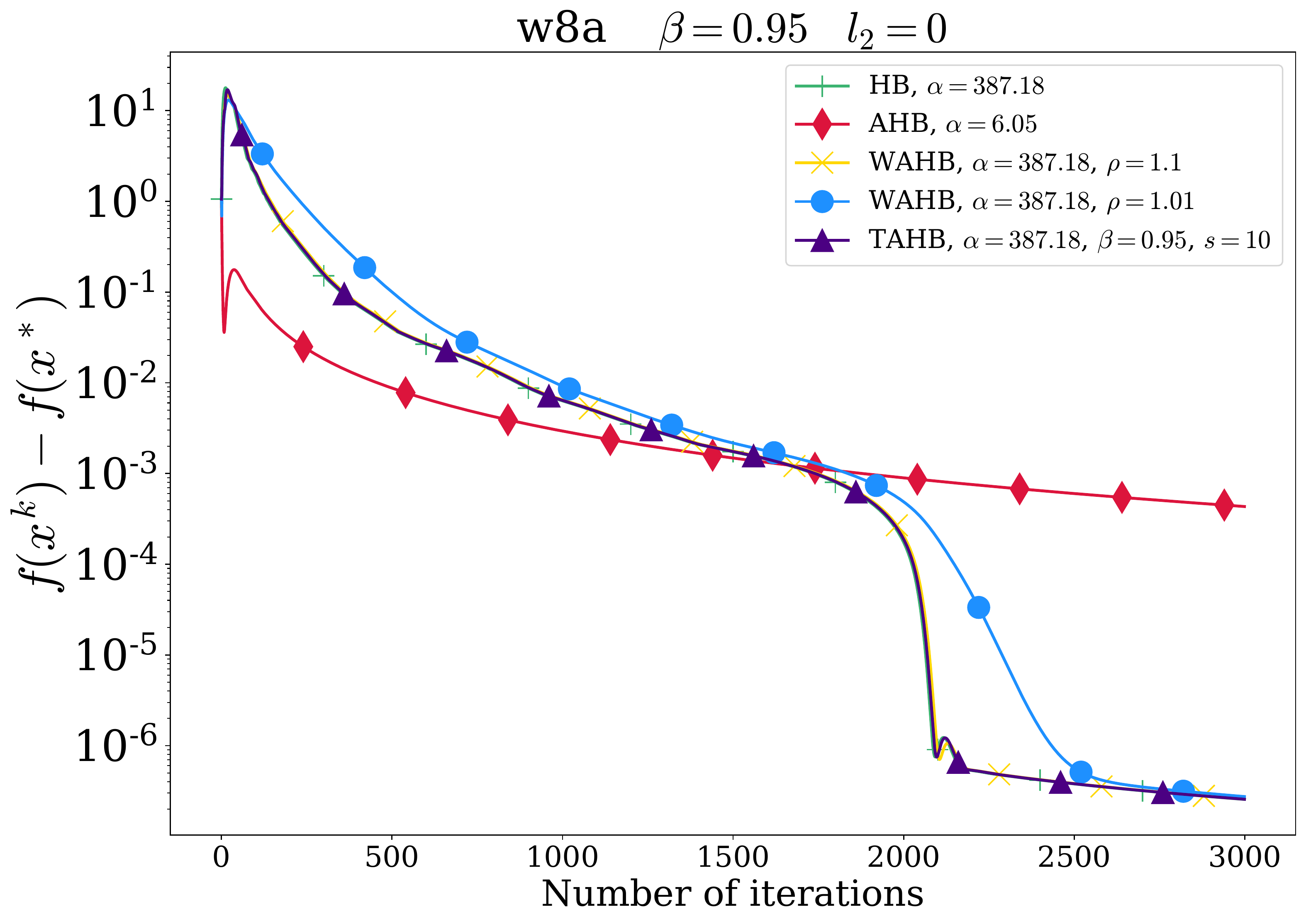}
    \includegraphics[width=0.325\textwidth]{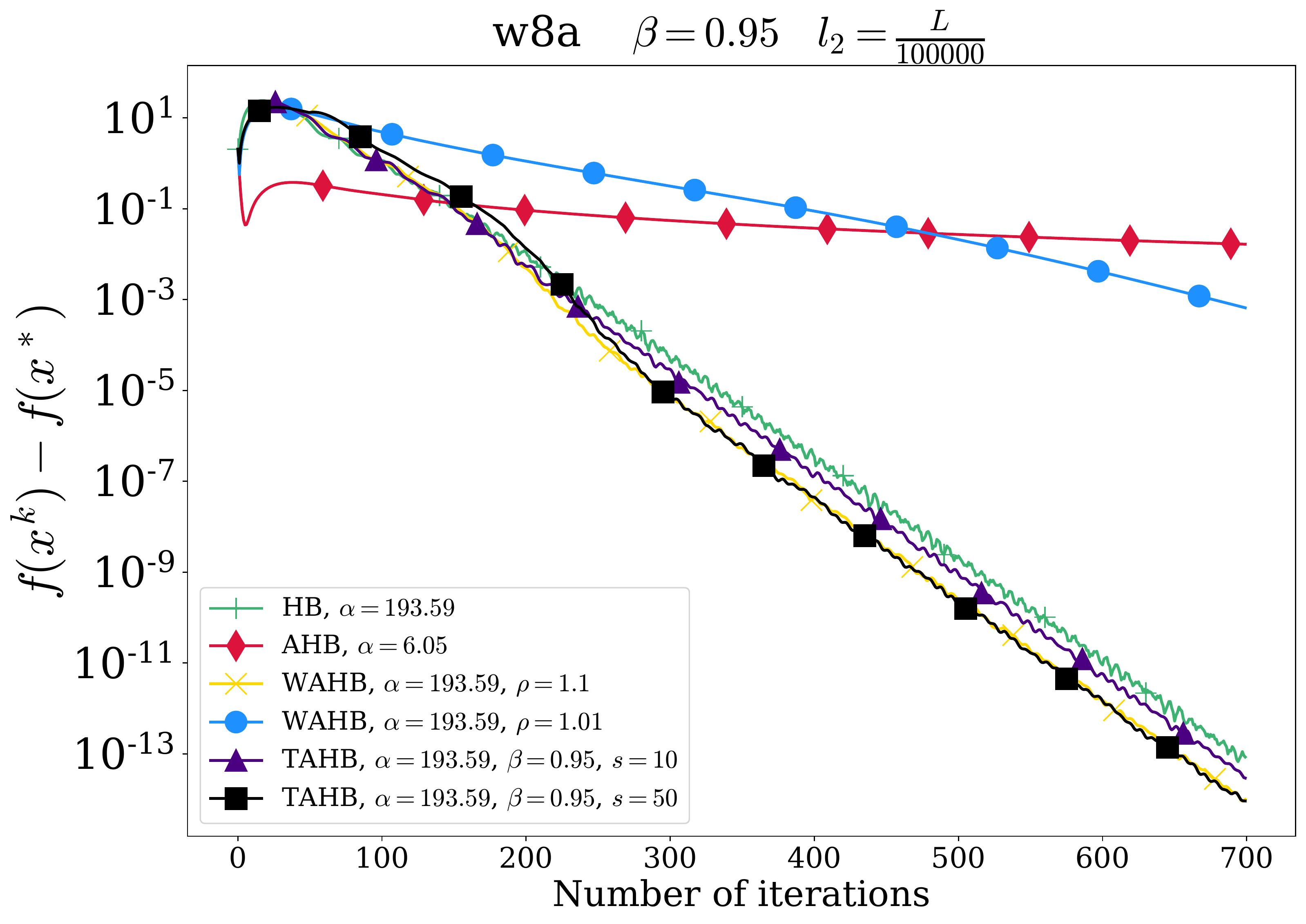}
    \includegraphics[width=0.325\textwidth]{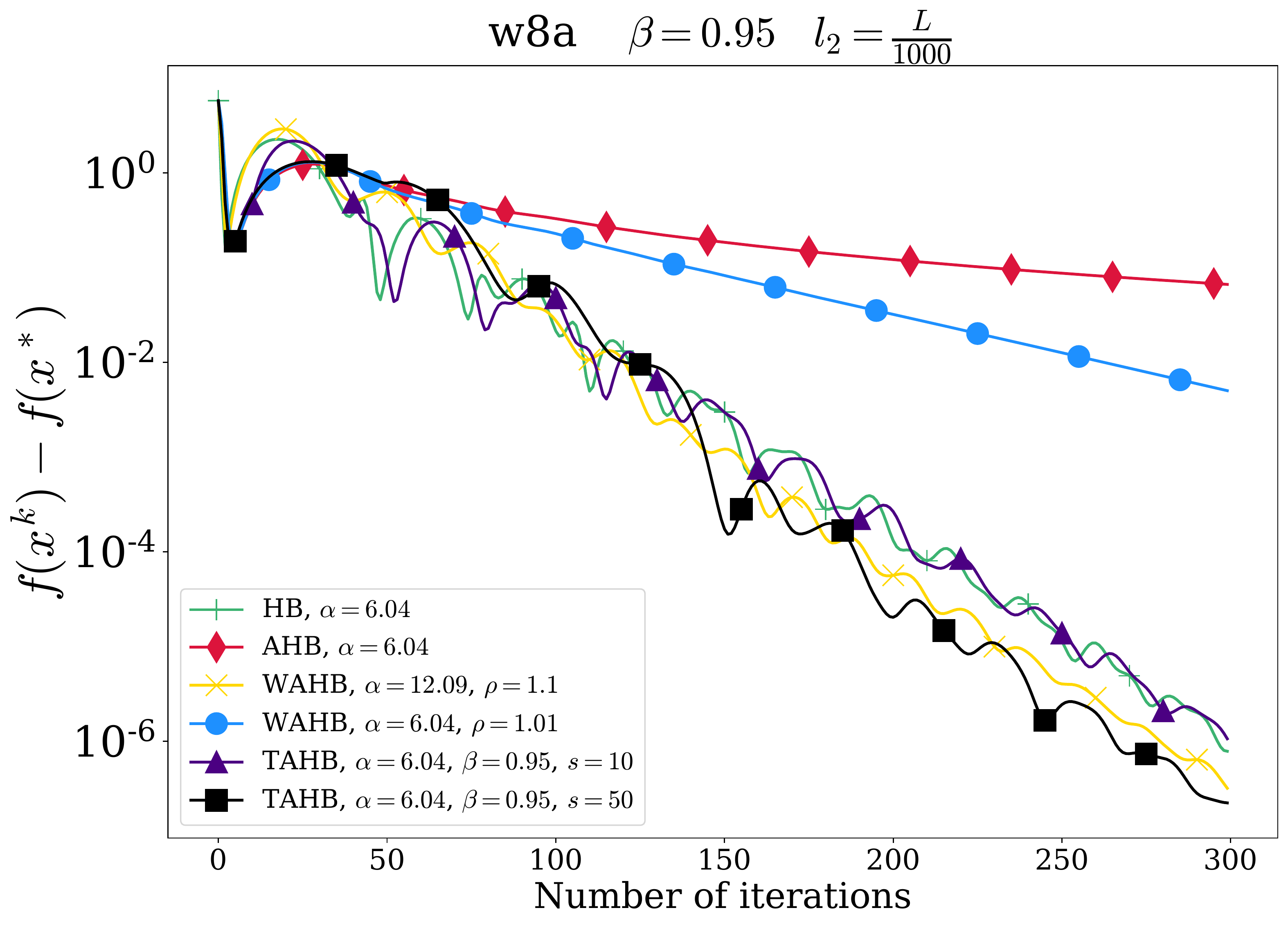}\\
    \includegraphics[width=0.325\textwidth]{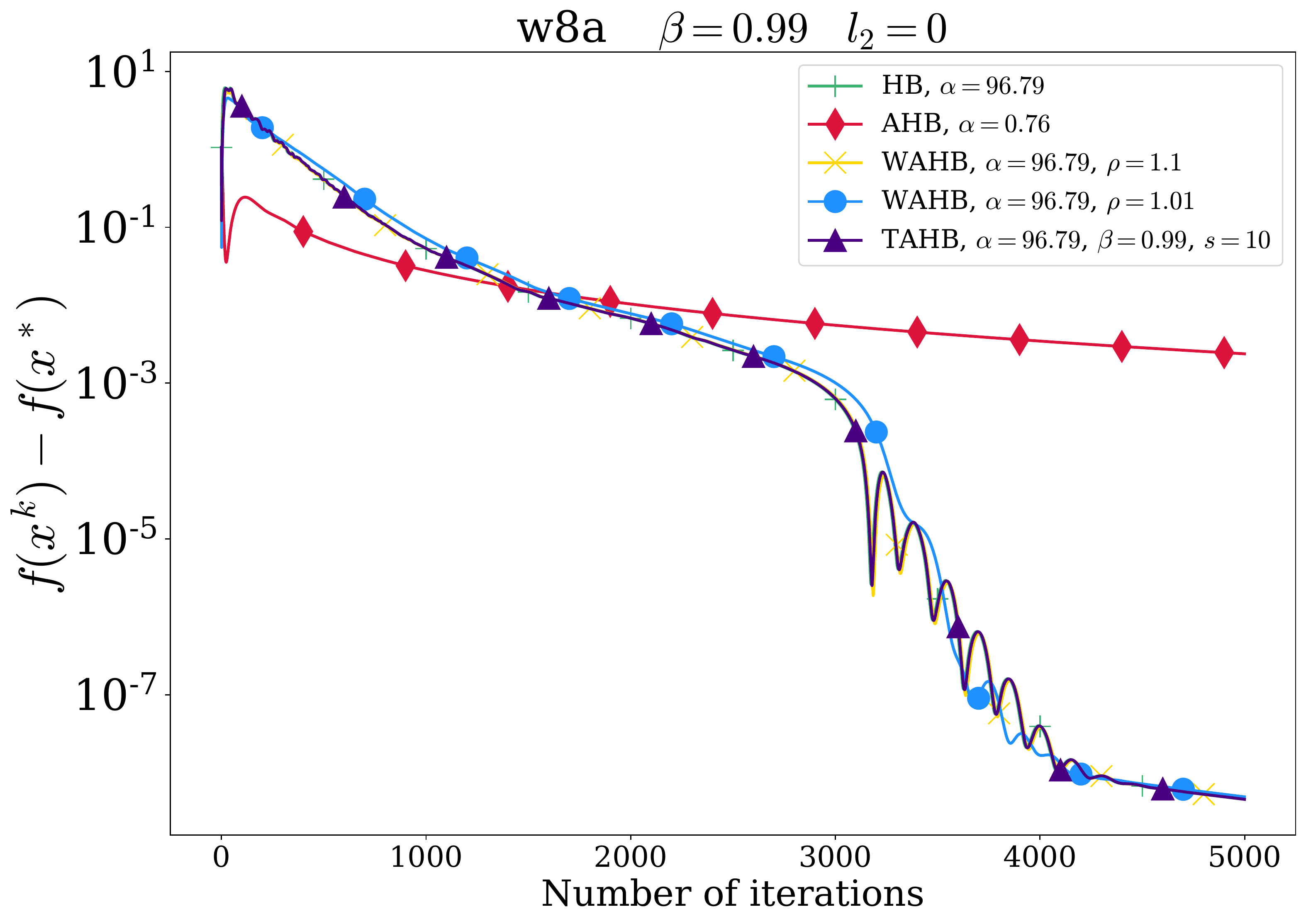}
    \includegraphics[width=0.325\textwidth]{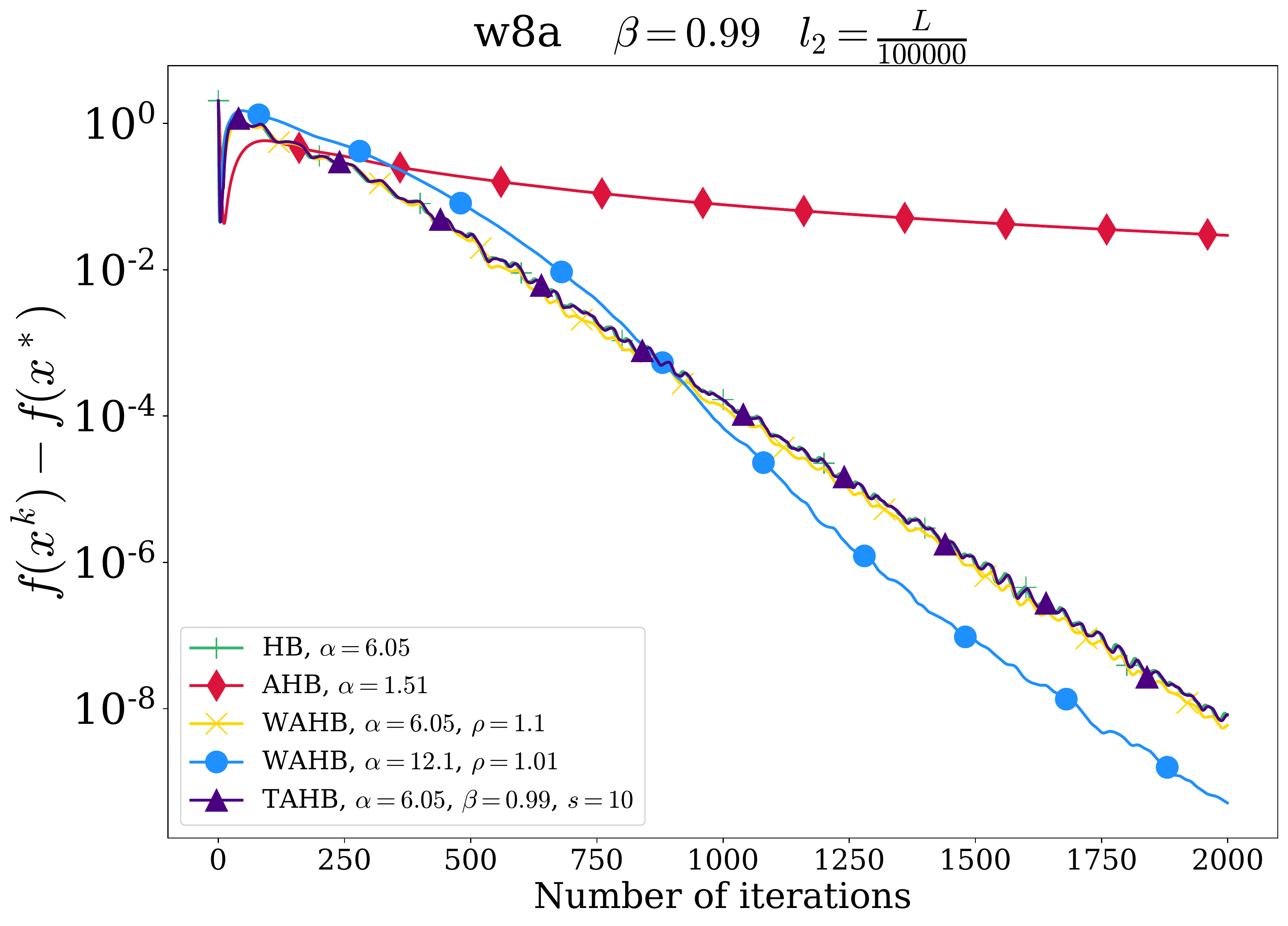}
    \includegraphics[width=0.325\textwidth]{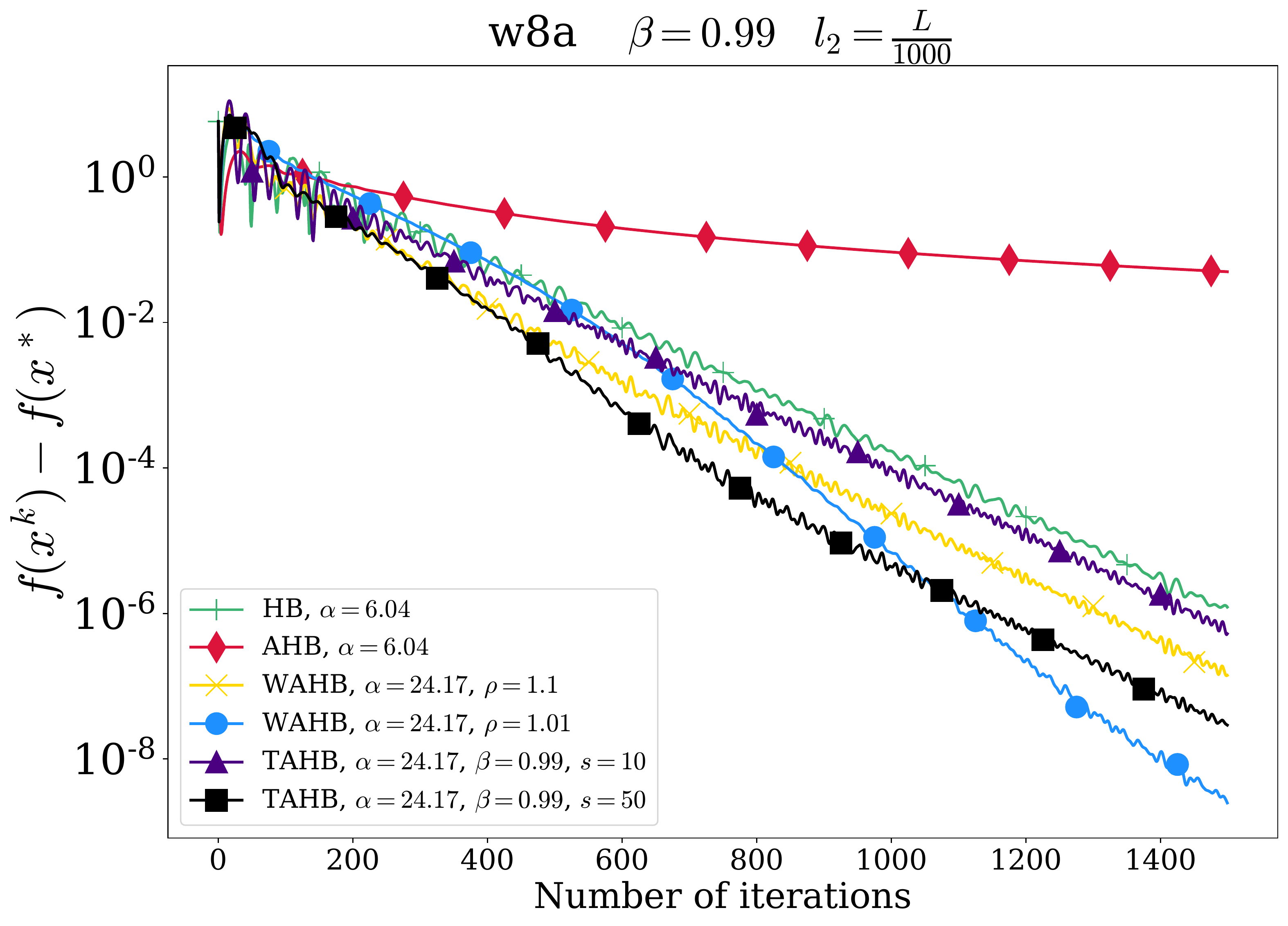}\\
    \includegraphics[width=0.325\textwidth]{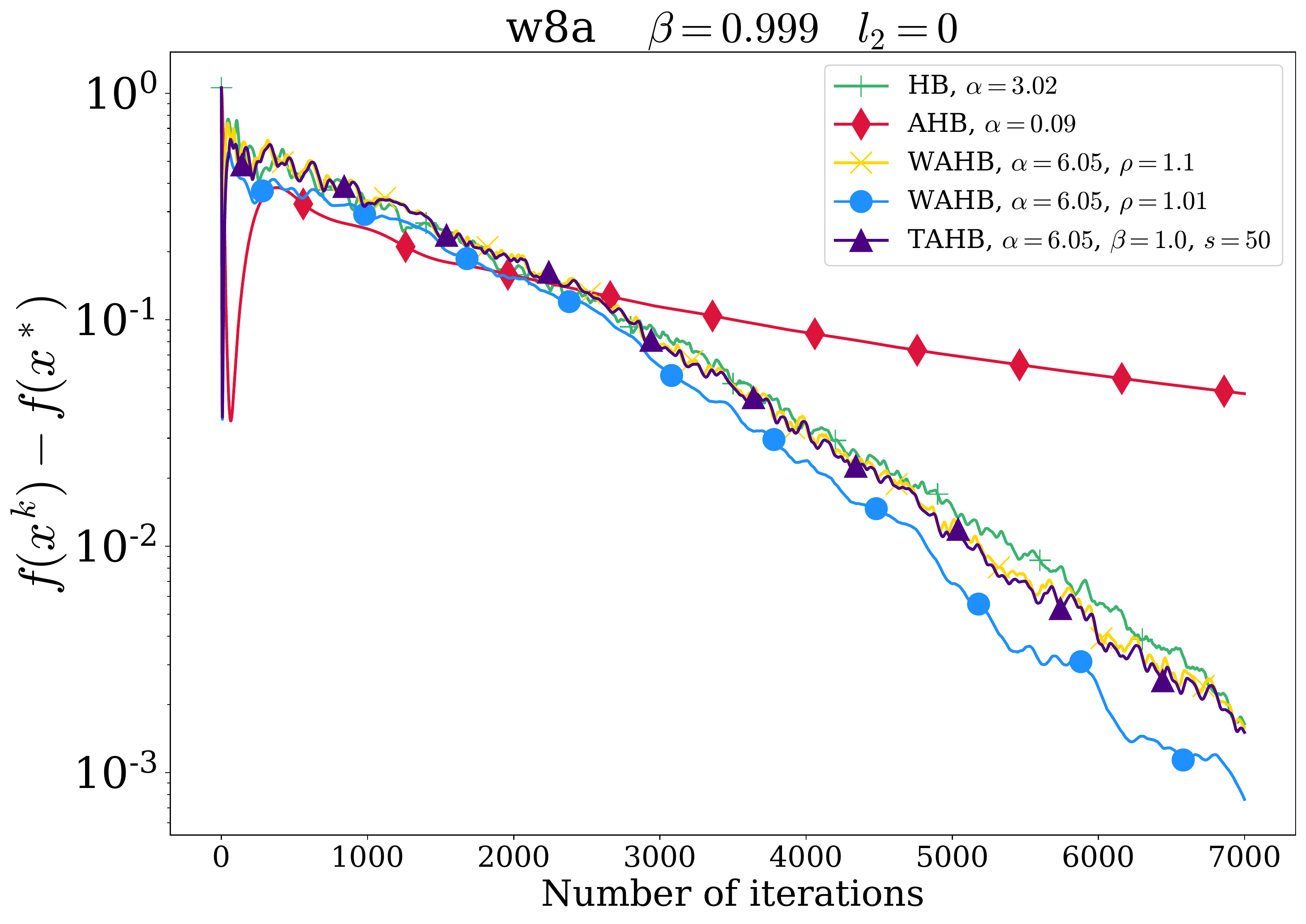}
    \includegraphics[width=0.325\textwidth]{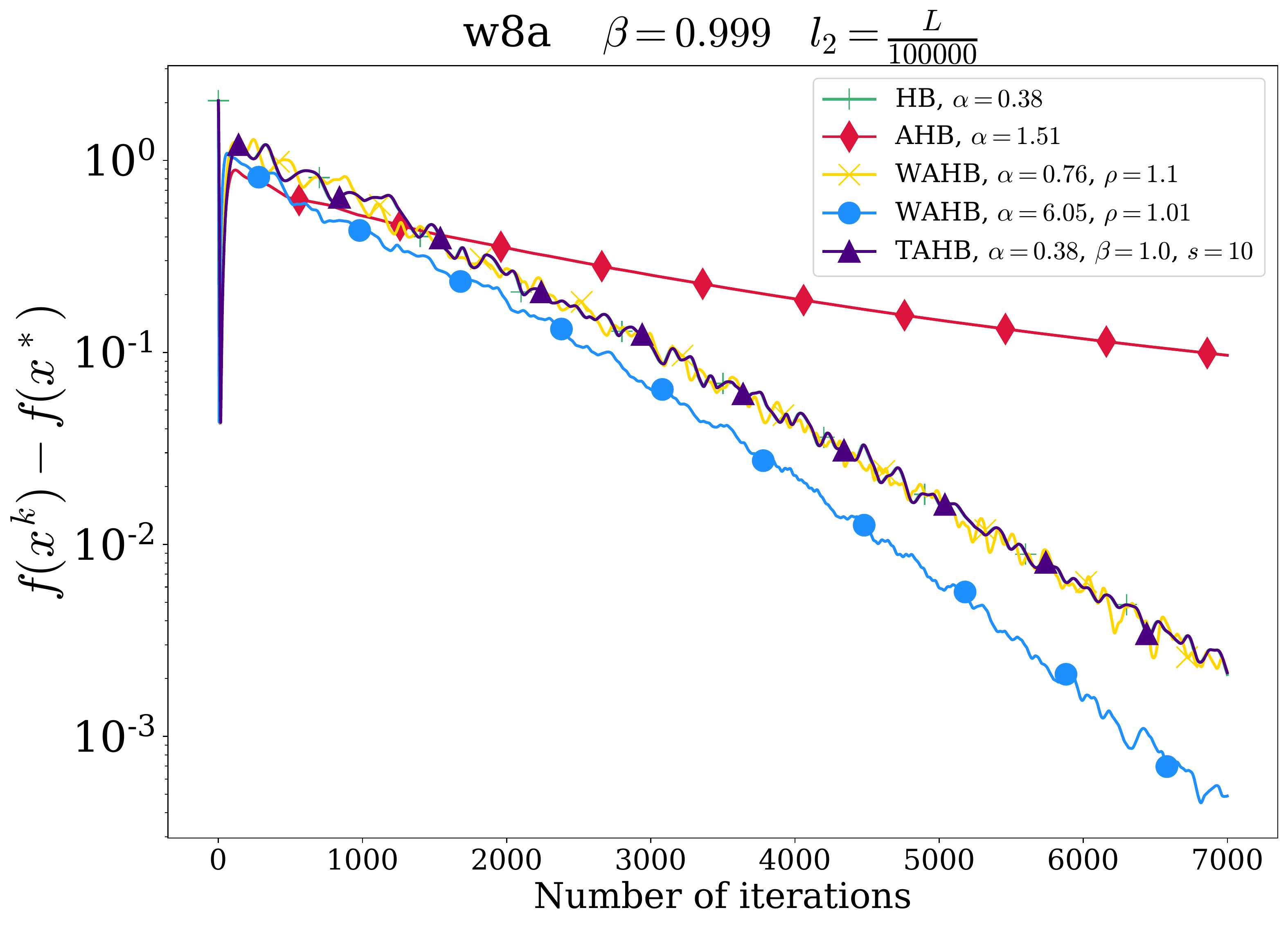}
    \includegraphics[width=0.325\textwidth]{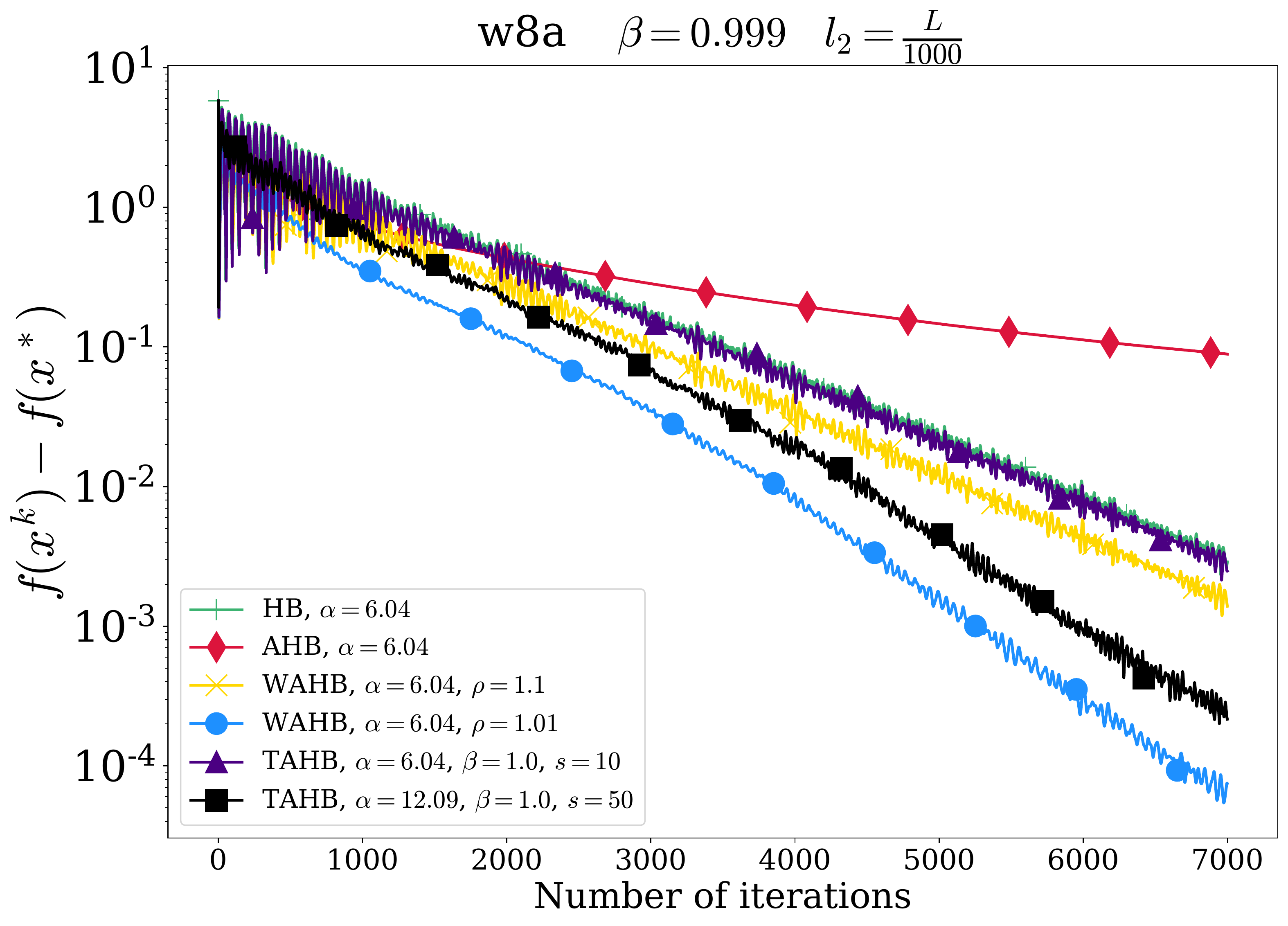}
    \caption{Trajectories of \algname{HB}, \algname{AHB}, \algname{WAHB}, and \algname{TAHB} with different momentum parameters $\beta$ applied to solve logistic regression problem with $\ell_2$-regularization for {\tt phishing} dataset. Stepsize $\alpha$ was tuned for each method and each choice of $\beta$ (and $\rho$, $s$) separately.}
    \label{fig:w8a_different_betas}
\end{figure}

\begin{figure}
    \centering
    \includegraphics[width=0.325\textwidth]{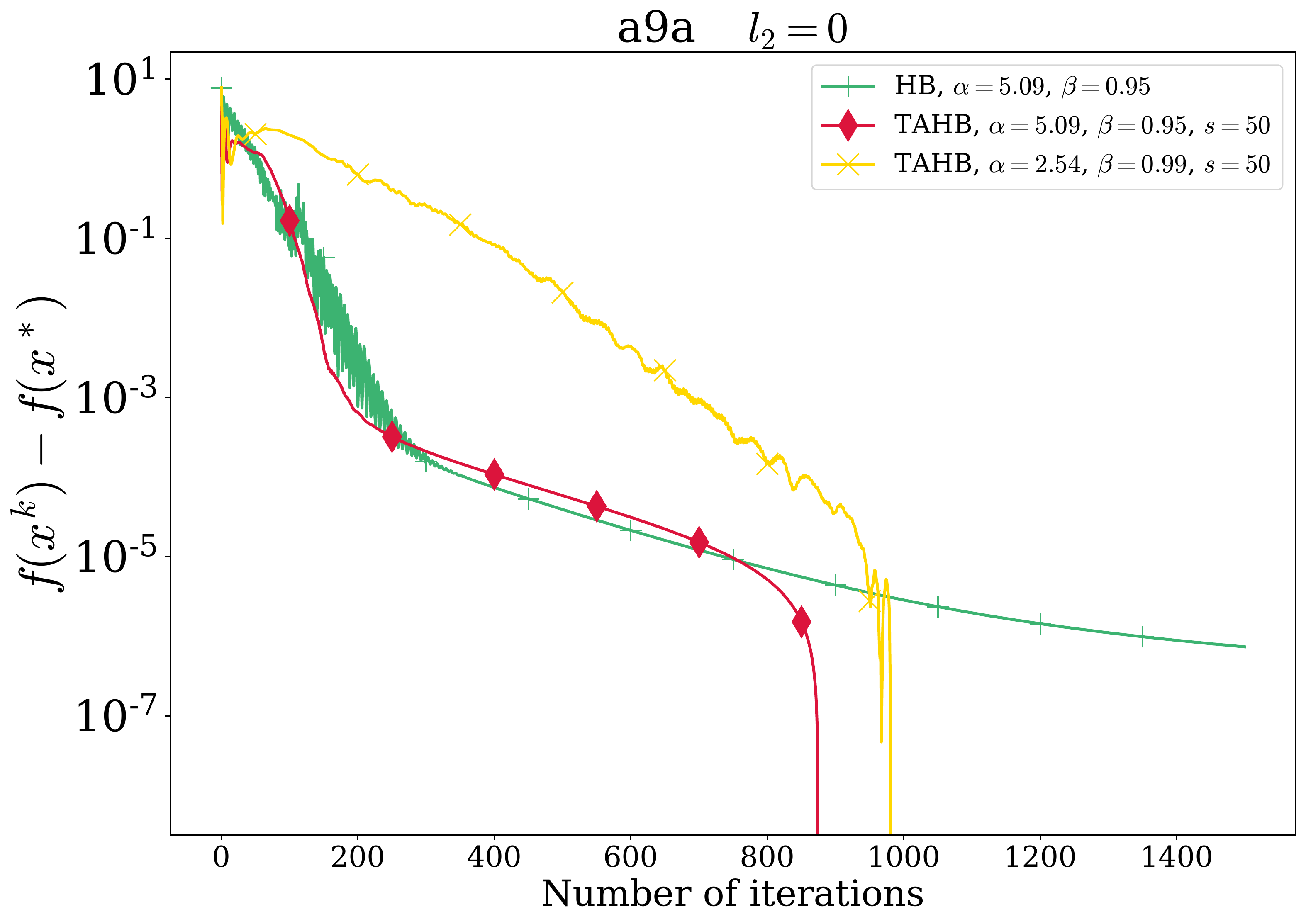}
    \includegraphics[width=0.325\textwidth]{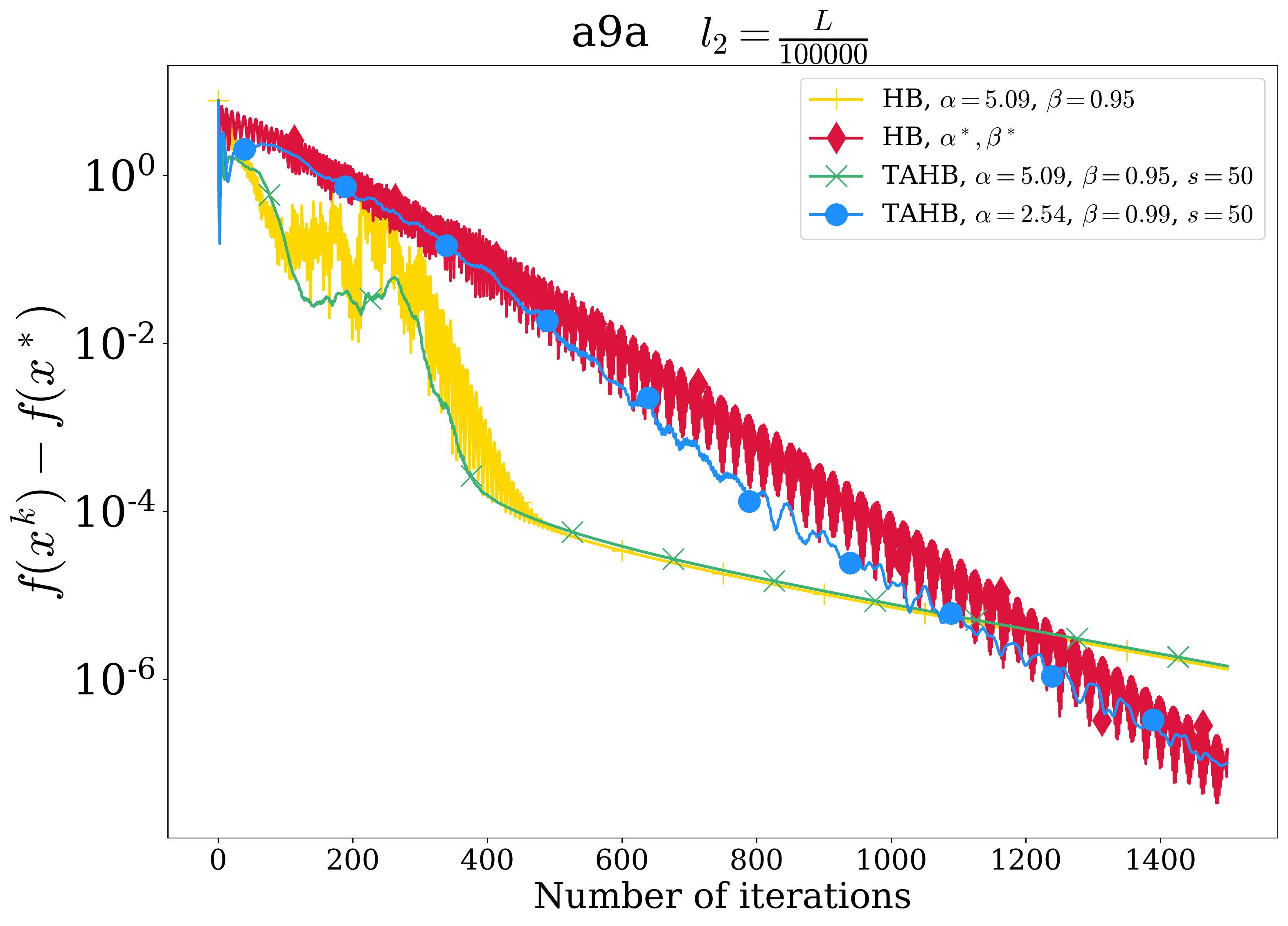}
    \includegraphics[width=0.325\textwidth]{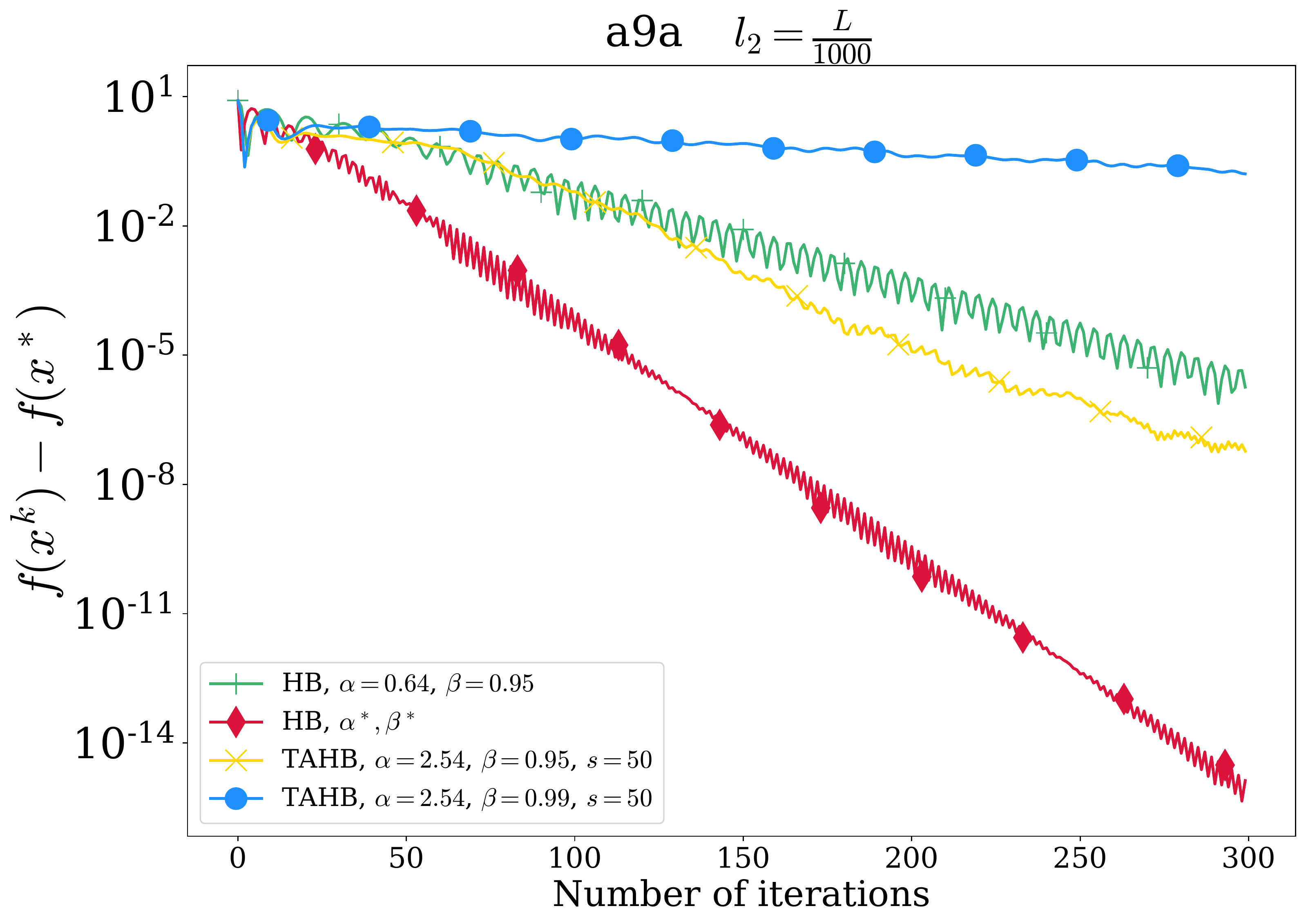}
    \includegraphics[width=0.325\textwidth]{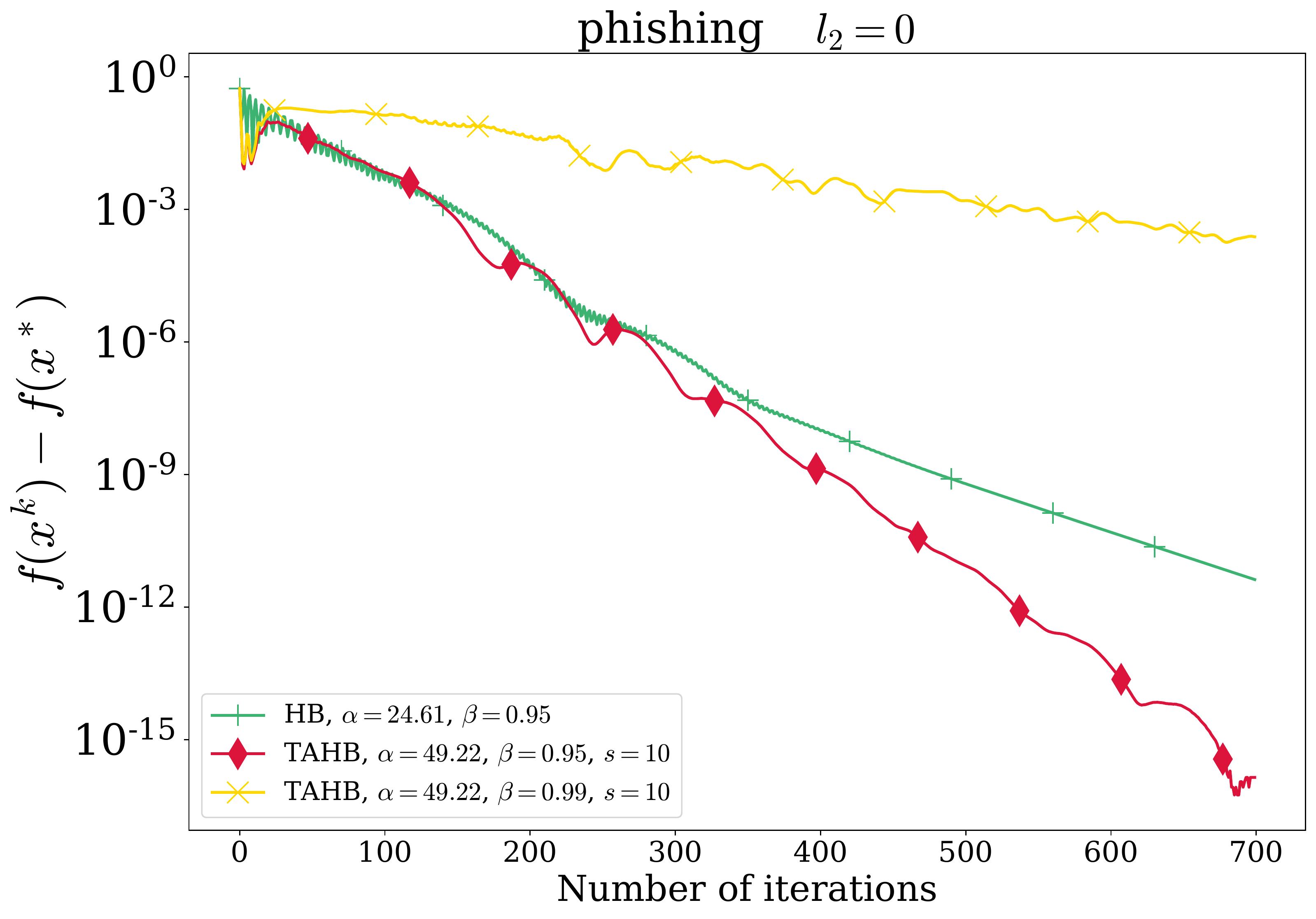}
    \includegraphics[width=0.325\textwidth]{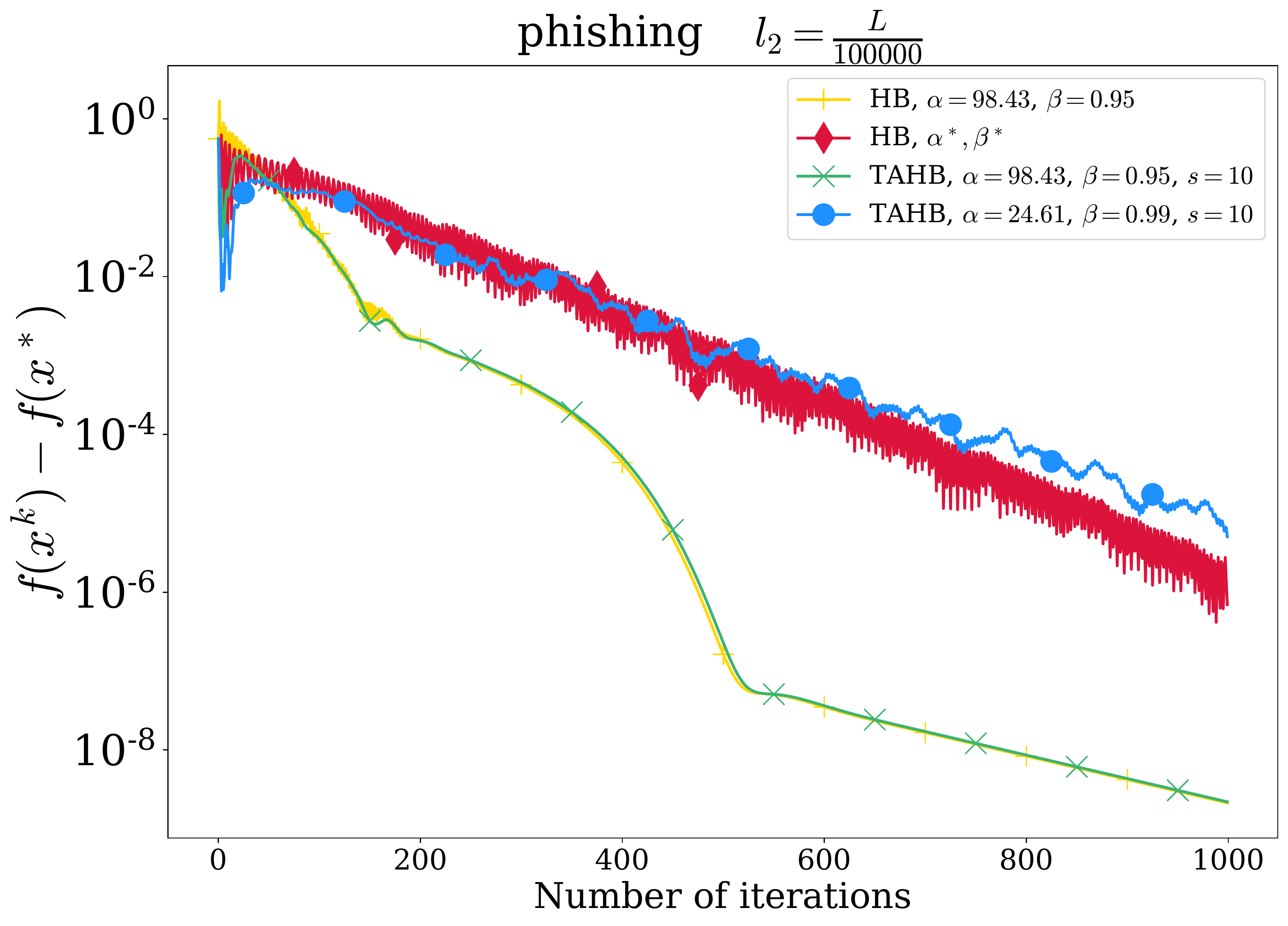}
    \includegraphics[width=0.325\textwidth]{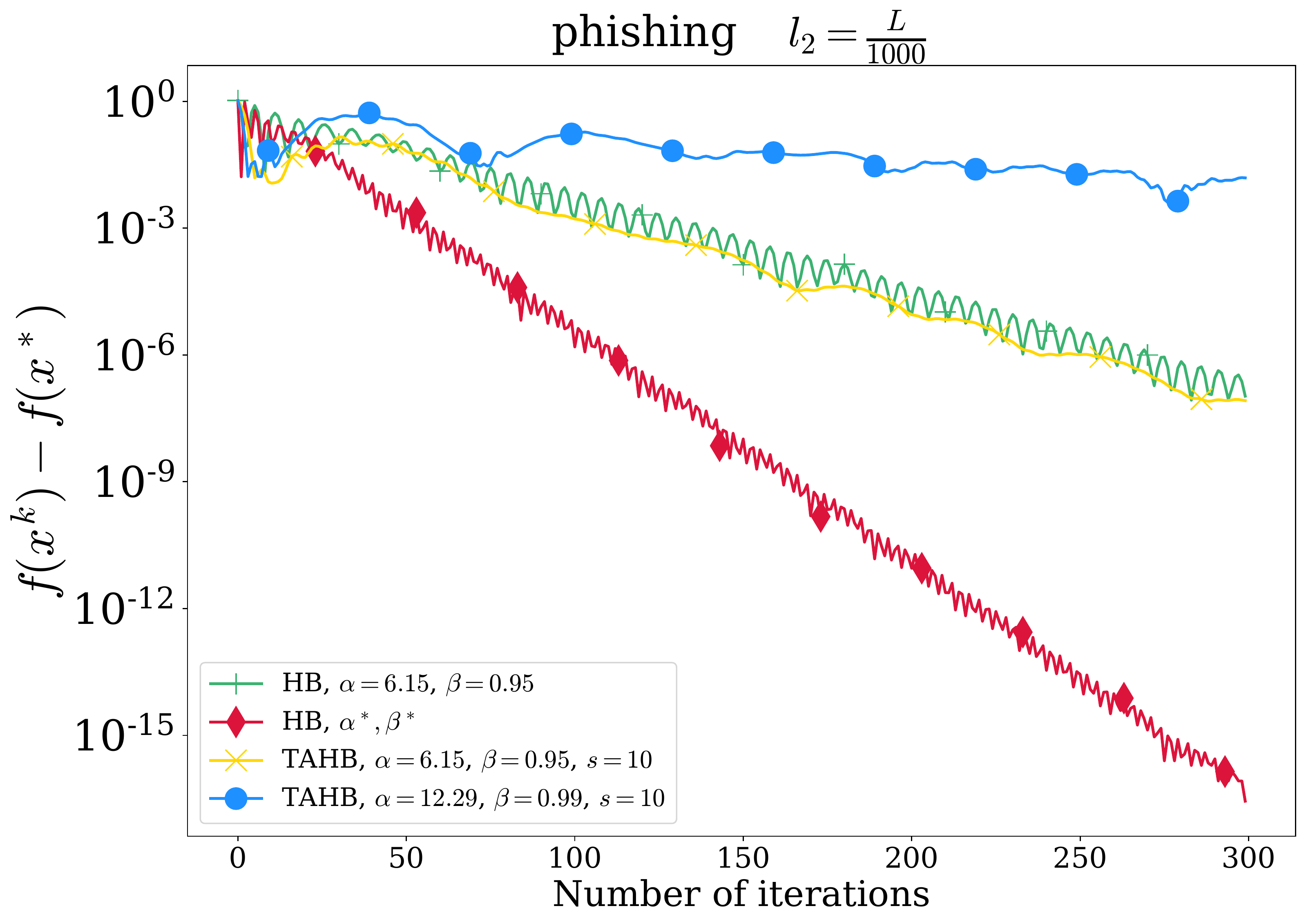}
    \includegraphics[width=0.325\textwidth]{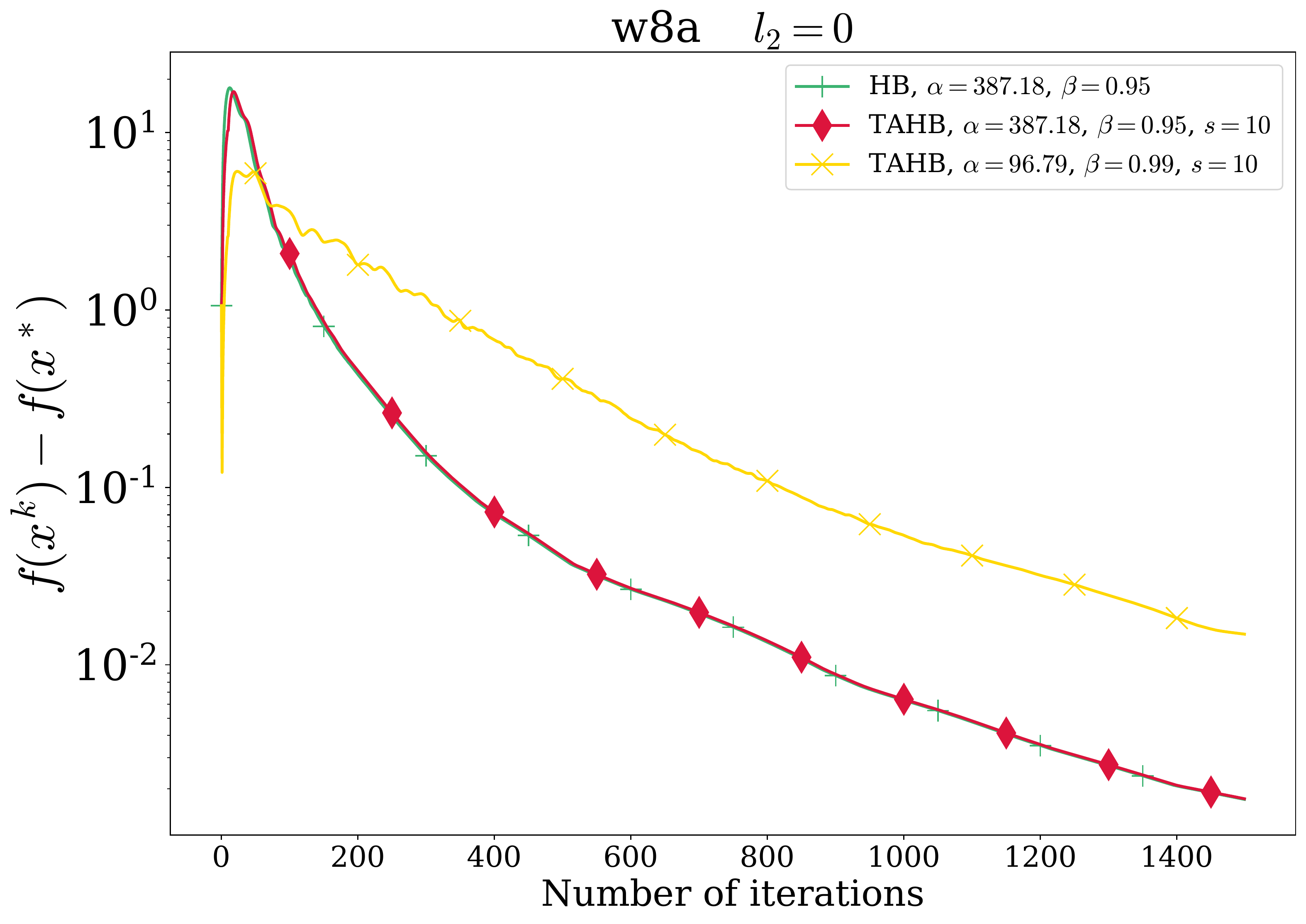}
    \includegraphics[width=0.325\textwidth]{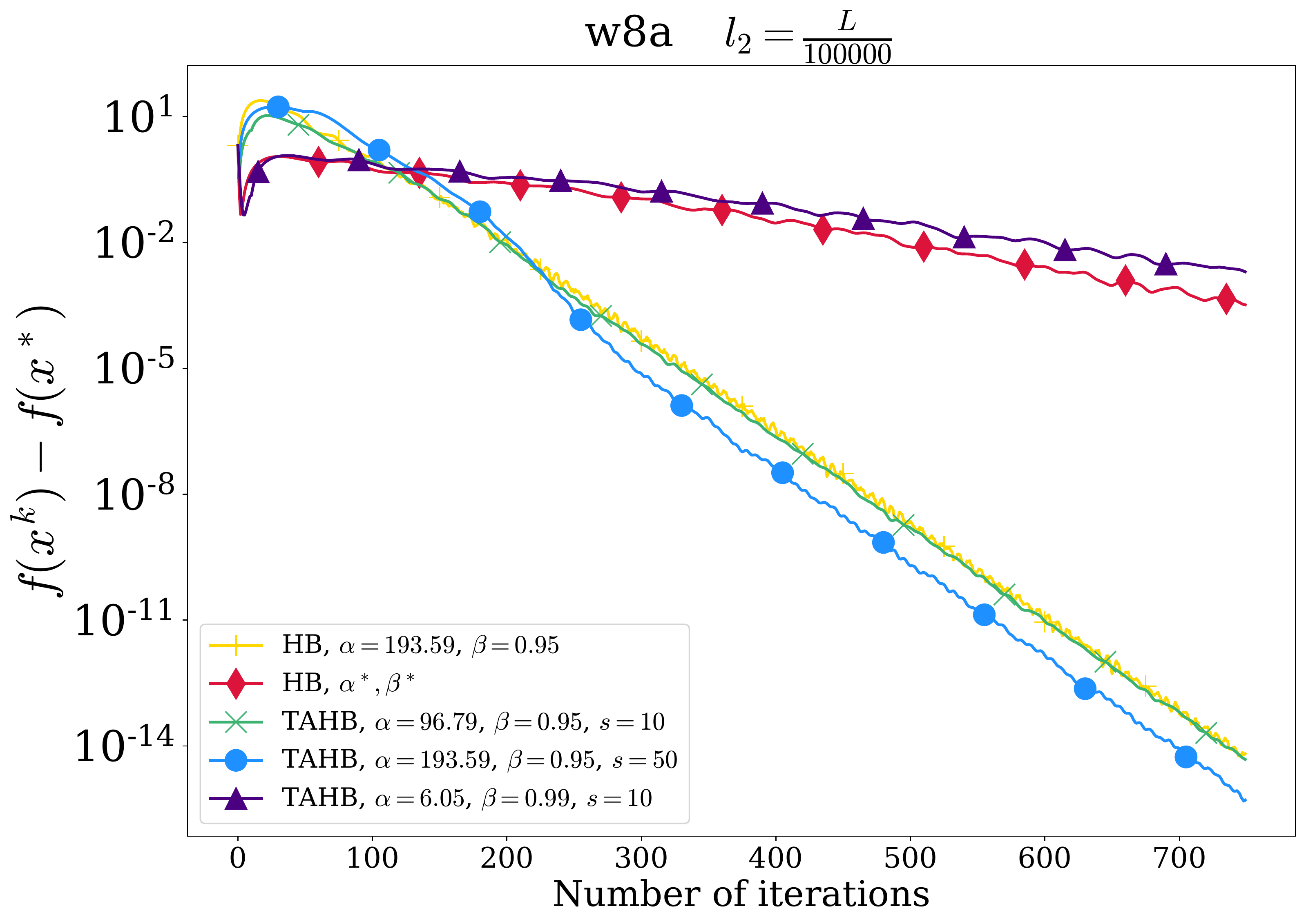}
    \includegraphics[width=0.325\textwidth]{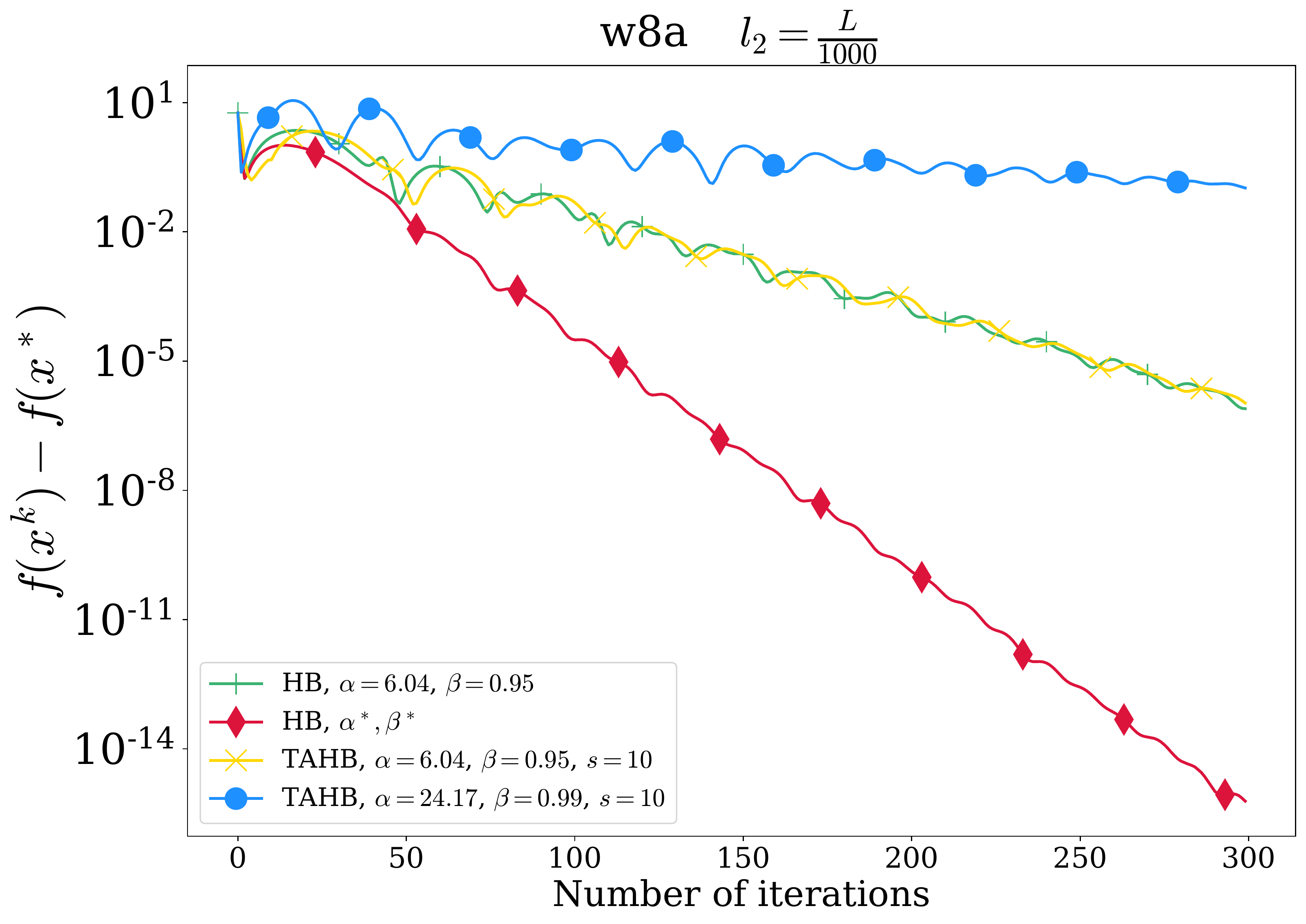}
    \caption{Trajectories of \algname{HB} with $\beta = 0.95$ (standard choice of $\beta$) and \algname{TAHB} with $\beta = 0.95$ and $\beta = 0.99$ (large $\beta$) applied to solve logistic regression problem with $\ell_2$-regularization for dataset from Table~\ref{tab:summary_data}. Stepsize parameter $\alpha$ was tuned for each method separately. For $\ell_2 > 0$ we also show the trajectories of \algname{HB} with optimal parameters $\alpha = \alpha^*$ and $\beta = \beta^*$ from \eqref{eq:optimal_params}.}
    \label{fig:best_params_logreg}
\end{figure}

We run \algname{HB}, \algname{AHB} and \algname{WAHB} with different momentum parameters $\beta$ solve this problem. Moreover, we also tested a modification of \algname{AHB} called Tail-Averaged Heavy-Ball method (\algname{TAHB}, see Algorithm~\ref{alg:TAHB}) with $s \in \{10, 50\}$\footnote{In our experiments, \algname{TAHB} with $s \ge 100$ performed significantly worse than \algname{TAHB} with $s = 50$. Therefore, we report only the resuts for $s \in \{10, 50\}$.}. The weights for \algname{WAHB} were chosen as $w_k = \rho^k$ for $\rho \in \{1.1., 1.01\}$. Next, we chose parameter $\beta$ from the set $\{0.9, 0.95, 0.99, 0.999\}$, and tuned stepsize parameter $\alpha \in \{2^{-4}, 2^{-3}, 2^{-2}, 2^{-1}, 1, 2, 4, 8, 16, 32, 64, 128, 256\}\cdot \nicefrac{1}{L}$ for each method separately for given $\beta$ (and for given $\rho$ in case of \algname{WAHB}, for given $s$ for \algname{TAHB}). The result are shown in Figures~\ref{fig:a9a_different_betas}-\ref{fig:best_params_logreg}.

\paragraph{Figures~\ref{fig:a9a_different_betas}-\ref{fig:w8a_different_betas}.} The plots show that for small $\beta$, i.e., $\beta = 0.9, 0.95$, \algname{HB} does not have significant oscillations and \algname{WAHB} and \algname{TAHB} have comparable performance. However, for larger $\beta$, i.e., $\beta = 0.99, 0.999$, the behavior of \algname{HB} is signigicantly non-monotone and oscillations are quite large. In contrast, \algname{WAHB} and \algname{TAHB} have much smaller oscillations and converge faster than \algname{HB}. These facts illustrate the advantages of using proper averaging scheme for \algname{HB} (either in form of \algname{WAHB} or \algname{TAHB}).

\paragraph{Figure~\ref{fig:best_params_logreg}.} In these plots, we highlight the effect of averaging for large $\beta$. That is, we compare \algname{HB} with standard and commonly used choice of $\beta$ ($\beta = 0.95$) and \algname{TAHB} with $\beta = \{0.95, 0.99\}$. Moreover, for $\ell_2 > 0$ we also tested \algname{HB} with optimal parameters from \eqref{eq:optimal_params}. The results for all considered datasets show that \algname{TAHB} with $\beta = 0.95$ has comparable performance with \algname{HB} and oscillates smaller, while \algname{TAHB} with $\beta = 0.99$ is always slower than \algname{TAHB} with $\beta = 0.95$. Next, when $\ell_2 = \nicefrac{L}{100000}$ (ill-conditioned problems), \algname{TAHB} with $\beta = 0.99$ is as fast as \algname{HB} with optimal parameters but has smaller oscillations. Finally, when $\ell_2 = \nicefrac{L}{1000}$ (well-conditioned problems), \algname{HB} with optimal parameters has negligible oscillations and shows the best performance. Such behavior is natural since for the well-conditioned problems \algname{HB} does not suffer significantly from the non-monotone behavior and peak-effect.

\section{Conclusion}\label{sec:conclusion}
This paper shows the advantages of using averaging for Heavy-Ball method both in theory and practice. That is, our theory and experiments imply that averaging helps to reduce the oscillations of \algname{HB}. Although the derived theoretical convergence guarantees for \algname{HB} with averaging are not better than existing ones for \algname{HB}, in our experiments, we observe that  \algname{HB} with properly adjusted averaging scheme can converge faster than \algname{HB} without averaging. In particular, we observe this phenomenon when momentum parameter $\beta$ for averaged versions of \algname{HB} is chosen to be large enough, e.g., larger than the standard choice of $\beta = 0.95$ and sometimes larger than the optimal choice of $\beta$ from \eqref{eq:optimal_params}.

\section*{Acknowledgments} 
Marina Danilova was supported by Russian Foundation
for Basic Research (Theorems~\ref{thm:AHB_bad_example_deviations}~and~\ref{thm:AHB_deviation_arbitrary_init}, project No. 20-31-90073)
and by Russian Science Foundation (Theorems~\ref{thm:R-AHB}~and~\ref{thm:WAHB_main_result}, project
No. 21-71-30005).

\bibliography{references}

\begin{thebibliography}{10}

\bibitem{chang2011libsvm}
Chih-Chung Chang and Chih-Jen Lin.
\newblock Libsvm: a library for support vector machines.
\newblock {\em ACM transactions on intelligent systems and technology (TIST)},
  2(3):1--27, 2011.

\bibitem{danilova2020recent}
Marina Danilova, Pavel Dvurechensky, Alexander Gasnikov, Eduard Gorbunov,
  Sergey Guminov, Dmitry Kamzolov, and Innokentiy Shibaev.
\newblock Recent theoretical advances in non-convex optimization.
\newblock {\em arXiv preprint arXiv:2012.06188}, 2020.

\bibitem{danilova2018non}
Marina Danilova, Anastasiia Kulakova, and Boris Polyak.
\newblock Non-monotone behavior of the heavy ball method.
\newblock In {\em International Conference on Difference Equations and
  Applications}, pages 213--230. Springer, 2018.

\bibitem{defazio2020understanding}
Aaron Defazio.
\newblock Understanding the role of momentum in non-convex optimization:
  Practical insights from a lyapunov analysis.
\newblock {\em arXiv preprint arXiv:2010.00406}, 2020.

\bibitem{ghadimi2015global}
Euhanna Ghadimi, Hamid~Reza Feyzmahdavian, and Mikael Johansson.
\newblock Global convergence of the heavy-ball method for convex optimization.
\newblock In {\em 2015 European control conference (ECC)}, pages 310--315.
  IEEE, 2015.

\bibitem{gorbunov2020stochastic}
Eduard Gorbunov, Adel Bibi, Ozan Sener, El~Houcine Bergou, and Peter Richtarik.
\newblock A stochastic derivative free optimization method with momentum.
\newblock In {\em International Conference on Learning Representations}, 2020.

\bibitem{lessard2016analysis}
Laurent Lessard, Benjamin Recht, and Andrew Packard.
\newblock Analysis and design of optimization algorithms via integral quadratic
  constraints.
\newblock {\em SIAM Journal on Optimization}, 26(1):57--95, 2016.

\bibitem{mania2017perturbed}
Horia Mania, Xinghao Pan, Dimitris Papailiopoulos, Benjamin Recht, Kannan
  Ramchandran, and Michael~I Jordan.
\newblock Perturbed iterate analysis for asynchronous stochastic optimization.
\newblock {\em SIAM Journal on Optimization}, 27(4):2202--2229, 2017.

\bibitem{mishchenko2019distributed}
Konstantin Mishchenko, Eduard Gorbunov, Martin Tak{\'a}{\v{c}}, and Peter
  Richt{\'a}rik.
\newblock Distributed learning with compressed gradient differences.
\newblock {\em arXiv preprint arXiv:1901.09269}, 2019.

\bibitem{mohammadi2021transient}
Hesameddin Mohammadi, Samantha Samuelson, and Mihailo~R Jovanovi{\'c}.
\newblock Transient growth of accelerated first-order methods for strongly
  convex optimization problems.
\newblock {\em arXiv preprint arXiv:2103.08017}, 2021.

\bibitem{nemirovsky1983problem}
A.S. Nemirovsky and D.B. Yudin.
\newblock {\em Problem Complexity and Method Efficiency in Optimization}.
\newblock J. Wiley \& Sons, New York, 1983.

\bibitem{nesterov1983method}
Yurii Nesterov.
\newblock A method for unconstrained convex minimization problem with the rate
  of convergence {O}$(1/k^2)$.
\newblock In {\em Doklady an ussr}, volume 269, pages 543--547, 1983.

\bibitem{nesterov2018lectures}
Yurii Nesterov.
\newblock {\em Lectures on convex optimization}, volume 137.
\newblock Springer, 2018.

\bibitem{polyak1987introduction}
Boris Polyak.
\newblock {\em Introduction to Optimization}.
\newblock New York, Optimization Software, 1987.

\bibitem{polyak1964some}
Boris~T Polyak.
\newblock Some methods of speeding up the convergence of iteration methods.
\newblock {\em USSR Computational Mathematics and Mathematical Physics},
  4(5):1--17, 1964.

\bibitem{taylor2019stochastic}
Adrien Taylor and Francis Bach.
\newblock Stochastic first-order methods: non-asymptotic and computer-aided
  analyses via potential functions.
\newblock In {\em Conference on Learning Theory}, pages 2934--2992. PMLR, 2019.

\bibitem{taylor2018lyapunov}
Adrien Taylor, Bryan Van~Scoy, and Laurent Lessard.
\newblock Lyapunov functions for first-order methods: Tight automated
  convergence guarantees.
\newblock In {\em International Conference on Machine Learning}, pages
  4897--4906. PMLR, 2018.

\bibitem{taylor2017performance}
Adrien~B Taylor, Julien~M Hendrickx, and Fran{\c{c}}ois Glineur.
\newblock Performance estimation toolbox (pesto): automated worst-case analysis
  of first-order optimization methods.
\newblock In {\em 2017 IEEE 56th Annual Conference on Decision and Control
  (CDC)}, pages 1278--1283. IEEE, 2017.

\bibitem{yang2016unified}
Tianbao Yang, Qihang Lin, and Zhe Li.
\newblock Unified convergence analysis of stochastic momentum methods for
  convex and non-convex optimization.
\newblock {\em arXiv preprint arXiv:1604.03257}, 2016.

\bibitem{yu2019linear}
Hao Yu, Rong Jin, and Sen Yang.
\newblock On the linear speedup analysis of communication efficient momentum
  sgd for distributed non-convex optimization.
\newblock In {\em International Conference on Machine Learning}, pages
  7184--7193. PMLR, 2019.

\end{thebibliography}

\appendix

\section{Basic Inequalities}\label{sec:basic_facts}
For all $a,b\in\R^n$ and $\lambda > 0, q\in (0,1]$
\begin{equation}
    |\la a, b\ra| \le \frac{\|a\|_2^2}{2\lambda} + \frac{\lambda\|b\|_2^2}{2},\label{eq:fenchel_young_inequality}
\end{equation}
\begin{equation}
    \|a+b\|_2^2 \le 2\|a\|_2^2 + 2\|b\|_2^2,\label{eq:squared_norm_sum}
\end{equation}
\begin{equation}
    \|a+b\|^2 \le (1+\lambda)\|a\|^2 + \left(1 + \frac{1}{\lambda}\right)\|b\|^2, \label{eq:1+lambda}
\end{equation}
\begin{equation}
    \la a, b\ra = \frac{1}{2}\left(\|a+b\|_2^2 - \|a\|_2^2 - \|b\|_2^2\right), \label{eq:inner_product_representation}
\end{equation}
\begin{equation}
    \left(1 - \frac{q}{2}\right)^{-1} \le 1+q, \label{eq:1-q/2}
\end{equation}
\begin{equation}
    \left(1 + \frac{q}{2}\right)(1-q) \le 1-\frac{q}{2}.  \label{eq:1+q/2}
\end{equation}

\section{Auxiliary Results}
\begin{lemma}{Lemma 1 from \cite{mohammadi2021transient}}\label{lem:aux_lemma_1}
    Let $\rho_1$ and $\rho_2$ be the eigenvalues of the matrix $\mM = \begin{bmatrix}
         a & b\\
         1 & 0
    \end{bmatrix}$ and let $k$ be a positive integer. If $\rho_1 \neq \rho_2$, then we have
    \begin{equation*}
        \mM^k = \frac{1}{\rho_2 - \rho_1} \begin{bmatrix}
             \rho_2^{k+1} - \rho_1^{k+1} & \rho_1\rho_2(\rho_1^k - \rho_2^k)\\
             \rho_2^k - \rho_1^k & \rho_1\rho_2(\rho_1^{k-1} - \rho_2^{k-1})
        \end{bmatrix}.
    \end{equation*}
    Moreover, if $\rho_1 = \rho_2 = \rho$, the matrix $\mM^k$ satisfies
    \begin{equation*}
        \mM^k = \begin{bmatrix}
            (k+1)\rho^k & -k\rho^{k+1}\\
            k\rho^{k-1} & (1-k)\rho^k
        \end{bmatrix}.
    \end{equation*}
\end{lemma}

\section{Missing Proofs from Section~\ref{sec:quadratic_case}}

In this section, for $x$, we use the upper index for an iteration counter, and the lower index denotes the component of the vector.

\subsection{Proof of Theorem~\ref{thm:AHB_bad_example_deviations}}

Rewriting the update rule of \algname{HB} for $f(x) = \frac12 x^{\top} \mA x$ with $ \mA = \mathrm{diag}\left(\mu, \ \lambda_2, \ \ldots, \ \lambda_{n-1}, L\right)$ with $\alpha = \frac{1}{L}$ we get
    \begin{eqnarray*}
        x_1^{k+1} &=& \left(1- \frac{\mu}{L} + \beta\right)x_1^k - \beta x_1^{k-1},\\
        x_2^{k+1} &=& \left(1- \frac{\lambda_2}{L} + \beta\right)x_2^k - \beta x_2^{k-1},\\
        &\vdots&\\
        x_{n-1}^{k+1} &=& \left(1- \frac{\lambda_{n-1}}{L} + \beta\right)x_{n-1}^k - \beta x_{n-1}^{k-1},\\
        x_n^{k+1} &=& \beta x_n^k - \beta x_n^{k-1}.
    \end{eqnarray*}
    To solve these recurrences we consider the corresponding characteristic equations:
    \begin{eqnarray*}
        \rho^2 &=& \left(1- \frac{\mu}{L} + \beta\right)\rho - \beta,\\
        \rho^2 &=& \left(1- \frac{\lambda_2}{L} + \beta\right)\rho - \beta,\\
        &\vdots&\\
        \rho^2 &=& \left(1- \frac{\lambda_{n-1}}{L} + \beta\right)\rho - \beta,\\
        \rho^2 &=& \beta \rho - \beta.
    \end{eqnarray*}
    
    Since $\beta \le (1-2\sqrt{\nicefrac{\mu}{L}})^2 < (1-\sqrt{\nicefrac{\mu}{L}})^2$ the roots of the first equation are
    \begin{eqnarray*}
        \rho_1(\mu) &=& \frac{1+\beta - \nicefrac{\mu}{L} + \sqrt{\left(1+\beta - \nicefrac{\mu}{L}\right)^2 - 4\beta}}{2},\\
        \rho_2(\mu) &=& \frac{1+\beta - \nicefrac{\mu}{L} - \sqrt{\left(1+\beta - \nicefrac{\mu}{L}\right)^2 - 4\beta}}{2}.
    \end{eqnarray*}
    Moreover, we have $\sqrt{\left(1+\beta - \nicefrac{\mu}{L}\right)^2 - 4\beta} \le 1 - \beta + \nicefrac{\mu}{L}$, and, as a consequence, $0 < \rho_2(\mu) < \rho_1(\mu) < 1$. Next, the first components of iterates produced by \algname{HB} satisfy
    \begin{equation}
        x_1^k = C_1\rho_1^k(\mu) + C_2\rho_2^k(\mu) \notag
    \end{equation}
    with some constants $C_1,C_2 \in \R$. This equation and the choice of the starting points $x^0=x^1 = (1,1,\ldots,1)^{\top}$ imply
    \begin{equation*}
        \begin{cases}
            C_1 + C_2 &= 1,\\
            C_1\rho_1(\mu) + C_2\rho_2(\mu) &= 1,
        \end{cases}
    \end{equation*}
    whence
    \begin{equation*}
        C_1 = \frac{1 - \rho_2(\mu)}{\rho_1(\mu) - \rho_2(\mu)},\quad C_2 = 1 - C_1 = \frac{\rho_1(\mu) - 1}{\rho_1(\mu) - \rho_2(\mu)}. 
    \end{equation*}
    Using the formula for $C_1$ and $\beta \in [(1-3\sqrt{\nicefrac{\mu}{L}})^2,  (1-2\sqrt{\nicefrac{\mu}{L}})^2]$ we derive that $C_1 > 0$ and
    \begin{eqnarray*}
        C_1 &=& \left(1 - \frac{1+\beta - \nicefrac{\mu}{L} - \sqrt{\left(1+\beta - \nicefrac{\mu}{L}\right)^2 - 4\beta}}{2}\right)\frac{1}{\sqrt{\left(1+\beta - \nicefrac{\mu}{L}\right)^2 - 4\beta}}\\
        &=& \frac{1 - \beta + \nicefrac{\mu}{L} + \sqrt{\left(1+\beta - \nicefrac{\mu}{L}\right)^2 - 4\beta}}{2\sqrt{\left(1+\beta - \nicefrac{\mu}{L}\right)^2 - 4\beta}}\\
        &=& \frac{1}{2} + \frac{1 - \beta + \nicefrac{\mu}{L}}{2\sqrt{\left(1+\beta - \nicefrac{\mu}{L}\right)^2 - 4\beta}}\\
        &\le& \frac{1}{2} + \frac{1 - (1-3\sqrt{\nicefrac{\mu}{L}})^2 + \nicefrac{\mu}{L}}{2\sqrt{\left(1+(1-2\sqrt{\nicefrac{\mu}{L}})^2 - \nicefrac{\mu}{L}\right)^2 - 4(1-2\sqrt{\nicefrac{\mu}{L}})^2}}\\
        &=& \frac{1}{2} + \frac{3\sqrt{\nicefrac{\mu}{L}} - 4\nicefrac{\mu}{L}}{\sqrt{\left(2 - 4\sqrt{\nicefrac{\mu}{L}} + 3\nicefrac{\mu}{L}\right)^2 - (2-4\sqrt{\nicefrac{\mu}{L}})^2}}\\
        &=& \frac{1}{2} + \frac{3\sqrt{\nicefrac{\mu}{L}} - 4\nicefrac{\mu}{L}}{\sqrt{3\nicefrac{\mu}{L}\left(4 - 8\sqrt{\nicefrac{\mu}{L}} + 3\nicefrac{\mu}{L}\right)}}.
    \end{eqnarray*}
    Since $L \ge 100\mu$ we can further upper bound the right-hand side of the previous inequality and get
    \begin{eqnarray*}
        C_1 &\le& \frac{1}{2} + \frac{3\sqrt{\nicefrac{\mu}{L}}}{\sqrt{3\nicefrac{\mu}{L}\left(4 - 8\sqrt{\nicefrac{\mu}{L}} \right)}} \le \frac{1}{2} + \frac{\sqrt{3}}{\sqrt{4 - \nicefrac{4}{5}}} = \frac{1}{2} + \frac{1}{2}\sqrt{\frac{15}{4}} \le \frac{3}{2}.
    \end{eqnarray*}
    Taking into account that $C_1 > 0$ and $C_2 = 1 - C_1$ we derive that $|C_2| = \max\{1-C_1, C_1 - 1\} \le \nicefrac{1}{2}$. Putting all together, we obtain
    \begin{equation*}
        |x_1^k| = |C_1\rho_1^k(\mu) + C_2\rho_2^k(\mu)| \le |C_1| + |C_2| \le 2 \quad \forall k\ge 0
    \end{equation*}
    
    In the remaining part of the proof, we handle the characteristic equations
    \begin{eqnarray*}
        \rho^2 &=& \left(1- \frac{\lambda_2}{L} + \beta\right)\rho - \beta,\\
        &\vdots&\\
        \rho^2 &=& \left(1- \frac{\lambda_{n-1}}{L} + \beta\right)\rho - \beta,\\
        \rho^2 &=& \beta \rho - \beta.
    \end{eqnarray*}
    Without loss of generality, we consider the equation
    \begin{equation}
        \rho^2 = \left(1- \frac{\lambda}{L} + \beta\right)\rho - \beta\label{eq:technical_char_eq}
    \end{equation}
    with $\lambda \in [\lambda_2, L]$. This equation serves as a characteristic equation for the sequence $\{y_k\}_{k\ge 0}\subseteq \R$ satisfying
    \begin{equation*}
        y_{k+1} = \left(1- \frac{\lambda}{L} + \beta\right)y_k - \beta y_{k-1}.
    \end{equation*}
    Since $\lambda \ge \lambda_2 \ge 10\mu$ and $\beta \ge (1-3\sqrt{\nicefrac{\mu}{L}})^2$ we conclude that $\beta \ge (1-\sqrt{\nicefrac{\lambda}{L}})^2$ and the characteristic equation has the complex roots with non-zero imaginary parts:
    \begin{eqnarray*}
        \rho_1(\lambda) &=& \frac{1+\beta - \nicefrac{\lambda}{L} + i\sqrt{4\beta - \left(1+\beta - \nicefrac{\mu}{L}\right)^2}}{2},\\
        \rho_2(\lambda) &=& \frac{1+\beta - \nicefrac{\mu}{L} - i\sqrt{4\beta - \left(1+\beta - \nicefrac{\mu}{L}\right)^2}}{2}.
    \end{eqnarray*}
    This implies that $|\rho_1(\lambda)| = |\rho_1(\lambda)| = \sqrt{\beta} < 1$ and
    \begin{equation*}
        y_k = C_1\rho_1^k(\lambda) + C_2\rho_2^k(\lambda)
    \end{equation*}
    for some complex numbers $C_1, C_2$. Let $y_0 = y_1 = 1$. Then,
    \begin{equation*}
        \begin{cases}
            C_1 + C_2 &= 1,\\
            C_1\rho_1(\lambda) + C_2\rho_2(\lambda) &= 1,
        \end{cases}
    \end{equation*}
    whence
    \begin{equation*}
        C_1 = \frac{1 - \rho_2(\lambda)}{\rho_1(\lambda) - \rho_2(\lambda)},\quad C_2 = 1 - C_1 = \frac{\rho_1(\lambda) - 1}{\rho_1(\lambda) - \rho_2(\lambda)}. 
    \end{equation*}
    Using the formula for $C_1$ and $\beta \in [(1-3\sqrt{\nicefrac{\mu}{L}})^2,  (1-2\sqrt{\nicefrac{\mu}{L}})^2]$ we derive that
    \begin{eqnarray*}
        C_1 &=& \left(1 - \frac{1+\beta - \nicefrac{\lambda}{L} - i\sqrt{4\beta - \left(1+\beta - \nicefrac{\lambda}{L}\right)^2}}{2}\right)\frac{1}{i\sqrt{4\beta - \left(1+\beta - \nicefrac{\lambda}{L}\right)^2}}\\
        &=& \frac{1 - \beta + \nicefrac{\lambda}{L} + i\sqrt{4\beta - \left(1+\beta - \nicefrac{\lambda}{L}\right)^2}}{2i\sqrt{4\beta - \left(1+\beta - \nicefrac{\lambda}{L}\right)^2}}\\
        &=& \frac{1}{2} - i\frac{1 - \beta + \nicefrac{\lambda}{L}}{2\sqrt{4\beta - \left(1+\beta - \nicefrac{\lambda}{L}\right)^2}}.
    \end{eqnarray*}
    Then, for the absolute value of $C_1$ we have
    \begin{eqnarray*}
        |C_1| &=& \frac{1}{2}\sqrt{1 + \frac{(1+\nicefrac{\lambda}{L} - \beta)^2}{4\beta - (1-\nicefrac{\lambda}{L} + \beta)^2}}\\
        &\le& \frac{1}{2}\sqrt{1 + \frac{\left(1+\nicefrac{\lambda}{L} - \left(1-3\sqrt{\nicefrac{\mu}{L}}\right)^2\right)^2}{4\left(1-2\sqrt{\nicefrac{\mu}{L}}\right)^2 - \left(1-\nicefrac{\lambda}{L} + \left(1-2\sqrt{\nicefrac{\mu}{L}}\right)^2\right)^2}}\\
        &=& \frac{1}{2}\sqrt{1 + \frac{\left(\nicefrac{\lambda}{L} - 9\nicefrac{\mu}{L} + 6\sqrt{\nicefrac{\mu}{L}} \right)^2}{\left(2-4\sqrt{\nicefrac{\mu}{L}}\right)^2 - \left(2 - \nicefrac{\lambda}{L} - 4 \sqrt{\nicefrac{\mu}{L}} + 4\nicefrac{\mu}{L}\right)^2}}\\
        &=& \frac{1}{2}\sqrt{1 + \frac{\left(\nicefrac{\lambda}{L} - 9\nicefrac{\mu}{L} + 6\sqrt{\nicefrac{\mu}{L}} \right)^2}{\left(4- \nicefrac{\lambda}{L} - 8\sqrt{\nicefrac{\mu}{L}} + 4\nicefrac{\mu}{L}\right)\left(\nicefrac{\lambda}{L} - 4\nicefrac{\mu}{L}\right)}}\\
        &\le& \frac{1}{2}\sqrt{1 + \frac{\left(\nicefrac{\lambda}{L} + 6\sqrt{\nicefrac{\mu}{L}} \right)^2}{\left(3 - \nicefrac{4}{5}\right)\left(\nicefrac{\lambda}{L} - \nicefrac{2\lambda}{5L}\right)}} = \frac{1}{2}\sqrt{1 + \frac{25}{33}\frac{\nicefrac{\lambda^2}{L^2}+ 12\nicefrac{(\lambda\sqrt{\mu})}{(L\sqrt{L})} + 36\nicefrac{\mu}{L}}{\nicefrac{\lambda}{L}}}\\
        &=& \frac{1}{2}\sqrt{1 + \frac{25}{33}\left(\frac{\lambda}{L} + 12\sqrt{\frac{\mu}{L}} + 36\frac{\mu}{\lambda}\right)} \le \frac{1}{2}\sqrt{1 + \frac{25}{33}\left(1 + \frac{6}{5} + \frac{8}{25}\right)} \le 1.
    \end{eqnarray*}
    Since
    \begin{equation*}
        C_2 = 1 - C_1 = \frac{1}{2} + i\frac{1 - \beta + \nicefrac{\lambda}{L}}{2\sqrt{4\beta - \left(1+\beta - \nicefrac{\lambda}{L}\right)^2}}
    \end{equation*}
    we also have $|C_2| = |C_1| \le 1$, and, as a consequence,
    \begin{equation*}
        |y_k| = |C_1\rho_1^k(\lambda) + C_2\rho_2^k(\lambda)| \le |C_1| + |C_2| \le 2\quad \forall k\ge 0.
    \end{equation*}
    This result implies that $|x_i^k| \le 2$ for all $k\ge 0$ and $i = 2,\ldots, n$.
    
    Finally, since $\overline{x}^k = \frac{1}{k+1}\sum_{t=0}^k x^t$ we conclude that
    \begin{equation*}
        |\overline{x}_i^k| \le \frac{1}{k+1}\sum\limits_{t=0}^k|x_i^t| \le 2 \quad \forall k \ge 0, \; i = 1,\ldots,n
    \end{equation*}
    that is equivalent to \eqref{eq:AHB_peak}.

\subsection{Proof of Theorem~\ref{thm:AHB_deviation_arbitrary_init}}

To estimate $\text{dev}_{\algname{HB}}(\alpha,\beta)$ we consider the spectral decomposition of matrix $\mA = \mU \mLambda \mU^\top \succ 0$, where $\mLambda = \mathrm{diag}(\lambda_1,\ldots,\lambda_n)$ is a diagonal matrix of the eigenvalues of $\mA$, $\lambda_1 \le \ldots \le \lambda_n$, and $U$ is a unitary matrix of the eigenvectors of $\mA$. Next, without loss of the generality we assume that $x^* = 0$. Applying the unitary transformation $\mU^\top$ to $x^k$ we obtain $\hat x^k = \mU^\top x^k$ and
\begin{equation*}
    \hat z^k := \begin{bmatrix}
         \hat x^{k+1}\\
         \hat x^k
    \end{bmatrix} = \hat \mT \begin{bmatrix}
         \hat x^{k}\\
         \hat x^{k-1}
    \end{bmatrix} =\ldots = \hat\mT^k \begin{bmatrix}
         \hat x^{1}\\
         \hat x^0
    \end{bmatrix},
\end{equation*}
where
\begin{equation}
    \hat\mT = \left[\begin{array}{@{}c|c@{}}
     \begin{matrix}
     (1 + \beta) \mI - \alpha \mLambda
     \end{matrix}
     & -\beta \mI \\
    \hline
     \mI &
     {\bf 0}
    \end{array}\right] = \left[\begin{array}{@{}c|c@{}}
     \mU^\top & {\bf 0} \\
    \hline
    {\bf 0} & \mU^\top 
    \end{array}\right]\mT. \notag
\end{equation}
In particular, these formulas imply
\begin{equation*}
    \begin{bmatrix}
         \hat x_j^{k+1}\\
         \hat x_j^k
    \end{bmatrix} =  \hat \mT_j \begin{bmatrix}
         \hat x_j^{k}\\
         \hat x_j^{k-1}
    \end{bmatrix} =\ldots = \hat\mT_j^k \begin{bmatrix}
         \hat x_j^{1}\\
         \hat x_j^0
    \end{bmatrix},
\end{equation*}
where
\begin{equation}
    \hat\mT_j = \begin{bmatrix}
         1+\beta - \alpha\lambda_j & -\beta\\
         1 & 0
    \end{bmatrix} \notag
\end{equation}
for all $j = 1,\ldots, n$. Moreover, $\|\mC\mT^k\|_2 = \max\limits_{j=1,\ldots,n}\|\mC_j\hat\mT_j^k\|_2$, where $\mC_j = \begin{bmatrix}0& 1\end{bmatrix}$.

It is easy to see that the eigenvalues of $\hat \mT_j$ are
\begin{equation*}
    \rho_{j,1} = \frac{1+\beta - \nicefrac{\lambda_j}{\lambda_n}+\sqrt{(1+\beta - \nicefrac{\lambda_j}{\lambda_n})^2 - 4\beta}}{2},\quad \rho_{j,2} = \frac{1+\beta - \nicefrac{\lambda_j}{\lambda_n} - \sqrt{(1+\beta - \nicefrac{\lambda_j}{\lambda_n})^2 - 4\beta}}{2}
\end{equation*}
for all $\lambda_j$ such that $(1+\beta - \nicefrac{\lambda_j}{\lambda_n})^2 - 4\beta > 0$ and
\begin{equation*}
    \rho_{j,1} = \frac{1+\beta - \nicefrac{\lambda_j}{\lambda_n}+i\sqrt{4\beta - (1+\beta - \nicefrac{\lambda_j}{\lambda_n})^2}}{2},\quad \rho_{j,2} = \frac{1+\beta - \nicefrac{\lambda_j}{\lambda_n} - i\sqrt{4\beta - (1+\beta - \nicefrac{\lambda_j}{\lambda_n})^2}}{2}
\end{equation*}
for all $\lambda_j$ such that $(1+\beta - \nicefrac{\lambda_j}{\lambda_n})^2 - 4\beta < 0$. Taking into account the assumptions of the theorem, we derive
\begin{equation*}
    \rho_{1,1} = \frac{1+\beta - \nicefrac{\lambda_1}{\lambda_n}+\sqrt{(1+\beta - \nicefrac{\lambda_1}{\lambda_n})^2 - 4\beta}}{2},\quad \rho_{1,2} = \frac{1+\beta - \nicefrac{\lambda_1}{\lambda_n} - \sqrt{(1+\beta - \nicefrac{\lambda_1}{\lambda_n})^2 - 4\beta}}{2}
\end{equation*}
and
\begin{equation*}
    \rho_{j,1} = \frac{1+\beta - \nicefrac{\lambda_j}{\lambda_n}+i\sqrt{4\beta - (1+\beta - \nicefrac{\lambda_j}{\lambda_n})^2}}{2},\quad \rho_{j,2} = \frac{1+\beta - \nicefrac{\lambda_j}{\lambda_n} - i\sqrt{4\beta - (1+\beta - \nicefrac{\lambda_j}{\lambda_n})^2}}{2}
\end{equation*}
for all $j = 2,\ldots, n$. Moreover, $|\rho_{j,1}| = |\rho_{j,2}| = \sqrt{\beta}$.

Next, using Lemma~\ref{lem:aux_lemma_1} we get
\begin{eqnarray}
    \|\mC_j \hat\mT_j^k\|_2 &=& \left\|\frac{1}{\rho_{j,2} - \rho_{j,1}} \begin{bmatrix}
         0 & 1
    \end{bmatrix}\begin{bmatrix}
             \rho_{j,2}^{k+1} - \rho_{j,1}^{k+1} & \rho_{j,1}\rho_{j,2}(\rho_{j,1}^k - \rho_{j,2}^k)\\
             \rho_{j,2}^k - \rho_{j,1}^k & \rho_{j,1}\rho_{j,2}(\rho_{j,1}^{k-1} - \rho_{j,2}^{k-1})
        \end{bmatrix}\right\|_2\notag\\
    &=& \sqrt{\left|\frac{\rho_{j,2}^k - \rho_{j,1}^k}{\rho_{j,2} - \rho_{j,1}}\right|^2 + \left|\frac{\rho_{j,1}\rho_{j,2}\left(\rho_{j,2}^{k-1} - \rho_{j,1}^{k-1}\right)}{\rho_{j,2} - \rho_{j,1}}\right|^2}\notag\\
    &\le& \sqrt{\left(\sum\limits_{t=0}^{k-1}|\rho_{j,1}|^{k-1-t}|\rho_{j,2}|^t\right)^2 + \left(|\rho_{j,1}||\rho_{j,2}|\sum\limits_{t=0}^{k-2}|\rho_{j,1}|^{k-2-t}|\rho_{j,2}|^t\right)^2}. \label{eq:max_dev_HB_technical_1}
\end{eqnarray}
Consider the expression above for $j = 1$. To bound the sums appearing in the right-hand side of the previous inequality we derive:
\begin{eqnarray*}
    \frac{|\rho_{1,2}|}{|\rho_{1,1}|} &=&  \frac{1+\beta - \nicefrac{\lambda_1}{\lambda_n}-\sqrt{(1+\beta - \nicefrac{\lambda_1}{\lambda_n})^2 - 4\beta}}{1+\beta - \nicefrac{\lambda_1}{\lambda_n}+\sqrt{(1+\beta - \nicefrac{\lambda_1}{\lambda_n})^2 - 4\beta}}\\
    &=& 1 - \frac{2\sqrt{(1+\beta - \nicefrac{\lambda_1}{\lambda_n})^2 - 4\beta}}{1+\beta - \nicefrac{\lambda_1}{\lambda_n}+\sqrt{(1+\beta - \nicefrac{\lambda_1}{\lambda_n})^2 - 4\beta}}\\
    &\le& 1 - \frac{2\sqrt{\left(1 + \left(1 - F\sqrt{\nicefrac{\lambda_1}{\lambda_n}}\right)^2 - \nicefrac{\lambda_1}{\lambda_n}\right)^2 - 4\left(1 - F\sqrt{\nicefrac{\lambda_1}{\lambda_n}}\right)^2}}{2 - \nicefrac{\lambda_1}{\lambda_n} + \sqrt{(1-\nicefrac{\lambda_1}{\lambda_n})^2}}\\
    &=& 1 - \frac{2\sqrt{\left(2 - 2F\sqrt{\nicefrac{\lambda_1}{\lambda_n}} + (F^2-1)\nicefrac{\lambda_1}{\lambda_n}\right)^2 - 4\left(1 - F\sqrt{\nicefrac{\lambda_1}{\lambda_n}}\right)^2}}{3 - 2\nicefrac{\lambda_1}{\lambda_n}}\\
    &\le& 1 - \frac{2\sqrt{(F^2 - 1)\nicefrac{\lambda_1}{\lambda_n}\left(4 - 4F\sqrt{\nicefrac{\lambda_1}{\lambda_n}} + (F^2-1)\nicefrac{\lambda_1}{\lambda_n}\right)}}{3}\\
    &=& 1 - \frac{2\sqrt{(F^2-1)\left(\left(2 - F\sqrt{\nicefrac{\lambda_1}{\lambda_n}}\right)^2 - \nicefrac{\lambda_1}{\lambda_n}\right)}}{3\sqrt{\varkappa}}\le  1 - \frac{\sqrt{F^2-1}}{\sqrt{3\varkappa}},
\end{eqnarray*}
where the first inequality follows from the fact the function $g(\beta) = (1+\beta - \nicefrac{\lambda_1}{\lambda_n})^2 - 4\beta$ is decreasing for $\beta \le (1 - \sqrt{\nicefrac{\lambda_1}{\lambda_n}})^2$, and in the last inequality we apply $1 - F\sqrt{\nicefrac{\lambda_1}{\lambda_n}} \ge 0$, $\nicefrac{\lambda_1}{\lambda_n} \le \nicefrac{1}{10000} < \nicefrac{1}{4}$, and $\varkappa = \nicefrac{\lambda_n}{\lambda_1}$. Therefore,
\begin{eqnarray*}
    \sum\limits_{t=0}^{k-1}|\rho_{1,1}|^{k-1-t}|\rho_{1,2}|^t &=& |\rho_{1,1}|^{k-1} \sum\limits_{t=0}^{k-1} \left(\frac{|\rho_{1,2}|}{|\rho_{1,1}|}\right)^t \le \sum\limits_{t=0}^{\infty}\left(1 - \frac{\sqrt{F^2-1}}{\sqrt{3\varkappa}}\right)^t =  \frac{\sqrt{3\varkappa}}{\sqrt{F^2-1}}
\end{eqnarray*}
and, similarly,
\begin{eqnarray*}
    |\rho_{j,1}||\rho_{j,2}|\sum\limits_{t=0}^{k-2}|\rho_{j,1}|^{k-2-t}|\rho_{j,2}|^t &\le& \sum\limits_{t=0}^{k-2} \left(\frac{|\rho_{1,2}|}{|\rho_{1,1}|}\right)^t \le \sum\limits_{t=0}^{\infty}\left(1 - \frac{\sqrt{F^2-1}}{\sqrt{3\varkappa}}\right)^t =  \frac{\sqrt{3\varkappa}}{\sqrt{F^2-1}}.
\end{eqnarray*}
Plugging these upper bounds in \eqref{eq:max_dev_HB_technical_1} we derive
\begin{equation}
    \|\mC_j \hat\mT_j^k\|_2 \le \frac{\sqrt{6\varkappa}}{\sqrt{F^2-1}}. \label{eq:max_dev_HB_technical_2}
\end{equation}

Next, we consider the right-hand side of \eqref{eq:max_dev_HB_technical_1} for $j=2,\ldots,n$. In this case, $|\rho_{j,1}| = |\rho_{j,2}| = \sqrt{\beta} \le 1 - \nicefrac{F}{\sqrt{\varkappa}}$. Therefore,
\begin{eqnarray*}
    \sum\limits_{t=0}^{k-1}|\rho_{j,1}|^{k-1-t}|\rho_{j,2}|^t &=& k\left(\sqrt{\beta}\right)^{k-1} \le k\left(1 - \frac{F}{\sqrt{\varkappa}}\right)^{k-1} \le (k-1)\exp\left(-(k-1)\frac{F}{\sqrt{\varkappa}}\right) + 1
\end{eqnarray*}
and, similarly,
\begin{eqnarray*}
    |\rho_{j,1}||\rho_{j,2}|\sum\limits_{t=0}^{k-2}|\rho_{j,1}|^{k-2-t}|\rho_{j,2}|^t &=& (k-1)\left(\sqrt{\beta}\right)^{k} \le (k-1)\exp\left(-(k-1)\frac{F}{\sqrt{\varkappa}}\right).
\end{eqnarray*}
Since the maximal value of the function $g(x) = x a^x$ for $x\ge 0$ equals $-\nicefrac{1}{(e\ln(a))}$, we have
\begin{equation*}
    (k-1)\exp\left(-(k-1)\frac{F}{\sqrt{\varkappa}}\right) \le - \frac{1}{e\ln\left(\exp\left(-\frac{F}{\sqrt{\varkappa}}\right)\right)} = \frac{\sqrt{\varkappa}}{eF}.
\end{equation*}
Putting all together we obtain for all $j=2,\ldots, n$
\begin{equation}
    \|\mC_j \hat\mT_j^k\|_2 \overset{\eqref{eq:max_dev_HB_technical_1}}{\le} \sqrt{\left(\frac{\sqrt{\varkappa}}{eF} + 1\right)^2 + \frac{\varkappa}{e^2F^2}} \le \frac{\sqrt{5\varkappa}}{eF}, \label{eq:max_dev_HB_technical_3}
\end{equation}
where we use $F \le \sqrt{\varkappa}$.

Finally, with \eqref{eq:max_dev_HB_technical_2} and \eqref{eq:max_dev_HB_technical_3} in hand we derive
\begin{equation*}
    \text{dev}_{\algname{AHB}}(\alpha,\beta) \le \text{dev}_{\algname{HB}}(\alpha,\beta) = \|\mC \mT^k\|_2 = \max\limits_{j=1,\ldots,n} \|\mC_j\hat \mT_j^k\|_2 \le \frac{\sqrt{6\varkappa}}{\sqrt{F^2-1}}.
\end{equation*}
Theorem 1 from \cite{danilova2018non} implies that
\begin{equation*}
    \text{dev}_{\algname{HB}}(\alpha^*,\beta^*) \ge  \frac{\sqrt{\varkappa}}{2e},
\end{equation*}
where $\alpha^*$ and $\beta^*$ are given in \eqref{eq:optimal_params}. Therefore,
\begin{equation*}
    \text{dev}_{\algname{AHB}}(\alpha,\beta) \le \text{dev}_{\algname{HB}}(\alpha,\beta) \le \frac{2e\sqrt{6}}{\sqrt{F^2-1}}\text{dev}_{\algname{HB}}(\alpha^*,\beta^*).
\end{equation*}

\section{Missing Proofs from Section~\ref{sec:non_quadratics}}

\subsection{Proof of Lemma~\ref{lem:one_iter_progress_HB}}

Using recursion \eqref{eq:virtual_iter_recurrence} for the virtual iterates defined in \eqref{eq:virtual_iterates_HB}, we derive 
    \begin{eqnarray}
        \|\widetilde{x}_{k+1} - x_*\|^2 &=& \|\widetilde{x}_k - x_*\|^2 - \frac{2\alpha}{1-\beta}\langle \widetilde{x}_k - x_*, \nabla f(x_k) \rangle + \frac{\alpha^2}{(1-\beta)^2}\|\nabla f(x_k)\|^2\notag\\
        &=& \|\widetilde{x}_k - x_*\|^2 - \frac{2\alpha}{1-\beta}\langle x_k - x_*, \nabla f(x_k) \rangle + \frac{2\alpha}{1-\beta}\langle x_k - \widetilde{x}_k, \nabla f(x_k) \rangle \notag\\
        &&\quad + \frac{\alpha^2}{(1-\beta)^2}\|\nabla f(x_k)\|^2 \label{eq:one_iter_HB_technical_1}.
    \end{eqnarray}
    From $\mu$-strong convexity and $L$-smoothness of $f$ we have (e.g., see \cite{nesterov2018lectures})
    \begin{eqnarray}
        \langle x_k - x_*, \nabla f(x_k) \rangle &\ge& f(x_k) - f(x_*) + \frac{\mu}{2}\|x_k - x_*\|^2 \notag\\
        \|\nabla f(x_k)\|^2 &\le& 2L\left(f(x_k) - f(x_*)\right). \label{eq:L_smoothness_cor}
    \end{eqnarray}
    Together with \eqref{eq:one_iter_HB_technical_1} these relations give
    \begin{eqnarray*}
        \|\widetilde{x}_{k+1} - x_*\|^2 &\le& \|\widetilde{x}_k - x_*\|^2 - \frac{\alpha\mu}{1-\beta}\|x_k - x_*\|^2 - \frac{2\alpha}{1-\beta}\left(1 - \frac{\alpha L}{1-\beta}\right)\left(f(x_k) - f(x_*)\right)\\
        &&\quad + \frac{2\alpha}{1-\beta}\langle x_k - \widetilde{x}_k, \nabla f(x_k) \rangle.
    \end{eqnarray*}
    Next, we estimate the second and the fourth terms in the inequality above. Since $\|a+b\|^2 \ge \frac{1}{2}\|a\|^2 - \|b\|^2$ for all $a,b\in \R^n$ (see also \eqref{eq:squared_norm_sum}), we can estimate the second term as
    \begin{equation*}
        - \frac{\alpha\mu}{1-\beta}\|x_k - x_*\|^2 \le -\frac{\alpha\mu}{2(1-\beta)}\|\widetilde{x}_k - x_*\|^2 + \frac{\alpha\mu}{1-\beta}\|x_k - \widetilde{x}_k\|^2.
    \end{equation*}
    Using Fenchel-Young inequality \eqref{eq:fenchel_young_inequality}, we derive
    \begin{eqnarray*}
        \frac{2\alpha}{1-\beta}\langle x_k - \widetilde{x}_k, \nabla f(x_k) \rangle &\le& \frac{2\alpha L}{1-\beta}\|x_k - \widetilde{x}_k\|^2 + \frac{2\alpha}{4L(1-\beta)}\|\nabla f(x_k)\|^2\\
        &\overset{\eqref{eq:L_smoothness_cor}}{\le}& \frac{2\alpha L}{1-\beta}\|x_k - \widetilde{x}_k\|^2 + \frac{\alpha}{1-\beta}\left(f(x_k) - f(x_*)\right).
    \end{eqnarray*}
    Putting all togetherm, we obtain
    \begin{eqnarray*}
        \|\widetilde{x}_{k+1} - x_*\|^2 &\le& \left(1 - \frac{\alpha\mu}{2(1-\beta)}\right)\|\widetilde{x}_k - x_*\|^2 - \frac{2\alpha}{1-\beta}\left(\frac{1}{2} - \frac{\alpha L}{1-\beta}\right)\left(f(x_k) - f(x_*)\right)\\
        &&\quad + \frac{\alpha}{1-\beta}\left(2L + \mu\right)\|x_k - \widetilde{x}_k\|^2\\
        &\overset{\eqref{eq:virtual_iterates_HB},\eqref{eq:params_WAHB_1}}{\le}& \left(1 - \frac{\alpha\mu}{2(1-\beta)}\right)\|\widetilde{x}_k - x_*\|^2 - \frac{\alpha}{1-\beta}\left(f(x_k) - f(x_*)\right) + \frac{3L\alpha\beta^2}{(1-\beta)^3}\|m_{k-1}\|^2
    \end{eqnarray*}
    that finishes the proof.

\subsection{Proof of Lemma~\ref{lem:weighted_sum_of_momentums}}

Using the update rule for $m_k$, we get
    \begin{eqnarray*}
        \|m_k\|^2 &=& \|\beta m_{k-1} + \alpha\nabla f(x_k)\|^2\\
        &\overset{\eqref{eq:1+lambda}}{\le}& \beta^2\left(1 + \frac{1-\beta}{\beta}\right)\|m_{k-1}\|^2 + \alpha^2\left(1 + \frac{\beta}{1-\beta}\right)\|\nabla f(x_k)\|^2\\
        &\overset{\eqref{eq:L_smoothness_cor}}{\le}& \beta\|m_{k-1}\|^2 + \frac{2L\alpha^2}{1-\beta}\left(f(x_k) - f(x_*)\right)
    \end{eqnarray*}
    implying
    \begin{equation*}
        \|m_{k-1}\|^2 \le \frac{2L\alpha^2}{1-\beta}\sum\limits_{l=0}^{k-1}\beta^{k-1-l}\left(f(x_{l}) - f(x_*)\right).
    \end{equation*}
    Summing up these inequalities for $k=0,1,\ldots, K$ with weights $w_k = \left(1 - \frac{\alpha\mu}{2(1-\beta)}\right)^{-(k+1)}$, we derive
    \begin{eqnarray}
        \frac{3L\alpha\beta^2}{(1-\beta)^3}\sum\limits_{k=0}^{K}w_k\|m_{k-1}\|^2 &\le& \frac{6L^2\alpha^3\beta^2}{(1-\beta)^4}\sum\limits_{k=0}^{K}\sum\limits_{l=0}^{k-1}w_k\left(f(x_l) - f(x_*)\right)\beta^{k-1-l}\notag\\
        &\le& \frac{6L^2\alpha^3\beta}{(1-\beta)^4}\sum\limits_{k=0}^{K}\sum\limits_{l=0}^{k}w_k\left(f(x_l) - f(x_*)\right)\beta^{k-l}. \label{eq:weighted_sum_momentums_technical_1}
    \end{eqnarray}
    Next, we upper bound $w_k$ in the following way: for all $l = 0,1,\ldots,k$
    \begin{equation}
        w_k = \left(1 - \frac{\alpha\mu}{2(1-\beta)}\right)^{-(k-l)}w_{l} \overset{\eqref{eq:1-q/2}}{\le} \left(1 + \frac{\alpha\mu}{1-\beta}\right)^{k-l}w_l \overset{\eqref{eq:params_WAHB_2}}{\le} \left(1 +\frac{1-\beta}{2}\right)^{k-l}w_l. \notag
    \end{equation}
    Plugging this inequality into \eqref{eq:weighted_sum_momentums_technical_1} we get
    \begin{eqnarray}
        \frac{3L\alpha\beta^2}{(1-\beta)^3}\sum\limits_{k=0}^{K}w_k\|m_{k-1}\|^2 &\le& \frac{6L^2\alpha^3\beta}{(1-\beta)^4}\sum\limits_{k=0}^{K}\sum\limits_{l=0}^{k}w_l\left(f(x_l) - f(x_*)\right)\left(1 +\frac{1-\beta}{2}\right)^{k-l}\beta^{k-l} \notag\\
        &\overset{\eqref{eq:1-q/2}}{\le}& \frac{6L^2\alpha^3\beta}{(1-\beta)^4}\sum\limits_{k=0}^{K}\sum\limits_{l=0}^{k}w_l\left(f(x_l) - f(x_*)\right)\left(1 - \frac{1-\beta}{2}\right)^{k-l}\notag\\
        &\le& \frac{6L^2\alpha^3\beta}{(1-\beta)^4}\left(\sum\limits_{k=0}^K w_k \left(f(x_k) - f(x_*)\right)\right)\left(\sum\limits_{k=0}^\infty \left(1 - \frac{1-\beta}{2}\right)^k\right) \notag\\
        &=& \frac{12L^2\alpha^3\beta}{(1-\beta)^5}\sum\limits_{k=0}^K w_k \left(f(x_k) - f(x_*)\right).\notag
    \end{eqnarray}
    Note that our choice of $\alpha$ \eqref{eq:params_WAHB_2} implies
    \begin{equation*}
        \frac{12L^2\alpha^3\beta}{(1-\beta)^5} \le \frac{\alpha}{4(1-\beta)}.
    \end{equation*}
    Together with previous inequality it gives \eqref{eq:weighted_sum_of_momentums}.

\subsection{Proof of Theorem~\ref{thm:WAHB_main_result}}

From Lemma~\ref{lem:one_iter_progress_HB} we have
    \begin{equation*}
        \frac{\alpha}{2(1-\beta)}\left(f(x_k) - f(x_*)\right) \le \left(1 - \frac{\alpha\mu}{2(1-\beta)}\right)\|\widetilde{x}_k - x_*\|_2^2 - \|\widetilde{x}_{k+1} - x_*\|_2^2 + \frac{3L\alpha\beta^2}{(1-\beta)^3}\|m_{k-1}\|_2^2.
    \end{equation*}
    Summing up these inequalities for $k = 0,1,\ldots,K$ with weights $w_k = \left(1 - \frac{\alpha\mu}{2(1-\beta)}\right)^{-(k+1)}$, we get
    \begin{eqnarray*}
        \frac{\alpha}{2(1-\beta)}\sum\limits_{k=0}^K w_k\left(f(x_k) - f(x_*)\right) &\le& \sum\limits_{k=0}^K\left(w_k\left(1 - \frac{\alpha\mu}{2(1-\beta)}\right)\|\widetilde{x}_k - x_*\|_2^2 - w_k\|\widetilde{x}_{k+1} - x_*\|_2^2\right)\\
        &&\quad + \frac{3L\alpha\beta^2}{(1-\beta)^3}\sum\limits_{k=0}^K w_k\|m_{k-1}\|_2^2\\
        &\overset{\eqref{eq:weighted_sum_of_momentums}}{\le}& \sum\limits_{k=0}^K\left(w_{k-1}\|\widetilde{x}_k - x_*\|_2^2 - w_k\|\widetilde{x}_{k+1} - x_*\|_2^2\right)\\
        &&\quad + \frac{\alpha}{4(1-\beta)}\sum\limits_{k=0}^K w_k\left(f(x_k) - f(x_*)\right)\\
        &=& \|x_0 - x_*\|_2^2 + \frac{\alpha}{4(1-\beta)}\sum\limits_{k=0}^K w_k\left(f(x_k) - f(x_*)\right).
    \end{eqnarray*}
    Rearranging the terms and dividing both sides of the inequality by $W_K = \sum_{k=0}^K w_k$, we derive
    \begin{equation*}
        \frac{1}{W_K}\sum\limits_{k=0}^K w_k\left(f(x_k) - f(x_*)\right) \le \frac{4(1-\beta)\|x_0 - x_*\|_2^2}{\alpha W_K}.
    \end{equation*}
    Using Jensen's inequality, we obtain
    \begin{equation*}
        f(\overline{x}_K) \le \frac{1}{W_K}\sum\limits_{k=0}^K w_k f(x_k)
    \end{equation*}
    that implies \eqref{eq:WAHB_main_result}. Next, when $\mu > 0$ we have $W_K \ge w_{K-1} = \left(1 - \frac{\alpha\mu}{2(1-\beta)}\right)^{-K}$ that gives \eqref{eq:WAHB_main_result_str_cvx}. Finally, when $\mu = 0$ we have $W_K = {K+1} > K$ implying \eqref{eq:WAHB_main_result_cvx}.

\subsection{Proof of Theorem~\ref{thm:R-AHB}}

Theorem~\ref{thm:WAHB_main_result} for $\mu = 0$ implies that for $t = 1,2,\ldots,\tau$
    \begin{equation}
        f(\widehat x_t) - f(x_*) \le \frac{4(1-\beta)\widehat{R}_{t-1}^2}{\alpha N}, \label{eq:R-AHB_technical_1}
    \end{equation}
    where $\widehat{R}_{t} = \|\widehat{x}_t - x_*\|_2$ for $t = 0,1,\ldots,\tau$. In the remaining part of the prove, we derive via induction that for $t = 1,2,\ldots,\tau$
    \begin{equation}
        f(\widehat x_t) - f(x_*) \le \frac{\mu R_0^2}{2^{t+1}},\quad \widehat{R}_t \le \frac{R_0^2}{2^t}, \label{eq:R-AHB_technical_2}
    \end{equation}
    where $R_0 \ge \|x_0 - x_*\|_2 = \|\widehat{x}_0 - x_*\|_2$. First of all, for $t = 1$ we have
    \begin{equation*}
        f(\widehat{x}_1) - f(x_*) \overset{\eqref{eq:params_R-AHB},\eqref{eq:R-AHB_technical_1}}{\le} \frac{\mu R_0^2}{4}.
    \end{equation*}
    From $\mu$-strong convexity of $f$ we derive
    \begin{equation*}
        \frac{\mu \widehat{R}_1^2}{2} \le f(\widehat{x}_1) - f(x_*) \quad \Longrightarrow \quad \widehat{R}_1^2 \le \frac{R_0^2}{2}.
    \end{equation*}
    Next, assume that \eqref{eq:R-AHB_technical_2} holds for all $t = 1,2,\ldots,k < \tau$ and let us prove it for $t = k+1$. From \eqref{eq:R-AHB_technical_1} we have
    \begin{equation*}
        f(\widehat x_{k+1}) - f(x_*) \le \frac{4(1-\beta)\widehat{R}_{k}^2}{\alpha N} \overset{\eqref{eq:R-AHB_technical_1}}{\le} \frac{(1-\beta) R_0^2}{2^{k-2}\alpha N} \overset{\eqref{eq:params_R-AHB}}{\le} \frac{\mu R_0^2}{2^{k+2}}.
    \end{equation*}
    Again, applying $\mu$-strong convexity of $f$ we derive
    \begin{equation*}
        \frac{\mu \widehat{R}_{k+1}^2}{2} \le f(\widehat{x}_{k+1}) - f(x_*) \quad \Longrightarrow \quad \widehat{R}_{k+1}^2 \le \frac{R_0^2}{2^{k+1}}
    \end{equation*}
    that finishes the proof of \eqref{eq:R-AHB_technical_2}. Therefore, after $\tau = \max\{\left\lceil\log_2(\nicefrac{\mu R_0^2}{\varepsilon})\right\rceil - 1,1\}$ iterations \algname{R-AHB} finds such point $\widehat{x}_\tau$ that
    \begin{equation*}
        f(\widehat{x}_\tau) - f(x_*) \le \frac{\mu R_0^2}{2^{\tau+1}} \le \frac{\mu R_0^2}{2^{\log_2(\nicefrac{\mu R_0^2}{\varepsilon})}} = \varepsilon.
    \end{equation*}
    Finally, if
    \begin{equation*}
        \alpha = \min\left\{\frac{1-\beta}{4L}, \frac{(1-\beta)^2}{4L\sqrt{3\beta}} \right\},
    \end{equation*}
    then the total number of \algname{AHB} iterations equals
    \begin{equation*}
        N\tau = \cO\left(\left(\frac{L}{\mu} + \frac{L\sqrt{\beta}}{\mu(1-\beta)}\right)\log\frac{\mu R_0^2}{\varepsilon}\right).
    \end{equation*}

\end{document}